\documentclass[11pt,a4paper,british]{article}
\usepackage{authblk}
\def\isfinal{1}
\def\ifnofinal{\ifnum\isfinal=0}

\usepackage{etoolbox}
\patchcmd{\thebibliography}{*}{}{}{}
\pretocmd\thebibliography{\csname c@secnumdepth\endcsname=-2 }{}{}

\usepackage[T1]{fontenc}
\usepackage{csquotes}
\usepackage[french,main=british]{babel}
\usepackage{bold-extra}

\usepackage{relsize}

\usepackage[backend=biber,doi=false,url=false,isbn=false]{biblatex}
\DeclareFieldFormat[article]{volume}{\textbf{#1}}
\DeclareFieldFormat[article]{number}{-#1}
\DeclareFieldFormat[book]{volume}{\textbf{#1}}

\renewbibmacro{in:}{%
  \ifentrytype{article}{}{\printtext{\bibstring{in}\intitlepunct}}}

\renewbibmacro*{volume+number+eid}{%
  \printfield{volume}
  \printfield[parens]{number}%
  \setunit{\addcomma\space}%
  \printfield{eid}}
\patchcmd{\thebibliography}{*}{}{}{}
\pretocmd\thebibliography{\csname c@secnumdepth\endcsname=-2 }{}{}
\DeclareFieldFormat{postnote}{#1}
\DeclareFieldFormat{multipostnote}{#1}
\usepackage[titletoc]{appendix}

\usepackage[margin=6em]{geometry}

\usepackage{amssymb,amsmath,amsthm,aliascnt,pdfsync,stmaryrd,mathtools,tabularx}
\usepackage{nicematrix}
\usepackage{chemarrow,tikz-cd}
\usepackage{enumitem} 
\setlist[enumerate]{label=(\roman*)} 

\usepackage[nonewpage]{indextools}              
\indexsetup{level=\subsection*}
\makeindex[title=Index of notation]
\makeindex[name=ch1,title=Part I,columns=4]
\makeindex[name=ch2,title=Part II,columns=4]

\usepackage{denpastyle}

\ifnofinal
	\usepackage[inline]{showlabels}
\fi

\addto\extrasfrench{

}
\addto\extrasbritish{

}

\newcommand{\ubar}[1]{\underline{\smash{#1}}}

\newcommand{\idv}{\mathrm{red}}
\newcommand{\intt}{{\smallint}}

\newcommand{\omegaa}{{\mathsmaller{\Omega}}}
\newcommand{\dW}{\partial^{\mathsmaller{\mathsmaller{W}}}}

\renewcommand{\to}{\longrightarrow}

\title{Knizhnik--Zamolodchikov functor for degenerate double affine Hecke algebras : algebraic theory}
\date{\today}
\author{\sc Wille Liu\thanks{ \textsc{Wille Liu},\;
		Max-Planck-Institut f\"ur Mathematik,
	Vivatsgasse 7, 53111 Bonn, Germany\par\nopagebreak
  \textit{E-mail}: \href{mailto:wille@mpim-bonn.mpg.de}{wille@mpim-bonn.mpg.de}}
}

\begin{document}
\maketitle
\abstract{In this article, we define an algebraic version of the Knizhnik--Zamolodchikov functor for the degenerate double affine Hecke algebras (a.k.a. trigonometric Cherednik algebras). We compare it with the KZ monodromy functor constructed by Varagnolo--Vasserot. We prove the double centraliser property for our functor and give a characterisation of its kernel. We establish these results for a family of algebras, called quiver double Hecke algebras, which includes the degenerate double affine Hecke algebras as special cases.}

\section*{Introduction}
\subsection*{Degenerate double affine Hecke algebras}
The degenerate double affine Hecke algebras (dDAHA), also known as trigonometric Cherednik algebras, were introduced by I. Cherednik in his study of integration of the trigonometric form of the Knizhnik--Zamolodchikov equations (KZ)~\cite{cherednik92}.  \par
The degenerate double affine Hecke algebras, unlike their non-degenerate version and its rational degeneration, are not ``symmetric'': it contains a polynomial subalgebra and a Laurent polynomial subalgebra. Due to this asymmetry, one can adopt two different points of view to study the dDAHA: either viewing it 
\begin{enumerate}
	\item
		as the algebra generated by regular functions on a torus $T^{\vee}$ attached to a root system $R$, the Weyl group of $R$ acting the torus $T^{\vee}$ and the trigonometric Dunkl operators on it, or 
	\item
		as the algebra generated by Demazure-like difference operators on $E$, where $E$ is an affine space which carries an affine root system; this is the affine version of the graded affine Hecke algebras of G. Lusztig~\cite{lusztig89}.
\end{enumerate}
The former approach allows one to apply various techniques of $\calD$-modules, symplectic geometry and is closer to the theory of rational Cherednik algebras~\cite{EM10,BG15}; the latter approach allows one to apply cohomological, K-theoretic or sheaf-theoretic methods~\cite{CG,vasserot05}, and is closer to the (non-degenerate) double affine Hecke algebras. \par
In the present work, we will adopt the second approach most of the time. We show that with this point of view, the dDAHAs can be easily generalised and are quite flexible in the choice of parameters. We show also that some of the features from first approach can be recovered with the second approach, namely the integration of the KZ equations. \par
\vspace{10pt}
\subsection*{Quiver Hecke algebras}
The quiver Hecke algebras, also known as Khovanov--Lauda--Rouquier algebras, were introduced in~\cite{KhLa09} and~\cite{rouquier08b}. They were introduced in the purpose of categorifying the Drinfel'd--Jimbo quantum groups for Kac--Moody algebras as well as their integrable representations. \par

It was proven by Brundan--Kleshchev--McNamara~\cite{BKN14} and Kato~\cite{kato14} that quiver Hecke algebras for Dynkin quivers of finite ADE types have pretty nice homological properties. Retrospectively speaking, they proved that the categories of graded modules over these algebras carry an {\it affine highest weight structure} in the sense of~\cite{kleshchev15}. As a consequence, these algebras have finite global dimension. However, once one goes beyond the family of finite type, the quiver Hecke algebras often have infinite global dimension. The simplest example would be the cyclic quivers of length $\ge 2$. According to the result of Brundan--Kleshchev~\cite{BK09} and Rouquier~\cite{rouquier08b}, the quiver Hecke algebras of cyclic quivers are equivalent to affine Hecke algebras for $GL_n$ with parameter at roots of unity. The representation theory of affine Hecke algebras at roots of unity is known to share several features of the modular representation theory finite groups. Notably, there are fewer simple modules in the modular case than there are in the ordinary case. \par

One approach to the modular representation theory is to resolve this lack of simple objects by finding a larger, but better behaved category, of which the modular category is a quotient. In the case of modular representation theory of symmetric groups, one uses the Schur algebras as resolution via the Schur--Weyl duality. In the same spirit, for Hecke algebras of complex reflection groups, the rational Cherednik algebras provide resolution, as it was first established in~\cite{GGOR03}. For affine Hecke algebras, the resolution would be the degenerate double affine Hecke algebras. This perspective appeared in~\cite{VV04}, where degenerate DAHAs are viewed as replacement for affine q-Schur algebras in relation with affine Hecke algebras~\cfauto{subsec:dimglH} and \autoref{subsec:monodromie}. We will introduce a new family of algebras, called \emph{quiver double Hecke algebras}, which we believe to play the r\^ole of ``resolution'' for quiver Hecke algebras. \par

\subsection*{Results of the present article}
Let $(V, R)$ be an irreducible finite root system and let $(E, S)$ be its affinisation (the definition is recalled in~\autoref{subsec:sysaff}). In particular, $E \cong V$ is a euclidean affine space. We fix a basis $\Delta_0\subset R$, which extends in a standard way to an affine basis $\Delta\subset S$. The affine Weyl group $W_S$ is generated by affine simple reflections $s_a$ for $a\in \Delta$ and the finite Weyl group $W_R\subset W_S$ is the subgroup generated by $s_a$ for $a\in \Delta_0$. The extended affine Weyl group $\til W_S$ acts on $S$.  \par

The degenerate double affine Hecke algebra attached to $(E, S)$ is given by $\bbH = \bfC W_S\otimes \bfC[E]$ as vector space. The multiplication of $\bbH$ depends on a function $h: S  \to \bfC$, called parameters, see~\autoref{subsec:dDAHA} for the precise definition.  For $\lambda\in E$, let $\calO_{\lambda}(\bbH)$ denote the category of finitely generated $\bbH$-modules on which the subalgebra $\bfC[E]$ acts locally finitely with eigenvalues lying in the orbit $W_S\cdot \lambda\subset E$.

The affine Hecke algebra attached to $(V, R)$ given by $\bbK = H_R\otimes \bfC[T]$, where $H_R$ is the Iwahori-Hecke algebra of type $(W_R, \Delta_0)$ and $\bfC[T]$ is the group algebra of the weight lattice of the root system $(V, R)$. See~\autoref{subsec:AHA} for the precise definition. For $\ell\in V$, let $\calO_{\ell}(\bbK)$ denote the category of finite-dimensional $\bbK$-modules on which the subalgebra $\bfC[T]$ acts with eigenvalues lying in the orbit $W_R\cdot \ell\subset T$. \par

There is an exponential map $\exp: E \to T$. Fix $\lambda_0\in E$ and let $\ell_0 = \exp(\lambda_0)\in T$. Denote by $\bbV: \calO_{\lambda_0}(\bbH)\to \calO_{\ell_0}(\bbK)$ the monodromy functor for the Knizhnik--Zamolodchikov equations introduced by Varagnolo--Vasserot in~\cite{VV04}. We show in~\autoref{prop:quotM} that $M$ is a quotient functor. The first main result is the following:

\begin{theointro}[=\autoref{defi:KZalg}+\autoref{prop:compker}]\label{theo:A}
	There is a quotient functor $\bfV: \calO_{\lambda_0}(\bbH)\to \calO_{\ell_0}(\bbK)$ defined in algebraic terms such that
	\begin{align*}
		\ker \bfV = \ker\bbV.
	\end{align*}
\end{theointro}
We expect that there exists an isomorphism $\bfV\cong\bbV$. In order to construct $\bfV$, we introduce in~\autoref{subsec:Hlambda} and~\autoref{subsec:Kell} two auxiliary algebras $\bfH_{\lambda_0}$ and $\bfK_{\ell_0}$ and show in~\autoref{prop:equivHH} and~\autoref{prop:equivKK} that $\bfH_{\lambda_0}$ and $\bfK_{\ell_0}$ are Morita-equivalent respectively to $\calO_{\lambda_0}(\bbH)$ and $\calO_{\ell_0}(\bbK)$. By analysing the structure of the quiver-Hecke-like algebras $\bfH_{\lambda_0}$ and $\bfK_{\ell_0}$, we show in~\autoref{theo:isomKH} that there exists an idempotent $\bfe_{\gamma}\in \bfH_{\lambda_0}$ such that the idempotent subalgebra $\bfe_{\gamma}\bfH_{\lambda_0}\bfe_{\gamma}$ is isomorphic to $\bfK_{\ell_0}$. This allows us to define the functor $\bfV$ as the idempotent truncation by $\bfe_{\gamma}$. \par

The second main result concerns $\bfV$:
\begin{theointro}[=\autoref{theo:bicommutante}+\autoref{prop:catV}]\label{theo:B}
	The following statements hold:
	\begin{enumerate}
		\item
			The functor $\bfV$ satisfies the double centraliser property (i.e. fully faithful on projective objects) after passing to a suitable completion of $\calO_{\lambda_0}(\bbH)$ and $\calO_{\ell_0}(\bbK)$.
		\item
			The kernel $\ker \bfV$ is the Serre subcategory generated by simple objects $L\in \calO_{\lambda_0}(\bbH)$ such that the projective envelope of $L$ in the completion of $\calO_{\lambda_0}(\bbH)$ is not relatively injective with respect to the categorical centre $\rmZ(\calO_{\lambda_0}(\bbH))$.
	\end{enumerate}
\end{theointro}
Notice that by the comparison result~\autoref{theo:A}, the statements of~\autoref{theo:B} also hold for $\bbV$. The second statement of~\autoref{theo:B} implies in particular that the subcategory $\ker\bfV$ is an invariant of the category $\calO_{\lambda_0}(\bbH)$. In fact, we construct $\bfV$ and establish~\autoref{theo:B} for a greater family of algebras, \emph{quiver double Hecke algebras}, which are introduced in~\autoref{subsec:Aomega}. This family of algebras seems to be related to a localised Iwahori version of Coulomb branch algebras of Braverman--Finkelberg--Nakajima~\cite{BFNII} for semisimple groups. 

\subsection*{Related works}
As mentioned above, the algebra $\bfA^{\omega}$ that we introduce in~\autoref{partII} is expected to be related to Iwahori version of the quantised Coulomb branch algebras. There exist in the literature some works on the representation theory of such algebras with an approach similar to ours. \par
In~\cite{webster19representation}, B. Webster studied a module category of the rational Cherednik algebra for the complex reflection group $G(\ell, 1, n)$ whose objects admit a weight decomposition for the action of a polynomial subalgebra defined by Dunkl--Opdam~\cite{DO03}. He introduced an algebraic version of the KZ functor and he classified the simple objects of that category. The results were later generalised in~\cite{lepage20rational}, to the rational Cherednik algebra for $G(\ell, d, n)$. \par

Our construction of KZ functor $\bfV$ can be regarded as a variant of theirs. One can expect that their functor also satisfies the properties listed in~\autoref{theo:B}.

\subsection*{Organisation}
This paper is composed of two parts. The first part serves mainly as preliminary materials and motivation for the second part. The proof of most of the statements in the first part can be found in the literature~\cite{lusztig89,cherednik94,opdam00,VV04,VV09}. \par

We review briefly the affine root systems in~\autoref{subsec:sysaff}, the dDAHAs in~\autoref{subsec:dDAHA} and the affine Hecke algebras (AHA) in~\autoref{subsec:AHA}. We introduce the idempotent form of these algebras, each controlling a block of the category $\calO$ of both algebras. The definition of idempotent forms is a straightforward generalisation of the result of Brundan--Kleshchev~\cite{BK09} and Rouquier~\cite{rouquier08b} on the equivalence between affine Hecke algebras for $\GL_n$ and quiver Hecke algebras for linear and cyclic quivers. \par
We recall in~\autoref{subsec:monodromie} the monodromy functor $\bbV$ introduced in~\cite{VV04} as the trigonometric counterpart of the KZ functor of~\cite{GGOR03}. We prove that it is a quotient functor in the sense of Gabriel. \par
We discuss in~\autoref{sec:VV} the relations between the monodromy functor $\bbV$ and the functor $\bfV$, which will be defined in algebraic terms in~\autoref{subsec:V}.  \par

In the second part we introduce quiver double Hecke algebras (QDHA). They can be viewed as a generalisation of degenerate double affine Hecke algebras (dDAHA) or as an affinisation of quiver Hecke algebras (QHA). \par

In~\autoref{sec:agha}, we introduce the quiver double Hecke algebras $\bfA^{\omega}$ attached to an affine root system $(E, S)$ with spectrum being a $W_S$-orbit in $E$ and with parameter $\omega$. We define the filtration by length on $\bfA^{\omega}$ in~\autoref{subsec:filtlg} and prove the basis theorem in~\autoref{subsec:basis} with this filtration. We study the associated graded $\gr^F\bfA^{\omega}$ of the filtration by length in~\autoref{subsec:grA}. \par

In~\autoref{sec:Agmod}, we study the categories of graded and ungraded $\bfA^{\omega}$-modules. We introduce in~\autoref{subsec:indres} a functor of induction from the quiver Hecke algebras attached to the finite root system $(V, R)$ underlying $(E, S)$.\par
In~\autoref{sec:Amodfilt}, we study good filtrations on $\bfA^{\omega}$-modules and use it to define the Gelfand--Kirillov dimension of an $\bfA^{\omega}$-module. We prove that ``induced $\bfA^{\omega}$-modules'' are of maximal Gelfand--Kirillov dimension. \par
In~\autoref{sec:MQHA}, we introduce the quiver Hecke algebra $\bfB^{\omega}$ attached to a finite root system $(V, R)$ and with parameter $\omega$. We prove a basis theorem for $\bfB^{\omega}$ and we introduce a Frobenius form on $\bfB^{\omega}$. \par
In~\autoref{sec:KZ}, we prove that the algebra $\bfB^{\omega}$ is isomorphic to an idempotent subalgebra of $\bfA^{\omega}$. We use this isomorphism to define the Knizhnik--Zamolodchikov functor $\bfV$, which is a quotient functor. We give characterisations for the kernel of $\bfV$ in~\autoref{subsec:suppV} and~\autoref{subsec:catV}. The double centraliser property for $\bfV$ is proven in~\autoref{subsec:bicomm}. \par

In~\autoref{sec:pro}, we collect some basic facts about the category of pro-objects of abelian categories, which are used to construct completions of the categories $\calO_{\lambda_0}(\bbH)$ and $\calO_{\ell_0}(\bbK)$.

\subsection*{Acknowledgements}
This article is part of my doctoral thesis under the supervision of Eric Vasserot. I would like to express my gratitude toward him. I would also like to thank Ben Webster for pointing out the references~\cite{DO03,webster19koszul,webster19representation,lepage20rational} as well as helpful discussions

\part{Degenerate double affine Hecke algebras}\label{partI}

\section{Reminder on affine root systems}
We review the notion of affine root systems. The reference is~\cite{macdonald03}.
\subsection{Affine reflections on euclidean spaces}\label{subsec:sysaff}
Let $E$ \index[ch1]{E@$E$} be an affine euclidean space of dimension $n > 0$ and let $V$ \index[ch1]{V@$V$} be its vector space of translations. In particular, $V$ is equipped with a positive definite scalar product $\langle\relbar,\relbar\rangle: V\times V\to \bfR$\index[ch1]{@$\langle \relbar,\relbar\rangle$}. The dual space $V^*$ is identified with $V$ via the scalar product $\langle\relbar, \relbar\rangle$. Let $\bfR[E]^{\le 1}$\index[ch1]{RE@$\bfR[E]^{\le 1}$} be the space of affine functions on $E$. We have a map of differential $\partial: \bfR[E]^{\le 1}\to V^*$\index[ch1]{@$\partial$} whose kernel is the set of constant functions. The space $\bfR[E]^{\le 1}$ is equipped with a symmetric bilinear form $\langle f, g\rangle = \langle \partial f, \partial g\rangle$. For any non-constant function $f\in \bfR[E]^{\le 1}$, let $f^{\vee} = 2f / |f|^2$ and define the reflection with respect to the zero hyperplane of $f$:
\begin{align*}
	s_f:E\to E,\quad s_f(x) = x - f^{\vee}(x)\partial f
\end{align*}\index[ch1]{s@$s_a,s_{\alpha}$} 
and 
\begin{align*}
	s_f:\bfR[E]^{\le 1}\to \bfR[E]^{\le 1},\quad s_f(g) = g - \langle f^{\vee}, g\rangle f.
\end{align*}
It extends to an automorphism of the ring of $\bfC$-valued polynomial functions $s_f:\bfC[E]\to \bfC[E]$.

\subsection{Affine root systems}

An affine root system on $E$ is the pair $(E, S)$, where $S\subset \bfR[E]^{\le 1}$\index[ch1]{S@$S$} is a subset satisfying the following conditions:
\begin{enumerate}
	\item
		$S$ spans $\bfR[E]^{\le 1}$ and the elements of $S$ are non-constant functions on $E$;
	\item
		$s_a(b)\in S$ for all $a, b\in S$;
	\item
		$\langle a^{\vee}, b\rangle\in\bfZ$ for all $a, b\in S$;
	\item
		the group $W_S$\index[ch1]{WS@$W_S$} of auto-isometries on $E$ generated by $\{s_a\;;\; a\in S\}$ acts properly on $E$. 
\end{enumerate}
The group $W_S$ is called the affine Weyl group (or simply the Weyl group of $S$). An affine root system $(E, S)$ is called irreducible if there is no partition $S = S_1\sqcup S_2$ with $\langle \relbar, \relbar\rangle\mid_{S_1\times S_2} = 0$ and $S_1\neq \emptyset$ and $S_2\neq\emptyset$; it is called reduced if $a\in S$ implies $2a\notin S$. \par

Let $(E, S)$ be an affine root system. The set $R = \partial(S)\subset V^*$\index[ch1]{R@$R$} is a finite root system on $V$. Let $P = P_R\subset V$ denote the weight lattice, $Q_R = \bfZ R$ the root lattice, $P^{\vee} = P^\vee_R$ the coweight lattice and $Q^{\vee} = Q^\vee_R = \bfZ R^{\vee}$ the coroot lattice.\index[ch1]{P@$P, Q, P^{\vee}, Q^{\vee}$} \par

Conversely, let $(V, R)$ be an irreducible finite root system, reduced or not. Define $R_{\idv} = R \setminus 2R$\index[ch1]{R@$R_\idv$} to be the set of indivisible roots. Let $P = P_R$ be the weight lattice and $Q = Q_R$ the root lattice; we define the affinisation of $(V, R)$ to be the affine root system $(E, S)$ with $E = V$ and 
\begin{align*}
		S = \langle \alpha + n\;;\; n\in \bfZ, \alpha\in R_{\idv}\rangle \sqcup\langle \alpha + 2m + 1\;;\; m\in \bfZ, \alpha\in R\cap 2R\rangle.
\end{align*}
Given a basis $\Delta_0\subset R$, we form $\Delta = \Delta_0 \cup \left\{ a_0 \right\}$, where $a_0 = 1 - \theta$ with $\theta\in R$ being the highest root with respect to the basis $\Delta_0$.

\subsection{Affine Weyl group}\label{subsec:sysaff-fini}

Let $(E, S)$ be the affinisation of $(V, R)$, which is an irreducible reduced affine root system. A basis of $S$ is an $\bfR$-linearly independent subset $\Delta\subset S$ \index[ch1]{Delta@$\Delta,\Delta_0$}such that the following conditions are satisfied:
\begin{enumerate}
	\item
		$S\subset \bfN \Delta \cup -\bfN \Delta$;
	\item
		the set $\bigcap_{a\in \Delta}\left\{ x\in E\;;\; a(x) > 0 \right\}$ is non-empty.
\end{enumerate}
The $W_S$-action on $S$ induces a simple transitive $W_S$-action on the set of bases of $S$.  Upon fixing a basis $\Delta$ of $S$, let $S^+ = S \cap \bfN \Delta$ and $S^- = S\cap -\bfN \Delta$ \index[ch1]{S@$S^+,S^-$}denote the sets of positive and negative roots. \par

The parabolic Coxeter subgroup $W_R = \langle s_a \;;\; a\in \Delta_0\rangle$ of $W_S$ can be identified with the Weyl group of the finite root system $(V, R)$ \index[ch1]{WR@$W_R$}and there is an isomorphism
\begin{align*}
	Q^{\vee}_R  \rtimes W_R &\cong W_S  \\
	(\mu, w) &\mapsto X^\mu w,
\end{align*}
where the element $X^{\mu}$ \index[ch1]{X@$X^{\mu}$}acts on $S$ by $a\mapsto a - \langle \partial a, \mu\rangle$.  The extended affine Weyl group is defined to be $\tilde W_S = P^{\vee}\rtimes W_R$\index[ch1]{WS@$\til W_S$}.
It acts on $S$ by extending the $W_S$-action by the same formula $X^{\mu} a =  a - \langle \partial a, \mu\rangle$ for $\mu\in P^{\vee}$. \par

The length function is defined to be
\[
	\ell:\til W_S\to \bfN,\quad  \ell(w) = \#\left( S^+ \cap w^{-1} S^- \right).
\]
\index[ch1]{l@$\ell(w)$} It extends the usual length function on the Coxeter group $W_S$ with respect to the set of generators $\left\{s_a  \right\}_{a\in \Delta}$.  We will need the following formula for the length function.
\begin{prop}\label{prop:lg}
	For $\mu\in P^{\vee}$ and $w\in W_R$, we have
	\begin{align*}
			\ell(wX^{\mu}) = \sum_{\alpha\in R^+_{\idv}\cap w^{-1} R^-_{\idv}} |\langle \alpha,\mu\rangle + 1| + \sum_{\alpha\in R^+_{\idv}\cap w^{-1} R^+_{\idv}} |\langle \alpha,\mu\rangle| + \sum_{\alpha\in R^+\cap 2R }\frac{|\langle \alpha,\mu\rangle|}{2}. \\
		\ell(X^{\mu}w) = \sum_{\alpha\in R^+_{\idv}\cap w R^-_{\idv}} |\langle \alpha,\mu\rangle - 1| + \sum_{\alpha\in R^+_{\idv}\cap w R^+_{\idv}} |\langle \alpha,\mu\rangle| + \sum_{\alpha\in R^+\cap 2R }\frac{|\langle \alpha,\mu\rangle|}{2}. \\
	\end{align*}
\end{prop}
These formulae can be obtained  by counting the set $S^+\cap w^{-1}S^-$ along the fibres of the differential map $\partial: S\rightarrow R$. \qed \par

 \subsection{Alcoves}\label{subsec:alcoves}
For each affine root $a\in S$, let $H_a = \left\{ \lambda\in E\;;\; a(\lambda) = 0 \right\}$ be the vanishing locus of $a$. 
The affine hyperplanes $\{H_a\}_{a\in S}$ yield a simplicial cellular decomposition of $E$. The open cells are called \textbf{alcoves}. Thus the set of connected components
\begin{align*}
	\pi_0\left( E \setminus \bigcup_{a\in S} H_a \right)
\end{align*}
is the set of alcoves. The affine Weyl group $W_S$ acts simply transitively on it. When a basis $\Delta\subset S$ is given, the fundamental alcove is defined to be $\nu_0 = \bigcap_{a\in \Delta}\left\{ x\in E\;;\; a(x) > 0 \right\}$.

\section{Reminder on degenerate double affine Hecke algebra}\label{subsec:dDAHA}
Let $\left( E, S, \Delta \right)$ be an irreducible reduced affine root system with a basis. We define in this section the degenerate double affine Hecke algebra $\bbH$ attached to $(E, S, \Delta)$ and its idempotent form $\bfH_{\lambda_0}$, which is a block algebra for the category $\calO$ of $\bbH$.
\subsection{Degenerate double affine Hecke algebra \texorpdfstring{$\bbH$}{H}}
Let $h = \{ h_a \}_{a\in S}$\index[ch1]{h@$h_a$} be a $\til W_S$-invariant family of complex numbers. The degenerate double affine Hecke algebra with parameters $h$ attached to the affine root system $S$ is the associative unital $\bfC$-algebra on the vector space $\bbH = \bfC W_S \otimes \bfC[E]$\index[ch1]{H@$\bbH$} whose multiplication satisfies following properties:
\begin{itemize}
	\item
		Each of the subspaces $\bfC W_S$ and $\bfC[E]$ is given the usual ring structure, so that they are subalgebras of $\bbH$.
	\item
		$w\in \bfC W_S$ and $f\in \bfC[E]$ multiply by juxtaposition: $(w\otimes 1)(1\otimes f) = w\otimes f$. 
	\item
		$a\in \Delta$ and $f\in \bfC[E]$ satisfy the equation:
		\begin{align*}
			\left(s_{a} \otimes 1 \right)(1 \otimes f) - \left(1 \otimes s_a(f) \right)(s_{a} \otimes 1) = 1\otimes h_{a} \frac{f - s_a(f)}{a}.
		\end{align*}
\end{itemize}

\subsection{Global dimension of \texorpdfstring{$\bbH$}{H}}\label{subsec:dimglH}
Put a filtration $F$ on $\bbH$ as follows:
\begin{align*}
	F_{\le -1} \bbH = 0,\quad F_{\le 0} \bbH = \bfC W_S,\quad F_{\le 1} \bbH = \left(F_{\le 0}\bbH\right)\bfC[E]^{\le 1},\quad F_{\le n}\bbH = (F_{\le 1} \bbH)^n,\quad n\ge 2.
\end{align*}
Namely, $\bbH$ is filtered by its polynomial part $\bfC[E]$. The filtration $F$ is compatible with the multiplication and its associated graded ring is given by the skew tensor product $\gr^F \bbH \cong \bfC W_R \ltimes \left(\bfC Q^{\vee}\otimes\bfC[V]\right)$. Since $\mathrm{dim.gl}\;\bbH \le \mathrm{dim.gl}\;\gr^F \bbH$ (\cite[D.2.6]{hotta2007d}) and since $\mathrm{dim.gl}\;\bfC W_R \ltimes \left(\bfC Q^{\vee}\otimes\bfC[V]\right) = 2r$, where $r = \rk S = \dim E$, we have the following:
\begin{prop}
	The global dimension of $\bbH$ is at most $2r$.  \qed 
\end{prop}

\subsection{Category \texorpdfstring{$\calO$}{O}}\label{subsec:OH}
For each $\lambda\in E_{\bfC}$, let $\frakm_{\lambda}\subset \bfC[E]$ be the defining ideal of the closed point $\lambda\in E$.
Given any module $M\in \bbH\Mod$, for each $\lambda\in E_{\bfC}$ consider the generalised $\lambda$-weight space in $M$:
\begin{align*}
	M_{\lambda} = \bigcup_{N\ge 0}\left\{ a\in M\;;\;\frakm_{\lambda}^N a = 0\right\}.
\end{align*}
For any $\lambda_0\in E_{\bfC}$, we define $\calO_{\lambda_0}\left( \bbH \right)$\index[ch1]{O@$\calO_{\lambda_0}(\bbH)$} to be the full subcategory of finitely generated left $\bbH$-modules $\bbH\mof$ consisting of those $M$ such that
\begin{align*}
	M = \bigoplus_{\lambda\in W_S\, \lambda_0} M_{\lambda}.
\end{align*}
In other words, the polynomial subalgebra $\bfC[E]$ acts locally finitely on $M$ with eigenvalues in the $W_S$-orbit of $\lambda_0\in E_{\bfC}$. \par
From the triangular decomposition $\bbH = \bfC Q^{\vee}\otimes\bfC W_R\otimes\bfC[E]$, we deduce the following:
\begin{prop}\label{prop:coherenceT}
	For every $\lambda_0\in E$, every object of $\calO_{\lambda_0}(\bbH)$ is a coherent $\bfC Q^{\vee}$-module.  \qed
\end{prop}

\subsection{Block algebra \texorpdfstring{$\bbH^{\wedge}_{\lambda_0}$}{H}}\label{subsec:haH}
In order to study the category $\calO_{\lambda_0}(\bbH)$, it is often useful to consider a certain completion of the polynomial part $\bfC[E]$ at the orbit $W_S\lambda_0\subset E$. The completion of $\bbH$ that we will consider is similar to the one from~\cite{DFO94} in the context of Gelfand--Zetlin algebras. A similar construction has been employed in~\cite{VV09} for double affine Hecke algebras. \par
Fix once and for all $\lambda_0\in E_{\bfC}$\index[ch1]{lambda@$\lambda_0$}. Define for each $\lambda\in W_S \lambda_0$ a polynomial ring $\Pol_\lambda = \bfC[V]$\index[ch1]{Pol@$\Pol_{\lambda}$} and let $\Pol = \bigoplus_{\lambda\in W_S\lambda_0} \Pol_{\lambda}$. Define the completion
\begin{align*}
	\ha\Pol_\lambda = \varprojlim_{N}\Pol_{\lambda} / \frakm_{0}^N\Pol_{\lambda} = \bfC\llbracket V\rrbracket,\quad \ha\Pol = \bigoplus_{\lambda\in W_S \lambda_0} \ha\Pol_\lambda,
\end{align*}
where $\frakm_0\subset \Pol_{\lambda}$ is the defining ideal of $0\in V$. The completion $\ha\Pol_{\lambda}$ is equipped with the $\frakm_{0}$-adic topology and $\ha\Pol$ is equipped with the colimit topology. \par
For $\lambda\in W_S\lambda_0$, the translation $\lambda_*:V_{\bfC}\xrightarrow{\lambda+} E_{\bfC}$ yields an isomorphism
\begin{align*}
	\lambda^*: \bfC[E]\cong \bfC[V] = \Pol_{\lambda}.
\end{align*}
We define an action of $\bbH$ on $\Pol$:
\begin{align*}
	\psi = \left( \psi_\lambda \right)_{\lambda}:\bbH&\to \End^{\cont}(\ha\Pol),\quad \psi_{\lambda}: \bbH\to \Hom^{\cont}(\ha\Pol_{\lambda}, \ha\Pol). \\
\end{align*}
by setting, for $f\in \bfC[E]$ and $a\in \Delta$, 
\begin{align}\label{equa:Haction}
	\psi_{\lambda}(f) &= \lambda^* f \\
	\psi_{\lambda}(s_a - 1) &=  \begin{cases} -\frac{\lambda^*(a - h_a)}{\lambda^* a}(s_{\partial a} - 1)\in \Hom^{\cont}(\ha\Pol_{\lambda}, \ha\Pol_{\lambda}) & a(\lambda) = 0 \\  \frac{\lambda^*(h_a - a)}{\lambda^*a} -\frac{(s_a\lambda)^*(a - h_a)}{(s_a\lambda)^*a}s_{\partial a}\in \Hom^{\cont}(\ha\Pol_{\lambda}, \ha\Pol_{\lambda} + \ha\Pol_{s_a\lambda}) & a(\lambda)\neq 0\end{cases}.
\end{align}
\begin{lemm}
	The map $\psi$ defines a faithful continuous action of $\bbH$ on $\ha\Pol$.  \qed 
\end{lemm}

Let $\bbH^{\wedge}_{\lambda_0}\subset \End^{\cont}(\ha\Pol)$ be the closure of the image of $\psi$. It has a set of topological generators which reflects better than $\bbH$ the weight-space decomposition of objects of $\calO_{\lambda_0}(\bbH)$. For $\lambda\in W_S \lambda_0$, we define a function $\ord_{\lambda}:S^+ \to \bfZ_{\ge -1}$ by
\begin{align}\label{equa:ordH}
	\ord_{\lambda}(a) = \ord_{z = a(\lambda)}(z - h_a)z^{-1}.
\end{align}

\begin{lemm}\label{lemm:generateurH}
	The topological algebra $\bbH^{\wedge}_{\lambda_0}$ is topologically generated by the following elements:
	\begin{enumerate}
		\item
			for each $\lambda\in W_S \lambda_0$, the projector $\bfe(\lambda):\ha\Pol\to \ha\Pol_{\lambda}\subset \ha\Pol$, \\
		\item
			the polynomial ring $\bfC[V]$, which acts diagonally on $\ha\Pol$ by multiplication on each factor $\ha\Pol_{\lambda} = \bfC\llbracket V\rrbracket$,
		\item
			for each $a\in \Delta$ an operator $\tau_a = \sum_{\lambda\in W_S \lambda_0}\tau_a\bfe(\lambda):\ha\Pol \to \ha\Pol$, where
			\begin{align}\label{equa:tauH1}
				\tau_a\bfe(\lambda):\ha\Pol_{\lambda}\to \ha\Pol_{s_a\lambda}, \quad
				\tau_a f = \begin{cases} (\partial a)^{-1}(s_{\partial a}(f) - f)  & \ord_\lambda(a) = -1 \\ (\partial a)^{\ord_{\lambda}(a)}s_{\partial a}(f) & \ord_\lambda(a) \ge 0   
				\end{cases}
			\end{align}
			for $f\in \ha\Pol_{\lambda} = \bfC\llbracket V\rrbracket$, where $\partial a\in R$ is the differential of $a\in S$ and $s_{\partial a}: \bfC\llbracket V\rrbracket\to \bfC\llbracket V\rrbracket$ is the reflection with respect to the finite root $\partial a\in R$, see~\autoref{subsec:sysaff}.
	\end{enumerate}
\end{lemm}
\begin{proof}
	Let $A\subset \End^{\cont}(\ha\Pol)$ denote the closure of the subalgebra generated by the three set of operators $\bfe(\lambda)$, $\bfC[V]$ and $\tau_a$. We need to show that $A = \bbH^{\wedge}_{\lambda_0}$. \par
	Consider the restriction $\psi\mid_{\bfC[E]}$. It	factorises as
	\[
		\bfC[E]\to \prod_{\lambda\in W_S \lambda_0}\varprojlim_{k}  \bfC[E] / \frakm_{\lambda}^k\xrightarrow{\cong} \prod_{\lambda}\ha\Pol_{\lambda},
	\]
	where $\frakm_\lambda\subset\bfC[E]$ is the defining ideal of the closed point $\lambda\in E$. The Chinese remainder theorem implies that the map has dense image. In particular, $\bfe(\lambda)\in \ha\Pol_{\lambda}\subset \ha\Pol$ lies in the closure of image for each $\lambda\in W_S \lambda_0$. Therefore, $\bbH^{\wedge}_{\lambda_0}$ contains the closure of $\psi(\bfC[E])\bfe(\lambda)$ in $\Hom^{\cont}(\ha\Pol, \ha\Pol)$; the latter is equal to the algebra $\bfC\llbracket V\rrbracket\bfe(\lambda)$ which acts on $\ha\Pol_{\lambda}$ by multiplication. Thus we have $\ha\Pol \subset A$ and $\ha\Pol\subset \bbH^{\wedge}_{\lambda_0}$. It remains to show that $\{\psi(s_a)\}_{a\in \Delta}$ lies in $A$ and $\{\tau_a\}_{a\in \Delta}$ lies in $\bbH^{\wedge}_{\lambda_0}$.
	\par
	For each $a\in \Delta$ and $\lambda\in W_S \lambda_0$, by comparison of the formulae~\eqref{equa:Haction} and~\eqref{equa:tauH1}, we see that the elements $\bfe(s_a\lambda)\psi(s_a - 1)\bfe(\lambda)$ and $\tau_a\bfe(\lambda)$ generate the same cyclic left $\bfC\llbracket V\rrbracket$-submodule of $\Hom^{\cont}(\ha\Pol_{\lambda}, \ha\Pol_{s_a\lambda})$. In particular, $\tau_a\bfe(\lambda)$ lies in $\bbH^{\wedge}_{\lambda_0}$ and conversely, $\bfe(s_a\lambda)\psi(s_a - 1)\bfe(\lambda)$ lies in $A$. For $\lambda\in W_S \lambda_0$ such that $s_a \lambda = \lambda$, we have $\psi(s_a - 1)\bfe(\lambda) = \bfe(\lambda)\psi(s_a - 1)\bfe(\lambda)\in A$. For $\lambda\in W_S \lambda_0$ such that $s_a \lambda \neq \lambda$, we have
	\begin{align*}
		&\psi(s_a - 1)\bfe(\lambda) = \bfe(\lambda)\psi(s_a - 1)\bfe(\lambda) + \bfe(s_a\lambda)\psi(s_a - 1)\bfe(\lambda),\\
		&\bfe(\lambda)\psi(s_a - 1)\bfe(\lambda) = \frac{\lambda^*(h_a - a)}{\lambda^*a}, \quad
		\bfe(s_a\lambda)\psi(s_a - 1)\bfe(\lambda) = -\frac{(s_a\lambda)^*(a - h_a)}{(s_a\lambda)^*a}s_{\partial a};
	\end{align*}
	since $\lambda^*(h_a - a)/\lambda^*a\in \bfC\llbracket V\rrbracket\bfe(\lambda)\subset A$, it follows that $\psi(s_a - 1)\bfe(\lambda)\in A$. 
	Summing over the idempotents, we obtain
	\[
		\psi(s_a - 1) = \sum_{\lambda\in W_S \lambda_0}\psi(s_a - 1)\bfe(\lambda)\in A,\quad \tau_a = \sum_{\lambda\in W_S \lambda_0}\tau_a\bfe(\lambda)\in \bbH^{\wedge}_{\lambda_0}.
	\]
	The result follows.
\end{proof}

Let $\bbH^{\wedge}_{\lambda_0}\mof^{\sm}$ be the category of finitely generated $\bbH^{\wedge}_{\lambda_0}$-modules $M$ such that for each element $m\in M$, the annihilator $\ann_{\bbH^{\wedge}_{\lambda_0}}(m)$ is an open left ideal of $\bbH^{\wedge}_{\lambda_0}$. Notice that these conditions imply 
\[
	M = \bigoplus_{\lambda\in W_{S}\lambda_0} \bfe(\lambda)M\quad \text{and}\quad \dim \bfe(\lambda)M < \infty,\quad \text{for $M\in \bbH^{\wedge}_{\lambda_0}\mof^{\sm}$}.
\]

\begin{lemm}\label{lemm:equivHH}
	 The restriction $\psi^*$ yields an equivalence of categories
	 \[
			\bbH^{\wedge}_{\lambda_0}\mof^{\sm} \cong \calO_{\lambda_0}(\bbH).\pushQED{\qed}\qedhere\popQED
		\]
\end{lemm}

\subsection{Idempotent form \texorpdfstring{$\bfH_{\lambda_0}$}{H}}\label{subsec:Hlambda}
In view of~\autoref{lemm:equivHH} and~\autoref{lemm:generateurH}, in order to study the block $\calO_{\lambda_0}(\bbH)$, it is convenient to consider the subalgebra generated by the generators given in~\autoref{lemm:generateurH}. \par
Observe that the operators $\left\{ \bfe(\lambda) \right\}_{\lambda\in W_S \lambda_0}$, $\bfC[V]$ and $\{\tau_a\}_{a\in \Delta}$ preserve the dense submodule $\Pol\subset \ha\Pol$. Let $\bfH_{\lambda_0}$ \index[ch1]{H@$\bfH_{\lambda_0}$}be the associative (non-unital) subalgebra of $\End_{\bfC}(\Pol)$ generated these operators.  
Let $\bfH_{\lambda_0}\mof_0$ be the category of finitely generated $\bfH_{\lambda_0}\mof_0$-modules $M$ such that $M = \bigoplus_{\lambda\in W_S\lambda_0}\bfe(\lambda)M$ and such that the subspace $V^*\subset \bfC[V]$ acts locally nilpotently on $M$.
\begin{lemm}\label{lemm:equivHH2}
	There is a natural inclusion $\bfH_{\lambda_0}\hookrightarrow \bbH^{\wedge}_{\lambda_0}$ with dense image, which induces an equivalence of categories by pulling back the module-structure:
	\begin{align*}
			\bbH^{\wedge}_{\lambda_0}\mof^{\sm} &\xrightarrow{\cong} \bfH_{\lambda_0}\mof_0. \\
	\end{align*}
\end{lemm}
\begin{proof}
	By the density of the submodule $\Pol\subset \ha\Pol$ and~\autoref{lemm:generateurH}, there is a unique inclusion $\bfH_{\lambda_0}\hookrightarrow \bbH^{\wedge}_{\lambda_0}$ with dense image which fixes the generators $\left\{ \bfe(\lambda) \right\}_{\lambda\in W_S \lambda_0}$, $\bfC[V]$ and $\{\tau_a\}_{a\in \Delta}$. The assertion on the equivalence of category follows straightforward from the density.
\end{proof}

Combining the equivalences of~\autoref{lemm:equivHH} and~\autoref{lemm:equivHH2}, we obtain the following result:
\begin{prop}\label{prop:equivHH}
	There is an equivalence of categories 
	\[
		\calO_{\lambda_0}(\bbH) \cong \bfH_{\lambda_0}\mof_0.\pushQED{\qed}\qedhere\popQED
	\]
\end{prop}

\begin{rema}\label{rema:HetA}
	In~\autoref{sec:agha}, we will attach to each family of functions $\left\{ \omega_{\lambda} \right\}_{\lambda\in W_S \lambda_0}$ an algebra $\bfA^{\omega}$. We will study them in a larger generality. The algebra $\bfH_{\lambda_0}$ is the special case where $\omega_{\lambda} = \ord_{\lambda}$ for $\lambda\in W_S \lambda_0$.
\end{rema}

\subsection{Central subalgebra \texorpdfstring{$\calZ^{\wedge}$}{Z}}\label{subsec:Zhat}
For $\lambda\in W_S \lambda_0$, let $W_{\lambda}$ denote the stabiliser of $\lambda$ in $W_S$. The stabiliser $W_{\lambda}$ is a finite parabolic subgroup of the Coxeter group $W_S$. The affine Weyl group $W_S$ acts on the vector space $V_{\bfC}$ via the finite quotient\footnote{The notation is chosen so that $(\dW w)(\partial a) = \partial (wa)$ for $a\in S$ and $w\in W_S$ as well as $\dW s_a = s_{\partial a}$ for $a\in S$.} $\dW:W_S \to W_S / Q^{\vee}\cong W_R$. Let $\calZ^{\wedge} = \bfC\llbracket V\rrbracket^{W_{\lambda_0}}$\index[ch1]{Z@$\calZ^{\wedge},\frakm_{\calZ}$} be the ring of $W_{\lambda_0}$-invariant formal power series. Since $W_{\lambda_0}$ acts by reflections on $V$, the ring $\calZ^{\wedge}$ is a complete regular local ring. Let $\frakm_{\calZ}\subset \calZ^{\wedge}$ be the maximal ideal. \par
For each $\lambda\in W_S \lambda_0$, we define a homomorphism $\calZ^{\wedge}\to \ha\Pol_\lambda$: choosing a $w\in W_S$ such that $w\lambda_0 = \lambda$, we let $f\mapsto w(f)\in \bfC\llbracket V\rrbracket^{W_{w\lambda}}\subset \ha\Pol_{\lambda}$. This map is clearly independent of the choice of $w$ and it identify $\calZ^{\wedge}$ with the invariant subspace $\bfC\llbracket V\rrbracket^{W_{w\lambda}}$. The space $\ha\Pol$ is regarded as a $\calZ^{\wedge}$-module via the diagonal action. It is easy to observe that $\calZ^{\wedge}$ lies in the centre of $\bbH^{\wedge}_{\lambda_0}$. 
\begin{rema}
	One can show that $\calZ^{\wedge}$ coincides with the centre of $\bbH^{\wedge}_{\lambda_0}$; however, we do not need this fact.
\end{rema}

\section{Reminder on affine Hecke algebra}\label{sec:AHA}

We keep the notation $(V, R, \Delta_0)$, $\left( E,S, \Delta \right)$ and $h = \left\{ h_{a} \right\}_{a\in S}$ as above. 
\subsection{Extended affine Hecke algebras}\label{subsec:AHA}

Put
\begin{align*}
		v = \left\{v_{\alpha} \right\}_{\alpha\in R}, \quad v_{\alpha} = \begin{cases} \exp(\pi i h_{\alpha}) & \alpha\in R_{\idv} \\ \exp(\pi i h_{\alpha+1}) & \alpha\in R\cap 2R \end{cases}.
\end{align*}
Recall that $\til W^{\vee}_S = P \rtimes W_R$ is the dual extended affine Weyl group (we identify $W_R$ with $W_{R^{\vee}}$ via the correspondence $s_{\alpha}\leftrightarrow s_{\alpha^{\vee}}$). Define the extended affine braid group $\frakB_S$ \index[ch1]{B@$\frakB_S$}for the dual root system $(V^*, R^{\vee})$ to be the group generated by $T_w$ for $w\in \til W_S^{\vee}$ with the following relation for each $y,w\in \til W^{\vee}_S$:
\begin{align*}
	T_yT_w = T_{yw},\quad \text{if } \ell(yw) = \ell(y) + \ell(w).
\end{align*}

The extended affine Hecke algebra in parameters $v$, denoted by $\bbK$, is the quotient of the group algebra $\bfC \frakB_{S}$ by the following relations for $\alpha\in \Delta_0$, in the case where $R$ is reduced:
\begin{align*}
	&(T_{s_{\alpha}} - v^2_{\alpha})(T_{s_{\theta^{\vee}}} + 1) = 0,\quad
	(T_{s_0} - v^2_{\theta})(T_{s_0} + 1) = 0
\end{align*}
where $s_0\in \til W^{\vee}_S$ is the reflection with respect to the affine simple root and $\theta\in R^+$ is the highest root. In the case where $R$ is non-reduced, let $\beta\in \Delta_0$ be the simple root such that $2\beta\in R$. Let $\bbK$ be the quotient of $\bfC\frakB_S$ by the following relations for $\alpha\in \Delta_0\setminus\left\{ \beta \right\}$:
\begin{align*}
	&(T_{s_{\alpha}} - v^2_{\alpha})(T_{s_\alpha} + 1) = 0 \\
	&(T_{s_{\beta}} - v^2_{\beta}v_{\theta})(T_{s_\beta} + 1) = 0\\
	&(T_{s_{0}} - v^2_{\beta}v_{\theta}^{-1})(T_{s_{0}} + 1) = 0. \\
\end{align*}

\subsection{Bernstein--Lusztig presentation}
Choose a square root $v_{\theta}^{1/2}$ of $v_{\theta}$. Define a group homomorphism $v: \frakB_S \to \bfC^{\times}$ by setting $v(s_{\alpha}) = v_{\alpha}$ for $\alpha\in \Delta_0$ and $v(s_0) = v_{\theta}$ in the case where $R$ is reduced; $v(s_{\alpha}) = v_\alpha$, $v(s_{\beta}) = v_{\beta}v_{\theta}^{1/2}$ and $v(s_0) = v_{\beta}v_{\theta}^{-1/2}$ in the case where $R$ is non-reduced and $\beta\in \Delta_0$ with $2\beta\in R$ and $\alpha\in \Delta_0\setminus\left\{ \beta \right\}$.  \par
There is a subalgebra $\bfC P\subset \bbK$ given by $\mu\mapsto v(\mu)T_{\mu}$ for $\mu\in P\subset \til W^{\vee}_S$ dominant with respect to the basis $\Delta_0$. For $\beta\in P$ in general, we decompose it into $\beta = \beta' - \beta''$ with $\beta'$ and $\beta''$ dominant and set $Y^{\beta} = T_{\beta'}T_{\beta''}^{-1}$. Then there is a decomposition
\begin{align*}
	\bbK = H_{R} \otimes \bfC P,
\end{align*}
where $H_R$ is the subalgebra generated by $\{T_{s_\alpha}\}_{a\in \Delta_0}$ and $\bfC P$ is the subalgebra generated by $\left\{ Y^{\beta} \right\}_{\beta\in P}$, with the following commutation relations: for each $f\in \bfC P$,
\begin{align}\label{equa:bernlusz}
	&T_{s_{\alpha}}f - s_{\alpha}(f)T_{s_{\alpha}} = 
	(v^2_{\alpha} - 1)\frac{f - s_{\alpha}(f)}{1 - Y^{-\alpha}}, & \alpha\in \Delta_0,\; 2\alpha\notin R \\
	& T_{s_{\beta}}f - s_{\beta}(f)T_{s_{\beta}} = \left((v_{\beta}^2v_{\theta} - 1) + \left( v_{\beta}^2 - v_{\theta} \right)Y^{-\beta}\right)\frac{f - s_{\beta}(f)}{1 - Y^{-2\beta}}, & \beta\in \Delta_0,\; 2\beta\in R.
\end{align}

\subsection{Finite dimensional modules}
Let $T$ \index[ch1]{T@$T$}be the torus defined by $T = Q^{\vee}\otimes \bfC^{\times}$ so that $Q^{\vee} = \bfX_*(T)$ is its group of cocharacters and $P = \bfX^*(T)$ is its group of characters. We view $\bfC[T] = \bfC P$ as a subalgebra of $\bbK$. \par
For each $\ell\in T$, let $\frakm_{\ell}\subset \bfC P$ denote the defining ideal of the closed point $\ell$, which is generated by $Y^{\beta} - Y^{\beta}(\ell)\in \bfC P$ for all $\beta\in P$. Given any module $M\in \bbK\Mod$, consider for each $\ell\in T$ the generalised $\ell$-weight space in $M$ of the action of the subalgebra $\bfC P\subset \bbK$:
\begin{align*}
		M_{\ell} = \bigcup_{N\ge 0}\left\{ a\in M\;;\;\frakm_{\ell}^N a = 0\right\}.
\end{align*}
For any $\ell_0\in T$, we define $\calO_{\ell_0}\left( \bbK \right)$\index[ch1]{O@$\calO_{\ell_0}(\bbK)$} to be the full subcategory of $\bbK\mof$ consisting of those $M\in\bbK\mof$ which admit a decomposition by weight:
\begin{align*}
	M = \bigoplus_{\ell\in W \ell_0} M_{\ell}.
\end{align*}

\subsection{Idempotent form \texorpdfstring{$\bfK_{\ell_0}$}{K}}\label{subsec:Kell}
Fix $\ell_0\in T$. 
As in the case of $\bbH$, we define an algebra which is more adapted to the study of the block $\calO_{\ell_0}(\bbK)$. Define for each $\ell\in W_R \ell_0$ a polynomial ring $\Pol_\ell = \bfC[V]$ and let $\Pol = \bigoplus_{\ell\in W_R\ell_0} \Pol_{\ell}$. For each $\ell$, define $\bfe(\ell):\Pol\to \Pol$ to be the idempotent linear endomorphism of projection onto the factor $\Pol_\ell$. Let $R^+_{\idv} = R^+ \setminus 2R^+$ denote the set of indivisible positive roots.
In view of~\eqref{equa:bernlusz}, for $\ell\in W_R \ell_0$, we define a function $\ord_{\ell}:R^+_\idv \to \bfZ$: 
\begin{align*}\label{equa:ordK}
	\ord_{\ell}(\alpha) = 
	\begin{cases} 
		\ord_{z = Y^{\alpha}(\ell)}(z - v_{\alpha}^2)(z - 1)^{-1} & 2\alpha\notin R \\
		\ord_{z = Y^{\alpha}(\ell)}(z - v_\alpha^2)(z + v_{\theta})(z^2 - 1)^{-1} & 2\alpha\in R .
	\end{cases}
\end{align*}
\par

For each $\alpha\in \Delta_0$ and $\ell\in W_R \ell_0$, we define an operator $\tau_\alpha\bfe(\ell):\Pol_\ell \to \Pol_{s_{\alpha}\ell}$ by
\begin{align*}
	\tau_\alpha \bfe(\ell) = \begin{cases} \alpha^{-1}(s_{\alpha} - 1)  & \ord_{\ell}(\alpha) = -1\\ \alpha^{\ord_{\ell}(\alpha)} s_{\alpha}  & \ord_{\ell}(\alpha) \ge 0  \end{cases}.
\end{align*}
Here $s_{\alpha}: \bfC[V]\to \bfC[V]$ is the reflection with respect to $\alpha$.  \par

Let $\bfK_{\ell_0}$\index[ch1]{K@$\bfK_{\ell_0}$} be the associative subalgebra of $\End_{\bfC}(\Pol)$ generated by $f\bfe( \ell )$ and $\tau_{\alpha}\bfe(\ell)$ for $f\in \bfC[V]$, $\alpha\in \Delta_0$ and $\ell\in W_R \ell_0$. Let $\bfK_{\ell_0}\mof_0$ be the category of finitely generated $\bfK_{\lambda_0}\mof_0$-modules $M$ such that the subspace $V^*\subset \bfC[V]$ acts locally nilpotently on $M$.  Same arguments as~\autoref{lemm:equivHH} and~\autoref{lemm:equivHH2} show that:
\begin{prop}\label{prop:equivKK}
	There is an equivalence of categories
	\[
		\calO_{\ell_0}\left( \bbK \right) \cong \bfK_{\ell_0}\mof_0. \pushQED{\qed}\qedhere\popQED
	\]
\end{prop}

\begin{rema}\label{rema:KetB}
	In~\autoref{sec:MQHA}, we will attach to each family of functions $\left\{ \omega_{\ell} \right\}_{\ell\in W_R \ell_0}$ an algebra $\bfB^{\omega}$. The algebra $\bfK_{\lambda_0}$ is the special case of $\bfB^{\omega}$ with $\omega_{\ell} = \ord_{\ell}$ for $\ell\in W_R \ell_0$.
\end{rema}

\section{The monodromy functor \texorpdfstring{$\bbV$}{V}}\label{subsec:monodromie}
In this section, we review the construction of the monodromy functor of~\cite{VV04}, which is a trigonometric analogue of the Knizhnik--Zamolodchikov functor introduced in~\cite{GGOR03} for rational Cherednik algebras. We prove in~\autoref{prop:quotM} that this functor is a quotient functor. \par
Keep the notation $(E, S, \Delta)$ and $a_0\in \Delta$ as above. In addition, we fix $\lambda_0 \in E_{\bfC}$. Consider the following exponential map
\begin{align}\label{equa:exp}
	E_{\bfC}\cong V_{\bfC} = Q^{\vee}\otimes \bfC &\xrightarrow{\exp} Q^{\vee}\otimes \bfC^{\times} = T\\
	\mu\otimes r&\mapsto \mu\otimes e^{2\pi i r}.\nonumber
\end{align}\index[ch1]{exp@$\exp$} 
Put $\ell_0 = \exp(\lambda_0)$. 
For simplifying the notation, denote $\calC_0 = \calO_{\lambda_0}\left( \bbH \right)$ and $\calB_0 = \calO_{\ell_0}\left( \bbK \right)$.\index[ch1]{C@$\calC_0$}\index[ch1]{B@$\calB_0$} \par

\subsection{Dunkl operators}

Consider the dual torus $T^{\vee} = P\otimes \bfC^{\times}$\index[ch1]{T@$T^{\vee}$}. The ring of regular functions $\bfC[T^{\vee}]$ is isomorphic to the group algebra of the coroot lattice $\bfC Q^{\vee}$:
\begin{align*}
	\bfC Q^{\vee} &\xrightarrow{\cong} \bfC[T^{\vee}] \\
	Q^{\vee}\ni\mu &\mapsto X^{\mu}
\end{align*}
\par

For each $\xi\in V^*$, let $\partial_\xi\in \Gamma\left(T^{\vee}, \calT_{T^{\vee}}  \right)^{T^{\vee}}$ be the translation-invariant vector field on $T^{\vee}$ such that $\partial_\xi\mid_{e} = \xi$ under the isomorphism $\calT_{T^{\vee}}\mid_e \cong V$. We view $\partial_{\xi}$ as a linear differential operator on $T^{\vee}$, so that $\partial_{\xi}(X^{\mu}) = \langle \xi,\mu\rangle X^{\mu}$ for each $\mu\in Q^{\vee}$.
\par

The regular part of $T^{\vee}$ is defined as $T^{\vee}_{\circ} = \bigcap_{\alpha\in R^+}\left\{X^{\alpha^{\vee}} \neq 1 \right\}\subset T^{\vee}$. Let $\calD(T^{\vee}_{\circ})$~\index[ch1]{D@$\calD(T^{\vee}_{\circ})$} denote the ring of algebraic differential operators on $T^{\vee}_{\circ}$.  \par

For $\xi\in V^*$, the trigonometric Dunkl operator $D_{\xi}:\bfC[T^{\vee}]\to \bfC[T^{\vee}]$ is the $\bfC$-linear operator defined as follows: 
\begin{align*}
	D_{\xi}(f) = \partial_\xi(f) - \sum_{\alpha\in R^+} h_\alpha\langle \xi, \alpha^{\vee}\rangle \frac{f - s_{\alpha}(f)}{1 - X^{-\alpha^{\vee}}} + \langle\xi,\rho^{\vee}_h\rangle f,\quad \rho^{\vee}_h = \frac{1}{2}\sum_{\alpha\in \Delta^+}h_{\alpha}\alpha^{\vee}\in V_{\bfC}.
\end{align*}
We consider $D_{\xi}$ as an element of $\calD(T^{\vee}_{\circ})\rtimes W_R$.\par

According to~\cite[4.1]{VV04}, the following homomorphism of $\bfC$-algebras
\begin{align*}
	\bfC[T^{\vee}]\otimes \bfC W_R\otimes \bfC[V]  = \bbH&\to \calD(T^\vee_{\circ})\rtimes W_R \\
	X^{\mu} \otimes w\otimes 1 &\mapsto X^{\mu}\otimes w\\
	1\otimes 1\otimes\xi &\mapsto D_{\xi} \\
\end{align*}
extends to an isomorphism $\bfC[T^{\vee}_{\circ}]\otimes_{\bfC[T^{\vee}]}\bbH\cong \calD(T^\vee_{\circ})\rtimes W_R$.

\subsection{Monodromy functor \texorpdfstring{$\bbV$}{V}}

Let $[T^{\vee}_{\circ} / W_R]$ be the quotient stack. According to~\cite[2.5]{heckman90}, there is an isomorphism between the orbifold fundamental group $\pi_1([T^{\vee}_{\circ} / W_R])$ and the extended affine braid group $\frakB_{S}$ from~\autoref{subsec:AHA}.

If $M\in \calO_{\lambda_0}(\bbH)$, then 
\begin{align*}
		M_{\circ} = \bfC[T^{\vee}_{\circ}]\otimes_{\bfC[T^{\vee}]}M
\end{align*}
is a $W$-equivariant $\calD(T^{\vee}_{\circ})$-module, which is in fact an integrable connection with regular singularities. Therefore the monodromy representation on the vector space of flat sections of $M$ on (the universal covering of) the orbifold $[T^{\vee}_{\circ} / W]$ defines a $\frakB_S$-module, which is denoted by $\bbV(M)$. It is shown in~\cite[5.1]{VV04} that the $\frakB_S$-action on $\bbV(M)$ factorises through the surjective algebra homomorphism $\bfC\frakB_S\to \bbK$ and yields an exact functor 
\begin{align*}
	\bbV: \calO_{\lambda_0}(\bbH)\to  \calO_{\ell_0}(\bbK).
\end{align*}
\index[ch1]{V@$\bbV$}

\subsection{Central actions of \texorpdfstring{$\calZ^{\wedge}$}{Z} intertwined by \texorpdfstring{$\bbV$}{V}}\label{subsec:ZM}
For convenience, we denote $\calC_0 = \calO_{\lambda_0}(\bbH)$ and $\calB_0 = \calO_{\ell_0}(\bbK)$. Recall the central subalgebra $\calZ^{\wedge} = \bfC\llbracket V\rrbracket^{W_{\lambda_0}}$ defined in~\autoref{subsec:Zhat}. Let $\rmZ(\calC_0) = \End(\id_{\calC_0})$ and  $\rmZ(\calB_0) = \End(\id_{\calB_0})$ denote the categorical centres. \par
Let $W_{\lambda_0}$ be the stabiliser of $\lambda_0\in E_{\bfC}$ in $W_S$ and let $W_{\ell_0}$ be the stabiliser of $\ell_0\in T$ in $W_R$. Let $\bar\lambda_0$ be the image of $\lambda_0$ in $E_{\bfC} / W_{\lambda_0}$ and let $\bar\ell_0$ be the image of $\ell_0$ in $T / W_{\ell_0}$. 
The exponential map~\eqref{equa:exp} induces an analytic map
\begin{align*}
	\exp^{\lambda_0}:E_{\bfC}/ W_{\lambda_0} \to T/  W_{\ell_0} , 
\end{align*}
which is locally biholomorphic near $\bar\lambda_0$. The push-forward along $\exp^{\lambda_0}$ at $\bar \lambda_0$ yields an isomorphism of complete local rings
\begin{align*}
		\exp^{\lambda_0}_*:\calO_{E_{\bfC}/ W_{\lambda_0}, \bar\lambda_0}^{\wedge}\xrightarrow{\sim}  \calO_{T/ W_{\ell_0}, \bar\ell_0}^{\wedge}. 
\end{align*}
Note that $\calZ^{\wedge} \cong \calO_{E_{\bfC}/ W_{\lambda_0}, \bar\lambda_0}^{\wedge}$.
For each $w\in W_S$, the action of $w$ on $E_{\bfC}$ and on $T$ (the latter via the quotient map $\dW$ from~\autoref{subsec:Zhat}) induces
\begin{align*}
		w_*: \calO_{E_{\bfC}/ W_{\lambda_0}, \bar\lambda_0}^{\wedge} \xrightarrow{\sim} \calO_{E_{\bfC}/ W_{w\lambda_0}, w\bar\lambda_0},\quad w_*:\calO_{T/ W_{\ell_0}, \bar\ell_0}^{\wedge}\cong \calO_{T/ W_{w\ell_0}, w\bar\ell_0}^{\wedge}.
\end{align*}
We define homomorphisms $\calZ^{\wedge}\to \rmZ\left( \calC_0 \right)$ and $\calZ^{\wedge}\to \rmZ\left( \calB_0 \right)$ as follows: for any $M\in \calC_0$, we decompose $M = \bigoplus_{\lambda\in W_S\, \lambda_0} M_{\lambda}$ and for each $\lambda = w \lambda_0$, an element $f\in \calZ^{\wedge}$ acts by $w_*f$ on $M_{\lambda}$. This depends only on the weight $\lambda$ but not on the choice of $w$. Similarly, for any $N\in \calB_0$, we decompose $N = \bigoplus_{\ell\in W_R\, \ell_0} N_{\ell}$. For each $\ell = w \ell_0$, an element $f\in \calZ^{\wedge}$ acts by multiplication by $w_*\exp^{\lambda_0}_*f$ on $N_{\ell}$.  \par
\begin{lemm}\label{lemm:ZM}
	The functor $\bbV:\calC_0\to \calB_0$ intertwines the $\calZ^{\wedge}$-actions on $\calC_0$ and $\calB_0$.
\end{lemm}
\begin{proof}
	Recall that the graded affine Hecke algebra is the subalgebra 
	\begin{align*}
		\ubar\bbH = \bfC W_R\otimes \Sym V^*_{\bfC} \subset \bbH. 
	\end{align*}
	For each weight $\lambda\in V_{\bfC}$, let $\calO_{\lambda}(\ubar \bbH)$ be the category of finite-dimensional $\ubar\bbH$-modules on which the action of the polynomial part $\Sym V_{\bfC}$ has weights lying in the orbit $W_S\lambda\subset V_{\bfC}$. \par
	There is a functor of induction 
	\begin{align*}
		\Ind^{\bbH}_{\ubar\bbH}:\ubar\bbH\mof\to \bbH\mof,\quad  \Ind^{\bbH}_{\ubar\bbH}M= \bbH\otimes_{\ubar\bbH} M
	\end{align*}
	and for each weight $\lambda\in E_{\bfC}$, it restricts to
	\begin{align*}
		\Ind^{\bbH}_{\ubar\bbH}:\calO_{\lambda}(\ubar \bbH)\to \calO_{\lambda}(\bbH)
	\end{align*}
	Let $\calI\subset \calC_0$ denote the essential image of $\Ind^{\bbH}_{\ubar\bbH}$. It is known that $\calI$ generates $\calC_0$ --- indeed, the module $P(\lambda)_n = \bbH / \bbH\cdot \frakm_{\lambda}^n$ lies in $\calI$ and the family $\left\{ P(\lambda)_n \right\}_{n\in \bfN,\; \lambda\in W_S \lambda_0}$ generate $\calC_0$. Therefore, it suffices to show that the restriction $\bbV\mid_{\calI}$ intertwines the actions of $\calZ^{\wedge}$.  We shall apply the deformation argument from~\cite[5.1]{VV04} to check this statement. \par

	Let $\calO = \bfC\llbracket\varpi\rrbracket$ and let $\calK = \bfC(\!(\varpi)\!)$. Let $\varepsilon \in V^{*}_{\bfC}$ be any regular coweight and put $\lambda_{0,\calO} = \lambda_0 + \varpi \varepsilon\in V^{*}_{\calO}$. Put $\ubar\bbH_{\calO} = \ubar\bbH\otimes \calO$ and $\bbK_{\calO} = \bbK\otimes \calO$. For each $\lambda_{\calO} \in W_S \lambda_{0,\calO}$ and for $n\in \bfZ_{\ge 1}$, let
	\begin{align*}
			\frakm_{\lambda_{\calO}} &= \left\langle \beta_{\calO} - \langle\beta_{\calO}, \lambda_{\calO}\rangle\;;\; \beta\in V_{\calO}\right\rangle\subset \Sym_{\calO} V^{*}_{\calO},\quad \frakm_{\lambda_{\calK}} = \frakm_{\lambda_{\calO}}[\varpi^{-1}],\\
			\bbS_{\lambda_{\calO}^n} &= \Sym_{\calO} V^{*}_{\calO} / \frakm_{\lambda_{\calO}}^n, \quad \bbS_{\lambda_{\calK}^n} =\bbS_{\lambda_{\calO}^n}[\varpi^{-1}], \\
		\ubar P(\lambda_{\calO})_n &= \ubar\bbH_{\calO}  \otimes_{\Sym_{\calO} V^{*}_{\calO}}\bbS_{\lambda_{\calO}^n},\quad \ubar P(\lambda_{\calK})_n = \ubar\bbH_{\calK}[\varpi^{-1}].
	\end{align*}
	Note that all these objects are flat over $\calO$.
	Let $\ubar P(\lambda_{\calO})_n^{\nabla}$ be the space of flat sections of the affine Knizhnik--Zamolodchikov equation (AKZ) on the constant vector bundle on $T^{\vee}_{\circ}$ of fibre $\ubar P(\lambda_{\calO})_n$. The monodromy representation yields a $\bbK_{\calO} = \bbK\otimes\calO$ action on $\ubar P(\lambda_{\calO})_n^{\nabla}$. \par
	Since the stabiliser of $\lambda_{\calO}$ in $W_S$ is trivial, there is an eigenspace decomposition
	\begin{align*}
			\ubar P(\lambda_{\calK})_n = \bigoplus_{w\in W_R}\left(\ubar P(\lambda_{\calK})_n\right)_{w\lambda_{\calK}}, \quad \left(\ubar P(\lambda_{\calK})_n\right)_{w\lambda_{\calK}} = b_{w}\bbS_{\lambda_{\calK}^n},
	\end{align*}
	where each $b_{w}\bbS_{\lambda_{\calK}^n} $ is a free $\bbS_{\lambda_{\calK}^n}$-module of rank $1$. Consider the boundary point of $T^{\vee}$:
	\[
		\lim_{n\rightarrow +\infty} \exp(n i\rho) = \left(X^{\alpha^{\vee}} = 0\right)_{\alpha\in R^+},\quad \text{where $\rho = (1/2)\sum_{\alpha\in R^+}\alpha$}.
	\]
	Applying the Frobenius method around this point, we obtain a fundamental solution $\left\{ b^{\nabla}_w \right\}_{w\in W_R}$ of the AKZ equation on $T^{\vee}_{\circ}$ which satisfies
	\begin{align*}
			b^{\nabla}_w(\exp(\mu)) &= e^{-2\pi i(\langle \mu, \rho^{\vee}_h\rangle + \mu)}\cdot \left(b_w + G(\mu)\right) \\
		&\text{for }\mu\in V^{*}_{\bfC} \text{ such that } \Im\langle \mu,\alpha^{\vee}\rangle \gg 0,\; \forall \alpha\in \Delta^+,
	\end{align*}
	where $G(\mu)$ is a $\ubar P(\lambda_{\calK})_n$-valued analytic function in $\mu$ with such that 
	\begin{align*}
		G(\mu)\to 0\quad \text{when }\Im\langle \alpha, \mu\rangle\to +\infty, \forall \alpha\in \Delta^+.
	\end{align*}
		\par
		The fundamental solution induces an $\bbS_{\lambda_{\calK}^n}$-linear isomorphism
		\begin{equation}\label{equa:isomP}
				\ubar P(\lambda_{\calK})_n \xrightarrow{\sim} \ubar P(\lambda_{\calK})^{\nabla}_n,\quad  b_w\mapsto b_w^{\nabla}.
		\end{equation}
		Under this isomorphism, the monodromy operator on the right-hand side corresponding to $\beta\in X$ is identified with $e^{2\pi i \beta}$ on the left-hand side. Put
		\begin{align*}
			\calZ^{\wedge}_{\calO} = \left(\left(\Sym V^*_{\calO}\right)^{W_{\lambda_0}}\right)^{\wedge}_{\bar\lambda_{0, \calO}},\quad \calZ^{\wedge}_{\calK} = \calZ^{\wedge}_{\calO}[\varpi^{-1}]\cong \left(\Sym V^*_{\calK}\right)^{\wedge}_{\bar\lambda_{0, \calO}}.
		\end{align*}
		We define the action of $\calZ^{\wedge}_{\calO}$ and $\calZ^{\wedge}_{\calK}$ on $\ubar\bbH_{\calO}$-modules and $\ubar\bbH_{\calK}$-modules in a similar way. \par
		Since the action of $\calZ^{\wedge}_{\calK}$ on $\ubar P(\lambda_{\calK})_n$ coincides with the action of the polynomial part $\Sym V^*_{\calK}\subset \ubar\bbH$ up to twists by elements of $W_R$, the induced action of $\calZ^{\wedge}_{\calK}$ on the $\bbK_{\calK}$-module $\ubar P(\lambda_{\calK})^{\nabla}_n$ is identified with the exponentiation of the action of $\calZ^{\wedge}_{\calK}$ on the $P(\lambda_{\calK})_n$ under~\eqref{equa:isomP}. \par
		Since the $\calO$-lattices $\ubar P(\lambda_{\calO})_n\subset \ubar P(\lambda_{\calK})_n$ and $\ubar P(\lambda_{\calO})^{\nabla}_n\subset \ubar P(\lambda_{\calK})^{\nabla}_n$ are stable under the action of the subring $\calZ^{\wedge}_{\calO}\subset \calZ^{\wedge}_{\calK}$, the functor $M\mapsto M^{\nabla}$ also intertwines the two $\calZ^{\wedge}_{\calO}$-actions. Put $\ubar P\left( \lambda \right)_n = \ubar P\left( \lambda_{\calO} \right)_n\otimes_{\calO}\bfC$. Then $\ubar P(\lambda)_n\mapsto \ubar P(\lambda)^{\nabla}_n = \bbV(P(\lambda)_n)$ also intertwines the two $\calZ^{\wedge}$-actions. Finally, since the family of modules $\ubar P(\lambda)_n$ for $\lambda\in W_S \lambda_0$ and $n\ge 1$ generates the category $\calO_{\lambda_0}\left( \ubar \bbH \right)$, the functor $\bbV$ restricted to $\calI$ intertwine the $\calZ^{\wedge}$-actions as asserted.
\end{proof}

\subsection{Completion of categories}
Since the affine Hecke algebra $\bbK$ is of finite rank over its centre, namely $(\bfC P)^W$, $\calB_0 = \calO_{\ell_0}(\bbK)$ is equivalent to the category of modules of finite length over some semi-perfect algebra. It is also the case for $\calC_0 = \calO_{\lambda_0}(\bbH)$.  In particular, they are both noetherian-artinian. Consider the category of pro-objects\footnote{The basic properties of categories of pro-objects are reviewed in~\autoref{sec:pro}.} $\Pro(\calC_0)$ and $\Pro(\calB_0)$. 
We have two central actions introduced in~\autoref{subsec:ZM}
	\begin{align*}
		\calZ^{\wedge}\to \End(\id_{\calC_0} ) \cong  \End(\id_{\Pro(\calC_0)} ) \\
		\calZ^{\wedge}\to \End\left( \id_{\calB_0} \right) \cong \End(\id_{\Pro(\calB_0)} ).
	\end{align*}
	By~\autoref{lemm:ZM}, the functor $\bbV:\calC_0\to \calB_0$ intertwines these $\calZ^{\wedge}$-actions. The extension $\bbV: \Pro(\calC_0)  \to  \Pro(\calB_0)$ still intertwines the $\calZ^{\wedge}$-actions.  \par

Define $\calC\subset \Pro(\calC_0)$ to be the subcategory consisting of objects $M\in \Pro(\calC_0)$ such that $M/\frakm_{\calZ}^k M\in \calC_0$ for all $k\ge 0$. Similarly we define $\calB\subset \Pro(\calB_0)$ to be the subcategory consisting of objects $N\in \Pro(\calB_0)$ such that $M/\frakm_{\calZ}^k M\in \calC_0$ for all $k\ge 0$.\index[ch1]{C@$\calC$}\index[ch1]{B@$\calB$} \par
	\begin{lemm}\label{lemm:PC}
		For each simple object $L\in \calC_0$ (resp. $L\in \calB_0$), its projective cover $\calP(L)\in \Pro(\calC_0)$ (resp. $\calP(L)\in \Pro(\calB_0)$) lies in $\calC$ (resp. $\calB$). 
	\end{lemm}
	\begin{proof}
		Notice that by the general result~\autoref{lemm:envproj}, the objects of $\calC_0$ (resp. $\calB_0$) admit projective covers in $\Pro(\calC_0)$ (resp. $\Pro(\calB_0)$).
		The statement holds obviously for $\calB_0$ because $\bbK$ is of finite rank over its centre. For $\calC_0$, by~\autoref{prop:equivHH}, there is an equivalence $\calC_0\cong \bfH_{\lambda_0}$ and the algebra $\bfH_{\lambda_0}$ is Morita-equivalent to an algebra of finite rank over its centre, \cfauto{subsec:Agmod}.
	\end{proof}

	\begin{lemm}\label{lemm:extM}
	The functor $\bbV:\Pro(\calC_0)\to\Pro(\calB_0)$ restricts to $\bbV:\calC\to \calB$.
\end{lemm}
\begin{proof}
	If $M\in \calC$, then $M / \frakm_{\calZ}^k M\in \calC_0$ and by~\autoref{lemm:ZM}, $\bbV(M) / \frakm_{\calZ}^k \bbV(M)\cong \bbV(M / \frakm_{\calZ}^k M)\in \calB_0$. It follows that $\bbV(M)\in \calB$.
\end{proof}

\subsection{Right adjoint of \texorpdfstring{$\bbV$}{V}}
Recall that $\calB_0 = \calO_{\ell_0}(\bbK)$ and $\calC_0 = \calO_{\lambda_0}(\bbH)$. 
\begin{lemm}\label{lemm:adjV}
	The functor $\bbV:\calC_0\to \calB_0$ admits a right adjoint functor $\bbV^{\top}: \calB_0\to \calC_0$.
\end{lemm}
\begin{proof}
	We first define a functor $\bbV^{\top}:\calB_0\to \Ind(\calC_0)$ with natural isomorphisms 
	\begin{align}\label{equa:adj}
		\Hom_{\calB_0}\left( \bbV(M), N \right)\cong\Hom_{\Ind(\calC_0)}\left( M, \bbV^{\top}(N) \right)
	\end{align}
	for $M\in \calC_0$ and $N\in \calB_0$. For any $N\in \calB_0$, let 
	\begin{align*}
		F_N: \calC_0^{\op}\to \bfC\Mod, \quad F_N:M\mapsto \Hom_{\calB_0}\left( \bbV(M), N \right)
	\end{align*}
	and let 
	\begin{align*}
		F_N(M)^{\min} = F_N(M)\setminus \bigcup_{0\neq M'\subset M} F_N(M/M').
	\end{align*}
	Here, we regard $F_N(M/M')$ as a subspace of $F_N(M)$ by the right exactness of $F_N$. Let $\calI_N$ be the category whose objects are pairs $(M, a)$, where $M\in \calC_0$ and $a\in F_N(M)^{\min}$, and whose morphisms are defined by
	\begin{align*}
		\Hom_{\calI_N}\left( (M, a), (M', a') \right) = \left\{ f\in \Hom_{\calC_0}(M, M')\;;\; F_N(f)(a') = a \right\}.
	\end{align*}
	We set
	\begin{align*}
		\bbV^{\top}(N) = \indlim{(M, a)\in \calI_N} M\in \Ind(\calC_0).
	\end{align*}
	According to~\cite[3.5, Lemma 6]{serre02}, $\bbV^{\top}(N)$ represents the functor $F_N$, so $\bbV^{\top}$ satisfies the desired adjoint property~\eqref{equa:adj}. \par
	Now we show that in fact the object $\bbV^{\top}(N)$ in $\Ind(\calC_0)$ lies in the subcategory $\calC_0$. Let $\calP_{\calC}\in \calC$ be the sum of all projective indecomposable objects (up to isomorphism) of $\calC$ so that for any $M\in \calC_0$, the dimension of $\Hom_{\calC}(\calP_{\calC}, M)$ is equal to the length of $M$. Since $\bbV(\calP_{\calC})\in \calB$ is a finitely generated $\bbK$-module, the vector space $\Hom_{\calB}\left( \bbV(\calP_{\calC}), N\right)$ is finite-dimensional. On the other hand, there are isomorphisms
	\begin{align}\label{equa:proind}
		\varinjlim_{\substack{M\subset \bbV^{\top}(N) \\ M\in \calC_0}}\Hom_{\calC}\left( \calP_{\calC},  M\right) &\cong \varinjlim_{\substack{M\subset \bbV^{\top}(N) \\ M\in \calC_0}}\varinjlim_{\substack{\calQ\subset \calP_{\calC} \\ \calP_{\calC}/\calQ\in \calC_0}}\Hom_{\calC_0}\left( \calP_{\calC}/\calQ,  M\right) \\
		&\cong \varinjlim_{\substack{\calQ\subset \calP_{\calC} \\ \calP_{\calC}/\calQ\in \calC_0}}\Hom_{\Ind(\calC_0)}\left( \calP_{\calC}/\calQ,  \bbV^{\top}(N)\right) \nonumber\\
		&\cong\varinjlim_{\substack{\calQ\subset \calP_{\calC} \\ \calP_{\calC}/\calQ\in \calC_0}}\Hom_{\calB_0}\left(  \bbV(\calP_{\calC}/\calQ), N\right)\nonumber\\
		& \cong \Hom_{\Pro(\calB_0)}{\Bigg(} \prolim{\substack{\calQ\subset \calP_{\calC} \\ \calP_{\calC}/\calQ\in \calC_0}} \bbV(\calP_{\calC}/\calQ), N{\Bigg)} \cong \Hom_{\calB}\left( \bbV(\calP_{\calC}), N\right).\nonumber
	\end{align}
	The first and the fourth isomorphisms are due to~\eqref{equa:hom-pro} of~\autoref{sec:pro}; the second one is exchanging the order of the two colimits and it holds due to the definition of morphisms between ind-objects; the third one is due to~\eqref{equa:adj}; the last one is due to~\autoref{lemm:extM}. \par
	Since $N\in \calB_0$, there is some integer $n$ such that $\frakm_{\calZ}^n N = 0$. Since $\bbV(\calP_{\calC})\in \calB$, the quotient $\bbV(P_{\calC}) / \frakm_{\calZ}^n\bbV(P_{\calC})$ lies in $\calB_0$. Thus the $\Hom$-space
	\begin{align*}
		\Hom_{\calB}(\bbV(P_{\calC}), N)\cong \Hom_{\calB_0}(\bbV(P_{\calC})/\frakm_{\calZ}^n\bbV(P_{\calC}), N)
	\end{align*}
	is finite-dimensional. The above isomorphisms~\eqref{equa:proind} imply that the length of the subobjects ${M\subset \bbV^{\top}(N)}$ such that $M\in \calC_0$ is bounded. It follows that $\bbV^{\top}(N)$ lies in $\calC_0$ by~\autoref{lemm:pro}~(iii). Thus $\bbV^{\top}: \calB_0\to \calC_0$ is a right adjoint to $\bbV$.
	
\end{proof}

\subsection{\texorpdfstring{$\bbV$}{V} is a quotient functor}\label{subsec:quotM}
\begin{prop}\label{prop:quotM}
The monodromy functor $\bbV:\calC_0\to \calB_0$ is a quotient functor. 
\end{prop}
\begin{proof}
	Recall that $\calD(T^\vee_{\circ})$ is the ring of algebraic linear differential operators on the regular part $T^{\vee}_{\circ}$ of the dual torus $T^{\vee} = P\otimes \bfC^{\times}$.
	By construction, the functor $\bbV$ factorises into the following
	\[
		\begin{tikzcd}[row sep=5mm]
			\bbH\Mod \arrow{r}{\mathrm{loc}} & \calD\left( T^{\vee}_{\circ} \right)\rtimes W_R\Mod &  \\
			\calC_0 \ar[rounded corners,"\bbV",to path={|- (\tikztotarget)[near end]\tikztonodes}]{drr} \arrow{r} \arrow[hookrightarrow]{u} & \conn_{W_R}^{\mathrm{rs}}\left( T^{\vee}_{\circ} \right) \arrow{r}{\mathrm{RH}} \arrow[hookrightarrow]{u} & \bfC \frakB_S\mof^{\rm{fini}} \\
			& & \calB_0 \arrow[hookrightarrow]{u}
		\end{tikzcd}
	\]
	where $\conn_{W_R}^{\mathrm{rs}}\left( T^{\vee}_{\circ} \right)$ is the subcategory of $\calD\left( T^{\vee}_{\circ} \right)\rtimes W_R\mof$ consisting of $W_R$-equivariant integrable connections on $T^{\vee}_{\circ}$ which have regular singularities along the boundary. 
	The arrow in the first line is the localisation functor $\mathrm{loc} = \bfC[T^{\vee}_{\circ}]\otimes_{\bfC[T^{\vee}]}\relbar$, whose right adjoint $\mathrm{loc}^{\top}$ is the restriction of the action of $\bbH_{\circ} = \calD\left( T^{\vee}_{\circ} \right)\rtimes W_R$ to $\bbH$. The restriction of $\mathrm{loc}$ to $\calC_0$ factorises through the inclusion of subcategory 
	\[
		\conn_{W_R}^{\mathrm{rs}}\left( T^{\vee}_{\circ} \right) \hookrightarrow \calD\left( T^{\vee}_{\circ} \right)\rtimes W_R\Mod
	\]
	and gives the first arrow of the second line. The functor $\mathrm{RH}$ is the Riemann--Hilbert correspondence (the Knizhnik--Zamolodchikov equations have regular singularities~\cite{opdam00}), due to Deligne~\cite[2.17+5.9]{deligne70}, between algebraic connections with regular singularities and finite-dimensional representations of the fundamental group $\pi_1\left( \left[ T^{\vee}_{\circ} / W_R\right] \right)\cong \frakB_S$. \par

	We show that $\bbV$ admits a section functor in the sense of Gabriel~\cite{gabriel62}. We have shown in~\autoref{lemm:adjV} that $\bbV$ admits a right adjoint functor $\bbV^{\top}$. The functor $\bbV^{\top}$ can be described as follows:
	\begin{align*}
		\calB_0 \hookrightarrow \bfC \frakB_S\mof^{\mathrm{fini}}\cong \conn^{\mathrm{rs}}_{W_R}\left( T^{\vee}_{\circ} \right)\to \calC_0,
	\end{align*}
	where the last arrow is the functor which sends an object $M\in \conn^{\mathrm{rs}}_{W_R}(T^{\vee}_{\circ})$ to the biggest $\bbH$-submodule of $M$ which lies in $\calC_0$, denoted by $M\mid_{\calC_0}\subset M$. We show that the adjunction counit $\bbV\circ\bbV^{\top} \to \id_{\calB_0}$ is an isomorphism. We first show that it is a monomorphism: for any $M\in \conn^{\mathrm{rs}}_{W_R}(T^{\vee}_{\circ})$, we have $\bfC[T^{\vee}_{\circ}]\otimes_{\bfC[T^{\vee}]}M\cong M$; by the flatness of $\bfC[T^{\vee}_{\circ}]$ over $\bfC[T^{\vee}]$, the inclusion $M\mid_{\calC_0}\hookrightarrow M$ gives rise to a monomorphism 
	\[
		\bfC[T^{\vee}_{\circ}]\otimes_{\bfC[T^{\vee}]}\left(M\mid_{\calC_0}\right)\to \bfC[T^{\vee}_{\circ}]\otimes_{\bfC[T^{\vee}]}M\cong M; 
	\]
	composing it with the Riemann--Hilbert correspondence, we see that $\bbV\circ \bbV^{\top}\to \id_{\calB_0}$ is a monomorphism. \par
	Let $N\in \calB_0$. By the exactness of $\bbV$, to show that the adjunction counit $\bbV\,\bbV^{\top}N\hookrightarrow N$ is an isomorphism, it suffices to find an $\bbH$-submodule of $\mathrm{RH}^{-1}(N)$ whose localisation to $T^{\vee}_{\circ}$ is equal to $\mathrm{RH}^{-1}(N)$. There exists a surjection
	\begin{align*}
		\bigoplus_{i\in \calI} P\left( \ell_i \right)_{n_i}\to N
	\end{align*}
	where $\calI$ is an index set and $P\left( \ell_{i} \right)_{n_{i}} = \bbK / \bbK\cdot \frakm_{\ell_{i}}^{n_{i}}$. By~\cite[5.1~(i)]{VV04}, for each $i\in \calI$ there is an induced module $P\left( \lambda_{i} \right)_{n_{i}} = \bbH / \bbH\cdot \frakm_{\lambda_{i}}^{n_{i}}\in \calC_0$ such that $\exp(\lambda_i) = \ell_i$ and $\bbV\left(P\left( \lambda_{i} \right)_{n_{i}}\right)\cong P\left( \ell_{i} \right)_{n_{i}}$. Hence the image of $P(\lambda_{i})_{n_{i}}$ in $\mathrm{RH}^{-1}(N)$ is an $\bbH$-submodule which satisfies the requirement. We conclude that $\bbV\circ\bbV^{\top}\cong \id_{\calB_0}$; therefore $\bbV^{\top}$ is a section functor for $\bbV$. \par

	By the criterion of Gabriel~\cite[3.2, Prop 5]{gabriel62}, $\bbV$ is a quotient functor. 
\end{proof}

\section{Comparison of \texorpdfstring{$\bbV$}{V} and \texorpdfstring{$\bfV$}{V}} \label{sec:VV}
\subsection{The functors \texorpdfstring{$\bbV$}{V} and \texorpdfstring{$\bfV$}{V}}
In~\autoref{partII}, we will study the idempotent forms $\bfH_{\lambda_0}$ and $\bfK_{\ell_0}$ in a broader context, ~\cfauto{rema:HetA} and~\autoref{rema:KetB}. Specifically, in~\autoref{subsec:V}, we will introduce a quotient functor for graded modules $\bfV:\bfH_{\lambda_0}\gmod\to \bfK_{\ell_0}\gmod$. It has an ungraded version $\bfV:\bfH_{\lambda_0}\mof_0 \to \bfK_{\lambda_0}\mof_0$. On the other hand, by~\autoref{prop:equivHH} and~\autoref{prop:equivKK}, we have equivalences of categories $\calO_{\lambda_0}(\bbH)\cong \bfH_{\lambda_0}\mof_{0}$ and $\calO_{\ell_0}(\bbK)\cong \bfK_{\lambda_0}\mof_{0}$. The situation can be depicted in a diagram:
\[
	\begin{tikzcd}[row sep=1em]
		\calO_{\lambda_0}(\bbH)\arrow{r}{\bbV}\arrow{d}{\cong} & \calO_{\ell_0}(\bbK) \arrow{d}{\cong} \\
		\bfH_{\lambda_0}\mof_0\arrow{r}{\bfV} & \bfK_{\ell_0}\mof_0 \\
	\end{tikzcd}
\]
\begin{conj}
There is an isomorphism of functors $\bbV\cong\bfV$. 
\end{conj}
In the rest of this section, we use results from~\autoref{partII} to prove a weaker version of this statement. 

\subsection{Comparison of the kernels}
By~\autoref{prop:quotM} and~\autoref{subsec:V}, the functors $\bbV$ and $\bfV$ are already known to be quotient functors. The following proposition generalises a result from~\cite{liu19}, where the geometric construction of the dDAHA was used.
\begin{prop}\label{prop:compker}
	The kernels $\ker \bbV$ and $\ker \bfV$ are identified via the equivalence $\calO_{\lambda_0}(\bbH)\cong \bfH_{\lambda_0}\mof_{0}$.
\end{prop}
\begin{proof}
	Let $F:\calO_{\lambda_0}(\bbH)\xrightarrow{\sim}\bfH_{\lambda_0}\mof_0$ denote the equivalence from~\autoref{prop:equivHH}. We show that for every object $M\in\calO_{\lambda_0}(\bbH)$, the condition~\autorefitem{theo:kerV}{iii} for $F M$ implies $\bbV M = 0$. Let $M = \bigoplus_{\lambda\in W_S \lambda_0} M_{\lambda}$ be the decomposition into generalised weight spaces of $\bfC[E]$ and let
	\begin{align*}
			M_{\le t} = \bigoplus_{\substack{\lambda\in W_S\lambda_0 \\ \|\lambda\| \le t}} M_{\lambda}, \quad \text{for } t\in \bfR_{\ge 0}.
	\end{align*}
	Note that under the equivalence $F$, the generalised weight space $M_{\lambda}$ is identified with $\bfe(\lambda) F(M_{\lambda})$. Following the same arguments as in the proof~\ref{theo:kerV-iii}$\Rightarrow$\ref{theo:kerV-iv} of~\autoref{theo:kerV}, we have $s_a M_{t}\le M_{t + \delta}$ for every $t\in \bfR_{\ge 0}$ and $a\in \Delta$. Let $U = \bfC[E]^{\le 1} + \sum_{a\in \Delta}\bfC\cdot s_a\subset\bbH$ so that $U$ generates $\bbH$ as $\bfC$-algebra. Then, by the assumption~\ref{theo:kerV-iii}, we see that for each finite-dimensional subspace $L\subset M$ and each $\varepsilon > 0$, 
	\begin{align*}
		\lim_{n\to \infty}\dim\left( U^n L \right) / n^{r - 1 + \varepsilon}  = 0,\quad r = \rk R.
	\end{align*}
	Hence we obtain $\dim_{\GK,\bbH} M \le r - 1$, and in particular $\dim_{\GK, \bfC[T^{\vee}]} M \le r - 1$ for the subalgebra $\bfC[T^{\vee}] = \bfC Q^{\vee}\subset \bbH$. As the algebra $\bfC[T^{\vee}]$ is commutative and by~\autoref{prop:coherenceT}, $M$ is coherent over $\bfC[T^{\vee}]$, the Gelfand--Kirillov dimension of $M$ coincides with the Krull dimension of the subvariety $\supp_{T^{\vee}} M\subset T^{\vee}$. As the localisation of $M$ on the regular part $T^{\vee}_{\circ}$ must be locally free, we see that it must be zero since $\dim T^{\vee}_{\circ} = r > \dim \supp M$. Hence $\bbV M = 0$ by the definition of $\bbV$. We see that $\ker \bfV\subset F(\ker \bbV)$. \par
	Since $\bfV$ and $\bbV$ are both quotient functors on noetherian-artinian categories, by comparison of the rank of the Grothendieck groups
	\begin{align*}
		\rk K_0\left( \ker \bbV \right) &= \rk K_0\left( \calO_{\lambda_0}(\bbH) \right) - \rk K_0\left( \calO_{\ell_0}(\bbK) \right)  \\
		&= \rk K_0\left( \bfH_{\lambda_0}\mof_0 \right) - \rk K_0\left( \bfK_{\ell_0}\mof_0 \right) = \rk K_0\left( \ker \bfV \right),
	\end{align*}
	we see that $\ker\bfV = F(\ker\bbV)$.
\end{proof}

\part{Quiver Hecke algebras}\label{partII}

\section{Quiver double Hecke algebra}\label{sec:agha}
Fix an irreducible based finite root system $(V, R, \Delta_0)$ and let $(E, S, \Delta)$ be its affinisation. In this section we will also abbreviate $P = P_R$, $Q = Q_R$, $P^{\vee} = P^{\vee}_R$ and $Q^{\vee} = Q^{\vee}_R$. 

\subsection{The polynomial matrix algebra \texorpdfstring{$\bfA^o$}{A}}
Fix once and for all $\lambda_0\in E$\index[ch2]{lambda@$\lambda_0$}. Define for each $\lambda\in W_S \lambda_0$ a polynomial ring $\Pol_\lambda = \bfC[V]$\index[ch2]{Pol@$\Pol_{\lambda}$} and let $\Pol_{W_S \lambda_0} = \bigoplus_{\lambda\in W_S\lambda_0} \Pol_{\lambda}$\index[ch2]{Pol@$\Pol_{W_S\lambda_0}$}. For each $\lambda$, define $\bfe(\lambda):\Pol_{W_S\lambda_0}\twoheadrightarrow \Pol_{\lambda}\subset \Pol_{W_S\lambda_0}$\index[ch2]{e@$\bfe(\lambda)$} to be the projection onto the factor $\Pol_\lambda$. \par
  
For each $a\in \Delta$, define an operator $\tau^o_a:\Pol_{W_S \lambda_0}\to \Pol_{W_S \lambda_0}$ by
\begin{align*}
		\tau^o_a = \sum_{\lambda\in W_S \lambda_0}\tau^o_a \bfe(\lambda),\quad \tau^o_a \bfe(\lambda):\Pol_{\lambda_0} \to \Pol_{s_a\lambda_0},\\
		\tau^o_a \bfe(\lambda) = \begin{cases} (\partial a)^{-1}(s_{\partial a} - 1)  & a(\lambda) = 0 \\ s_{\partial a} & a(\lambda) \neq 0   \end{cases}.
\end{align*}
Here $\partial a\in R$ is the differential of $a\in S$, \cfauto{subsec:sysaff}. \par
Let $\bfA^o = \bfA^o(E, S, \Delta, \lambda_0)$\index[ch2]{A@$\bfA^o$} be the associative (non-unital) subalgebra of $\End_{\bfC}(\Pol_{W_S\lambda_0})$ generated by $f\bfe( \lambda )$ and $\tau^o_{a}\bfe(\lambda)$ for $f\in \bfC[V]$, $a\in \Delta$ and $\lambda\in W_S \lambda_0$
	. \par

	\subsection{Centre \texorpdfstring{$\calZ$}{Z}}\label{subsec:Z}
For $\lambda\in W_S \lambda_0$, let $W_{\lambda}$\index[ch2]{W@$W_{\lambda}$} be the stabiliser of $\lambda$ in $W_S$. The stabiliser $W_{\lambda}$ is a finite parabolic subgroup of the Coxeter group $W_S$. The affine Weyl group $W_S$ acts on the vector space $V$ via the finite quotient $\dW:W_S \to W_S / Q^{\vee}\cong W_R$\index[ch2]{@$\partial$}. Let $\calZ = \bfC[V]^{W_{\lambda_0}}$\index[ch2]{Z@$\calZ,\frakm_{\calZ}$} be the ring of $W_{\lambda_0}$-invariant polynomials, graded by the degree of monomials. Since $W_{\lambda_0}$ acts by reflections on $V$, the ring $\calZ$ is a graded polynomial ring. Let $\frakm_{\calZ}\subset \calZ$ be the unique homogeneous maximal ideal. \par
	For each $\lambda\in W_S \lambda_0$, we define a homomorphism $\calZ\to \Pol_\lambda$: choosing a $w\in W_S$ such that $w\lambda_0 = \lambda$, we let $f\mapsto w(f)\in \bfC[V]^{W_{\lambda}}\subset \Pol_{\lambda}$. This map is clearly independent of the choice of $w$ and it identifies $\calZ$ with the invariant subspace $\bfC[V]^{W_{\lambda}}$. The infinite sum $\Pol_{W_S \lambda_0}$ is regarded as a $\calZ$-module via the diagonal action. \par
	The following are standard results from the invariant theory for reflection groups:
	\begin{prop}\label{prop:Ao}
		The following statements hold:
		\begin{enumerate}
			\item
				For each $\lambda\in W_S \lambda_0$, the $\calZ$-module $\Pol_{\lambda}$ is free of rank $\# W_{\lambda} = \# W_{\lambda_0}$.
			\item
				For any $w\in W_S$, choose a reduced expression $w = s_{a_l} \cdots s_{a_1}$ and put $\tau^{o}_w\bfe(\lambda) = \tau^{o}_{a_l}\cdots \tau^{o}_{a_1}\bfe(\lambda)$ for each $\lambda\in W_S\lambda_0$. Then the element $\tau^o_w\bfe(\lambda)$ is independent of the choice of the reduced expression for $w$ and moreover, there is a decomposition
				\begin{align*}
					\Hom_{\calZ}\left( \Pol_\lambda, \Pol_{W_S\lambda_0} \right) = \bigoplus_{w\in W_S} \tau^o_w \bfC[V]\bfe(\lambda).
				\end{align*}
			\item
				The $\bfA^o$-action on $\Pol_{W_S\lambda_0}$ commutes with $\calZ$ and yields an isomorphism
				\begin{align*}\label{equa:Mat}
						\bfA^o\xrightarrow{\sim} \bigoplus_{\lambda\in W_S \lambda_0}\Hom_{\calZ}(\Pol_\lambda, \Pol_{W_S \lambda_0})
				\end{align*}
				\qed 
		\end{enumerate}
	\end{prop}
	\subsection{Subalgebras \texorpdfstring{$\bfA^{\omega}$}{} of \texorpdfstring{$\bfA^o$}{A}}\label{subsec:Aomega}
	Let $\omega = \left\{ \omega_{\lambda}\right\}_{\lambda\in W_S\lambda_0}$\index[ch2]{omega@$\omega_\lambda$} be a family of functions $\omega_{\lambda}:S^+\to \bfZ_{\ge -1}$ satisfying the following properties:
\begin{enumerate}
	\item
		$\omega_\lambda(a) = -1$ implies $a(\lambda) = 0$;
	\item
		for $w\in W_S$ and $b\in S^+\cap w^{-1}S^+$ we have $\omega_{\lambda}(b) = \omega_{w\lambda}(wb)$. 
\end{enumerate}
One may extend $\omega_{\lambda}$ to a function $\til\omega_{\lambda}:S\to \bfZ_{\ge -1}$ by choosing $w\in W_S$ such that $wa\in S^+$ and setting $\til\omega_{\lambda}(a) = \omega_{w\lambda}(wa)$\index[ch2]{omega@$\til\omega_{\lambda}$}. We require $\omega$ to satisfy the following property:
\begin{enumerate}[resume*]
	\item
		For some (thus every) $\lambda\in W_S\lambda_0$, the extended function $\til\omega_{\lambda}:S\to \bfZ_{\ge -1}$ has finite support.
\end{enumerate}
We call the family $\left\{ \omega_{\lambda} \right\}_{\lambda \in W_S\lambda_0}$ a family of order functions. The order functions can be characterised as follows:
\begin{lemm}
	Every family of order functions $\left\{ \omega_\lambda \right\}_{\lambda\in W_S \lambda_0}$ is determined by the $W_{\lambda_0}$-invariant finitely supported function $\til \omega_{\lambda_0}: S\to \bfZ_{\ge -1}$ satisfying 
	\[
		\til\omega_{\lambda_0}(a) = -1\Rightarrow a(\lambda_0) = 0 \quad \forall a\in S.  \pushQED{\qed}\qedhere\popQED
	\]
\end{lemm}

Define an operator $\tau^{\omega}_a = \sum_{\lambda\in W_S \lambda_0}\tau^{\omega}_{a}\bfe(\lambda)\in \End_{\calZ}(\Pol_{W_S \lambda_0})$ with $\tau^\omega_a\bfe(\lambda):\Pol_{\lambda}\to \Pol_{s_a\lambda}$\index[ch2]{tau@$\tau^{\omega}_a$} by setting
\begin{align*}
	\tau^\omega_a\bfe(\lambda) = 
	\begin{cases} (\partial\alpha)^{-1}(s_{\partial a} - 1) & \omega_{\lambda}(a) = -1 \\
		(\partial\alpha)^{\omega_{\lambda}(a)} s_{\partial a} & \omega_{\lambda}(a) \ge 0 \\
	\end{cases}
\end{align*}
so that $\tau^{\omega}_a\bfe(\lambda)\in \bfA^{o}$. 

\begin{defi}\label{def:A}
	The quiver double Hecke algebra\footnote{ In this definition, the assumption that $\lambda_0\in E$ plays no essential role. We could have asked $\lambda_0$ to belong to some set on which $W_S$ acts transitively with finite parabolic stabiliser subgroups. However, the euclidean geometry of $E$ will facilitate some arguments.}  $\bfA^{\omega} = \bfA(E, S, \Delta, \lambda_0, \omega)$\index[ch2]{A@$\bfA^{\omega}$} is defined to be the subalgebra of $\bfA^o$ generated by $\bfC[V]\bfe(\lambda)$ and $\tau^{\omega}_{a}\bfe(\lambda)$ for $\lambda\in W_S$ and $a\in \Delta$. 
\end{defi}
We also introduce the rational function field and its matrix algebra: 
\begin{align*}
		\Rat_{\lambda} &= \Frac\Pol_{\lambda}= \Pol_{\lambda}\otimes_{\calZ}\Frac\calZ,\quad \Rat = \bigoplus_{\lambda\in W_S \lambda_0}\Rat_{\lambda}\\
		\bfA^{-\infty} &= \bigoplus_{\lambda\in W_S \lambda_0}\Hom_{\Frac\calZ}(\Rat_\lambda, \Rat)=\bfA^{o}\otimes_{\calZ}\Frac\calZ, \quad \tau_a^{-\infty} = s_a,
\end{align*}\index[ch2]{Rat@$\Rat_{\lambda},\Rat$} 
where $\Frac$ means the field of fractions.

\begin{exam}\label{exam:Aomega}
	\noindent
	\begin{enumerate}
		\item
			Let $o = \left\{ a\mapsto -\delta_{a(\lambda)=0}\right\}_{\lambda\in W_S\lambda_0}$ denote the smallest family of order functions. We recover the matrix algebra $\bfA^o$. \par
		\item
			Let $\omega = \left\{ 0 \right\}_{\lambda\in W_S \lambda_0}$ be the zero constant function. Then $\bfA^{\omega} = \Pol_{W_{S}\lambda_0}\rtimes W_S$ is the skew tensor product. If $W_{\lambda_0} = 1$, then $\bfA^{\omega} = \bfC[V]\wr W_S$ is the wreath product.
		\item
			Let $E = \bfR$, let $\epsilon$ be the coordinate function on $\bfR$ and let $S = \{\pm 2\epsilon\} + \bfZ$, so that $(E, S)$ is the affine root system of type $A^{(1)}_1$. Choose the basis $\Delta = \{a_1 = 2\epsilon, a_0 = 1 - 2\epsilon\}$. The affine Weyl group $W_S$ is generated by $s_0$ and $s_1$, where $s_1$ (resp. $s_0$) is the orthogonal reflection with respect to $0\in E$ (resp. $1/2\in E$). Set $\lambda_0 = 1/4\in E$, so that $W_S \lambda_0 = 1/4 + (1/2)\bfZ$ and $W_{\lambda_0} = 1$. It follows that $\Pol_{\lambda} = \bfC[\epsilon]$ for all $\lambda\in W_S \lambda_0$ and $\bfA^o$ is the matrix algebra over $\bfC[\epsilon]$ of rank $W_S \lambda_0$. \par
			Set
			\begin{align*}
				\til\omega_{\lambda_0}(a) = \begin{cases} 1 & a\in \Delta \\ 0 & a\in S \setminus\Delta \end{cases}
			\end{align*}
			and define the family of order functions $\omega = \left\{ \omega_{\lambda} \right\}_{\lambda\in W_S \lambda_0}$ by $\omega_{w\lambda_0}(a) = \til\omega_{\lambda_0}(w^{-1}a)$. It follows that $\bfA^{\omega}$ is equal to the idempotent form of the dDAHA $\bfH_{\lambda_0}$ introduced in~\autoref{subsec:Hlambda} with parameter $h_a = 1/2$ for all $a\in S$. We can depict the algebra $\bfA^{\omega}$ with the following diagram:
			\[
				\begin{tikzcd}
					\cdots & \Pol_{-3/4} \arrow[bend left]{r}{s} \arrow[phantom]{r}{\tau_1} & \arrow[bend left]{l}{s}\Pol_{3/4} \arrow[bend left]{r}{s} \arrow[phantom]{r}{\tau_0} & \arrow[bend left]{l}{-\epsilon\, s} \Pol_{1/4}\arrow[bend left]{r}{\epsilon\, s}\arrow[phantom]{r}{\tau_1} &\arrow[bend left]{l}{s} \Pol_{-1/4} \arrow[bend left]{r}{s}\arrow[phantom]{r}{\tau_0} & \arrow[bend left]{l}{s}\Pol_{5/4}& \cdots, \\
				\end{tikzcd}
			\]
			where $s:\bfC[\epsilon]\to \bfC[\epsilon]$ is given by the substitution $\epsilon \mapsto -\epsilon$. 
	\end{enumerate}
\end{exam}
\begin{rema}
	We may view $\bfA^{\omega}$ as an affinisation of the quiver Hecke algebra $R_{\beta}(\Gamma)$ attached to a certain quiver $\Gamma = (I, H)$ and a dimension vector $\beta\in \bfN I$, \cfauto{rema:KLR}. The parameter $\omega$ is an analogue of the polynomials $Q_{i,j}(u, v)$ in Rouquier's definition of quiver Hecke algebras.
\end{rema}
\begin{rema}
	Following~\cite[\S 2.3]{webster19koszul}, one can write down a complete list of relations between the generators $\tau^{\omega}_{a}\bfe( \lambda)$, $\bfC[V]\bfe(\lambda)$ for the algebra $\bfA^{\omega}$ in the manner of Khovanov--Lauda--Rouquier. The most sophisticated is the braid relation between pairs of generators from $\{\tau^{\omega}_{a}\bfe( \lambda)\}_{a\in \Delta,\lambda\in W_S \lambda_0}$. We will only prove a weaker version of it in~\autoref{lemm:braid}, which is enough for our needs. 
\end{rema}

\subsection{Filtration by length}\label{subsec:filtlg}

\begin{defi}\label{def:FA}
	We define the filtration by length $\{F_{\le n}\bfA^{\omega}\}_{n\in \bfN}$ on $\bfA^{\omega}$ by
\begin{align*}
	F_{\le n}\bfA^{\omega} = \sum_{\lambda\in W_S \lambda_0}\sum_{k =0}^n\sum_{\left( a_1, \ldots, a_k \right)\in \Delta^k} \bfC[V]\tau^{\omega}_{a_1}\cdots \tau^{\omega}_{a_k}\bfe(\lambda).
\end{align*}
\index[ch2]{FA@$F_{\le n}\bfA^{\omega}$}
\end{defi}

In general, it is hard to express the operators $\tau^{\omega}_{a_1}\cdots\tau^{\omega}_{a_k}$; however, the leading term is easy to describe.
\begin{lemm}\label{prop:formula}
	Let $w = s_{a_l} \cdots s_{a_1}$ be a reduced expression and let $\lambda\in W_{S}\lambda_0$. Then
	\begin{enumerate}
		\item\label{prop:formula-i}
			For any $f\in \bfC[V]$, and any family of order functions $\omega$ there is a commutation relation:
			\begin{align*}
				f\tau^{\omega}_{a_l}\cdots\tau^{\omega}_{a_1}\bfe(\lambda) \equiv \tau^{\omega}_{a_l}\cdots\tau^{\omega}_{a_1}\, w^{-1}(f)\bfe(\lambda)\mod F_{\le l - 1} \bfA^{\omega}.
			\end{align*}
		\item\label{prop:formula-ii}
			For any pair of families of order functions $\omega$ and $\omega'$ such that $\omega \le \omega'$ (pointwise), there is a congruence relation:
			\begin{align*}
				\tau^{\omega'}_{a_l}\cdots\tau^{\omega'}_{a_1}\bfe(\lambda) &\equiv \tau^{\omega}_{a_l}\cdots\tau^{\omega}_{a_1}\left(\prod_{b\in S^+\cap w^{-1}S^-}(-\partial b)^{\omega'_{\lambda}(b) - \omega_{\lambda}(b)}\right)\bfe(\lambda)\mod F_{\le l - 1} \bfA^{\omega}.
			\end{align*}
	\end{enumerate}
\end{lemm}
\begin{proof}
	We prove the statement~\ref{prop:formula-i} by induction on the length $l = \ell(w)$. It is trivial for $l = 0$. For $l = 1$:
	\begin{align}\label{equa:comm}
		(f\tau^{\omega}_a - \tau^{\omega}_a\, s_{\partial a}(f))\bfe(\lambda)  =
		\begin{cases}
			( \partial  a)^{-1}(s_{\partial a}(f) - f)\bfe(\lambda) & \omega_{\lambda}(a) = -1 \\
			0 & \omega_{\lambda}(a) \ge 0 \\
		\end{cases}
	\end{align}
	It belongs to $F_{\le 0}\bfA^{\omega}\bfe(\lambda) = \bfC[V]\bfe(\lambda)$ in both cases. \par
	For $l > 1$, by the induction hypothesis, we get
	\begin{align*}
		&(f\tau^{\omega}_{a_l}\cdots \tau^{\omega}_{a_1} - \tau^{\omega}_{a_l}\cdots \tau^{\omega}_{a_1}\, w^{-1}(f))\bfe(\lambda)  \\
		&=(f\tau^{\omega}_{a_l} - \tau^{\omega}_{a_l}\, s_{a_l}(f))\tau^{\omega}_{a_{l-1}}\cdots \tau^{\omega}_{a_1}\bfe(\lambda)  \\
		&+ \tau^{\omega}_{a_l}(s_{a_l}(f)\tau^{\omega}_{a_{l-1}}\cdots \tau^{\omega}_{a_1} - \tau^{\omega}_{a_{l-1}}\cdots \tau^{\omega}_{a_1}\, w^{-1}(f))\bfe(\lambda) \in F_{\le l-1}\bfA^{\omega},
	\end{align*}
	whence~\ref{prop:formula-i}. \par

	We prove~\ref{prop:formula-ii} by induction on $l = \ell(w)$. Put $w' = s_{a_{l-1}}\cdots s_{a_1}$ and $\lambda' = w' \lambda$. Then
	\begin{align*}
		&\tau^{\omega'}_{a_l}\cdots\tau^{\omega'}_{a_1}\bfe(\lambda) = (\partial a_l)^{\omega'_{\lambda'}(a_l) - \omega_{\lambda'}(a_l)}\tau^{\omega}_{a_l}\tau^{\omega'}_{a_{l-1}}\cdots\tau^{\omega'}_{a_1}\bfe(\lambda) \\
		&= \left((\partial a_l)^{\omega'_{\lambda'}(a_l) - \omega_{\lambda'}(a_l)}\tau^{\omega}_{a_l} - \tau^{\omega}_{a_l}(-\partial a_l)^{\omega'_{\lambda'}(a_l) - \omega_{\lambda'}(a_l)}\right)\tau^{\omega'}_{a_{l-1}}\cdots\tau^{\omega'}_{a_1}\bfe(\lambda)\\
		&+ \tau^{\omega}_{a_l}(-\partial a_l)^{\omega'_{\lambda'}(a_l) - \omega_{\lambda'}(a_l)}\tau^{\omega'}_{a_{l-1}}\cdots\tau^{\omega'}_{a_1}\bfe(\lambda).
	\end{align*}
	By~\eqref{equa:comm}, the first term belongs to $F_{\le l-1}\bfA^{\omega}$; the second term, by the statement~\ref{prop:formula-i} for $w' = a_{i_{l-1}}\cdots a_{i_1}$, satisfies
	\begin{align*}
		\tau^{\omega}_{a_l}(-\partial a_l)^{\omega'_{\lambda'}(a_l) - \omega_{\lambda'}(a_l)}\tau^{\omega'}_{a_{l-1}}\cdots\tau^{\omega'}_{a_1}\bfe(\lambda) &\equiv \tau^{\omega}_{a_l}\tau^{\omega'}_{a_{l-1}}\cdots\tau^{\omega'}_{a_1}\,w'^{-1}\left((-\partial a_l)^{\omega'_{\lambda'}(a_l) - \omega_{\lambda'}(a_l)}\right)\bfe(\lambda) \\
		&= 
		\tau^{\omega}_{a_l}\tau^{\omega'}_{a_{l-1}}\cdots\tau^{\omega'}_{a_1}\left(-\partial (w'^{-1}a_l)\right)^{\omega'_{\lambda}(w'^{-1} a_l) - \omega_{\lambda}(w'^{-1} a_l)}\bfe(\lambda).
	\end{align*}
	Here we have used the hypothesis that $\omega_{\lambda}(w'^{-1} a_l) = \omega_{\lambda'}(a_l)$. 
	Using  the induction hypothesis, we obtain
	\begin{align*}
		&\tau^{\omega'}_{a_l}\cdots\tau^{\omega'}_{a_1}\bfe(\lambda) \equiv \tau^{\omega}_{a_l}\tau^{\omega'}_{a_{l-1}}\cdots\tau^{\omega'}_{a_1}\left((-\partial (w'^{-1}a_l))^{\omega'_{\lambda}(w'^{-1} a_l) - \omega_{\lambda}(w'^{-1} a_l)}\right)\bfe(\lambda) \\
		&\equiv\tau^{\omega}_{a_l}\cdots\tau^{\omega}_{a_1}\left(\prod_{b\in S^+\cap w^{-1}S^-}(-\partial b)^{\omega'_{\lambda}(b) - \omega_{\lambda}(b)}\right)\bfe(\lambda).
	\end{align*}
	The last equation is  due to the relation $S^+\cap w^{-1} S^- = S^+\cap w'^{-1} S^-\cup\left\{w'^{-1} a_l \right\}$. This proves~\ref{prop:formula-ii}. 
\end{proof}

\subsection{Basis theorem}\label{subsec:basis}
We aim to prove an analogue of~\autoref{prop:Ao} for the subalgebra $\bfA^{\omega}\subset \bfA^{o}$. 
\begin{lemm}[braid relation]\label{lemm:braid}
	For any family of ordered functions $\left\{ \omega_{\lambda} \right\}_{\lambda\in W_S \lambda_0}$, the images of the operators $\tau^{\omega}_a\bfe(\lambda)$ in $\Gr^F\bfA^{\omega}$ satisfies the braid relations: for $a,b\in\Delta$ with $a\neq b$, let $m_{a,b}$ be the order of $s_as_b$ in $W_S$. If $m_{a,b}\neq \infty$, then
	\begin{align*}
		\underbrace{\tau^{\omega}_a\tau^{\omega}_b\tau^{\omega}_a\cdots}_{m_{a,b}}\bfe(\lambda) \equiv \underbrace{\tau^{\omega}_b\tau^{\omega}_a\tau^{\omega}_b\cdots}_{m_{a,b}}\bfe(\lambda)\quad \mod F_{\le m_{a,b} - 1}\bfA^{\omega}.
	\end{align*}
\end{lemm}
\begin{proof}
	The statement is empty for $m_{a, b}= \infty$, so we assume $m_{a, b}\neq \infty$. Let $W_{a, b}\subset W_S$ be the parabolic subgroup generated by $s_a$ and $s_b$, let $w_0\in W_{a, b}$ be the longest element and let $S_{a, b}\subset S$ be the sub-root system spanned by $a$ and $b$. Let $\bfA^{\omega}_{a, b}$ be the subalgebra of $\bfA^{\omega}$ generated by $\bfC[V]\bfe(\lambda)$, $\tau^{\omega}_a\bfe(\lambda)$ and $\tau^{\omega}_b\bfe(\lambda)$ for $\lambda\in W_S\lambda_0$ and let $F_{\le n}\bfA^{\omega}_{a, b}$ be the filtration by length defined as in~\autoref{def:FA}. It suffices to show the following
	\begin{align*}
		\underbrace{\tau^{\omega}_a\tau^{\omega}_b\tau^{\omega}_a\cdots}_{m_{a,b}}\bfe(\lambda) \equiv \underbrace{\tau^{\omega}_b\tau^{\omega}_a\tau^{\omega}_b\cdots}_{m_{a,b}}\bfe(\lambda)\quad \mod F_{\le m_{a,b} - 1}\bfA^{\omega}_{a, b}.
	\end{align*}
	because there is an inclusion $F_{\le m_{a, b} - 1}\bfA^{\omega}_{a, b}\subset F_{\le m_{a,b} - 1}\bfA^{\omega}$. An analogue of~\autoref{prop:formula} is valid for this subalgebra with the filtration $F_{\le n}\bfA^{\omega}_{a, b}$.
		\par
		We first prove the braid relation for the family $\omega' = \{\omega'_\lambda\}_{\lambda\in W_S \lambda_0}$, where $\omega'_\lambda(c) = \max\{\omega_{\lambda}(c), 0\}$. Since $\omega'_{\lambda}(c)\ge 0$ for all $c\in S^{+}_{a, b}$, the braid relation for $\tau^{\omega'}_a$ and $\tau^{\omega'}_b$ follows from the following formula (with similar proof as~\autorefitem{prop:formula}{ii}):
	\begin{align*}
		&\underbrace{\tau^{\omega'}_a\tau^{\omega'}_b\tau^{\omega'}_a\cdots}_{m_{a,b}}\bfe(\lambda) = \underbrace{s_{\partial a}s_{\partial b}s_{\partial b}\cdots}_{m_{a,b}}\prod_{c\in S^{+}_{a, b}}(-\partial c)^{\omega'_{\lambda}(c)}\bfe(\lambda).
	\end{align*}
	Let $\frakd = \prod_{\substack{c\in S^+_{a, b} \\ \omega_c(\lambda) = -1}}(\partial c)$.
	By~\autorefitem{prop:formula}{ii}, we have 
	\begin{align*}
		\underbrace{\tau^{\omega'}_a\tau^{\omega'}_b\tau^{\omega'}_a\cdots}_{m_{a,b}}\bfe(\lambda) \equiv \underbrace{\tau^{\omega}_a\tau^{\omega}_b\tau^{\omega}_a\cdots}_{m_{a,b}}\frakd\, \bfe(\lambda)  \mod F_{\le m_{a, b} - 1}\bfA^{\omega}_{a,b}.\\
	\end{align*}
	Write $X = (\underbrace{\tau^{\omega}_a\tau^{\omega}_b\tau^{\omega}_a\cdots}_{m_{a,b}} - \underbrace{\tau^{\omega}_b\tau^{\omega}_a\tau^{\omega}_b\cdots}_{m_{a,b}})\bfe(\lambda)$, so that $X\cdot \frakd\in \bfe(w_0\lambda)\left(F_{\le m_{a, b} -1}\bfA^{\omega}_{a, b}\right)\bfe(\lambda)$. Moreover, by~\autorefitem{prop:formula}{ii}, we have
	\begin{align*}
		X \equiv (\underbrace{\tau^{o}_a\tau^{o}_b\tau^{o}_a\cdots}_{m_{a,b}} - \underbrace{\tau^{o}_b\tau^{o}_a\tau^{o}_b\cdots}_{m_{a,b}})\prod_{c\in S^{+}_{a, b}}(-\partial c)^{\omega_{\lambda}(c) - o_{\lambda}(c)}\bfe(\lambda) \mod F_{\le m_{a, b} - 1}\bfA^o_{a, b}
	\end{align*}
	However, the elements $\tau^{o}_a\bfe(\lambda)$ satisfy the braid relations in $\bfA^o_{a, b}$ by~\autorefitem{prop:Ao}{ii}. It follows that $X\in F_{\le m_{a, b}-1}\bfA^o_{a, b}$ (notice that $\bfA^{\omega}_{a, b}\subseteq \bfA^{o}_{a, b}$).
	We claim that for $0\le j \le m_{a, b} - 1$, the quotient $F_{\le j}\bfA^{o}_{a, b}\bfe(\lambda) / F_{\le j}\bfA^{\omega}_{a, b}\bfe(\lambda)$ is right $\frakd$-torsion-free. This will imply that $X\in F_{\le m_{a, b} - 1}\bfA^{\omega}_{a, b}$ and complete the proof.
	\par
	We prove the claim by induction on $j$. For $j = 0$, this is obvious since $F_{\le 0}\bfA^o_{a, b} = F_{\le 0}\bfA^{\omega}_{a, b}$. Assume $j\in [1,m_{a,b}-1]$. The quotient $\Gr^{F}_{j} \bfA^{\omega}_{a, b}\bfe(\lambda)$ is spanned over $\bfC[V]$ by $\underbrace{\tau^{\omega}_a\tau^{\omega}_b\tau^{\omega}_a\cdots}_{j}\bfe(\lambda)$ and $\underbrace{\tau^{\omega}_b\tau^{\omega}_a\tau^{\omega}_b\cdots}_{j}\bfe(\lambda)$ since any non-reduced word in $a, b$ of length $ \le j$ contains consecutive letters $aa$ or $bb$ and since $(\tau^{\omega}_{a})^2, (\tau^{\omega}_{b})^2\in F_{\le 1}\bfA^{\omega}_{a, b}$. Similarly, $\Gr^{F}_{j} \bfA^{o}_{a, b}\bfe(\lambda)$ is spanned over $\bfC[V]$ by $\underbrace{\tau^{o}_a\tau^{o}_b\tau^{o}_a\cdots}_{j}\bfe(\lambda)$ and $\underbrace{\tau^{o}_b\tau^{o}_a\tau^{o}_b\cdots}_{j}\bfe(\lambda)$. Moreover, by~\autoref{prop:Ao}, $\Gr^{F}_{j} \bfA^{o}_{a, b}\bfe(\lambda)$ is free of rank $2$ over $\bfC[V]$. Denote $w = \underbrace{s_as_bs_a\cdots }_j$. Since $\omega \ge o$, by~\autorefitem{prop:formula}{ii}, we have
	\begin{align*}
		\underbrace{\tau^{\omega}_a\tau^{\omega}_b\tau^{\omega}_a\cdots}_{j}\equiv \underbrace{\tau^{o}_a\tau^{o}_b\tau^{o}_a\cdots}_{j} \left(\prod_{c\in S^+_{a, b}\cap w^{-1}S^-_{a, b}}(-\partial c)^{\omega_\lambda(c) - o_\lambda(c)}\right)\mod F_{\le j - 1}\bfA^o_{a, b}.
	\end{align*}
	The prime factors of $\frakd$ are $\partial c$ for $c\in S^+_{a, b}$ such that $\omega_{\lambda}(c) = -1$. Therefore $\frakd$ and the product
	\begin{align*}
		\prod_{c\in S^+_{a, b}\cap w^{-1}S^-_{a, b}}(-\partial c)^{\omega_\lambda(c) - o_\lambda(c)}
	\end{align*}
	are relatively prime. The same argument applies to the other product $\tau^{\omega}_b\tau^{\omega}_a\tau^{\omega}_b\cdots$. \par
	It follows that $\Gr^F_{j}\bfA^{\omega}_{a, b}\bfe(\lambda)$ and $\Gr^F_{j}\bfA^{o}_{a, b}\bfe(\lambda)$ are both free over $\bfC[V]$ of rank $2$, and the matrix representing the $\bfC[V]$-linear map $\varphi:\Gr^F_{j}\bfA^{\omega}_{a, b}\bfe(\lambda)\to \Gr^F_{j}\bfA^{o}_{a, b}\bfe(\lambda)$ (which is induced from the inclusion $\bfA^{\omega}_{a, b}\bfe(\lambda)\subset \bfA^{o}_{a, b}\bfe(\lambda)$) is diagonal with entries prime to $\frakd$. Hence $\coker \varphi$ is $\frakd$-torsion free. The snake lemma yields a short exact sequence
	\begin{align*}
		0\to \frac{F_{\le j-1}\bfA^{o}_{a, b}\bfe(\lambda)}{F_{\le j-1}\bfA^{\omega}_{a, b}\bfe(\lambda)} \to \frac{F_{\le j}\bfA^{o}_{a, b}\bfe(\lambda)}{F_{\le j}\bfA^{\omega}_{a, b}\bfe(\lambda)}\to \coker \varphi\to 0,
	\end{align*}
	in which the first term is also $\frakd$-torsion-free by induction hypothesis, and so is the middle term $\frakd$-torsion-free, whence the claim is proven.
\end{proof}

\begin{theo}\label{theo:basis}
	For each $w\in W_S$, choose a reduced expression $w = s_{a_l} \cdots s_{a_1}$ and put $\tau^{\omega}_w = \tau^{\omega}_{a_l}\cdots \tau^{\omega}_{a_1}$. Then there is a decomposition
	\begin{align*}
		\bfA^{\omega} = \bigoplus_{\lambda\in W_S \lambda_0}\bigoplus_{\substack{w\in W_S}} \bfC[V]\tau^{\omega}_w\bfe(\lambda).
	\end{align*}
\end{theo}
\begin{proof}
	By d\'evissage, it suffices to show that for each $n\in \bfN$, 
	\begin{align*}
		\gr^F_n\bfA^{\omega} = \bigoplus_{\lambda\in W_S \lambda_0}\bigoplus_{\substack{w\in W_S \\ \ell(w) = n}} \bfC[V]\tau^{\omega}_w\bfe(\lambda).
	\end{align*}
	It follows from the braid relations for $\tau^{\omega}_a$ in $\Gr^F\bfA^{\omega}$ proven in~\autoref{lemm:braid} and the fact that $(\tau^{\omega}_a)^2\bfe(\lambda)\in F_{\le 1}\bfA^{\omega}$, that these elements $\tau^{\omega}_w$ span $\Gr^F_n\bfA^{\omega}$. By the invariant theory of reflection groups, the family $\{ \tau^{o}_w\bfe(\lambda) \}_{w}$ is free over $\bfC[V]$ and forms a basis for $\End_{\calZ}\left( \Pol_{\lambda} \right)$. In view of \autorefitem{prop:formula}{ii}, the matrix of transition between the families $\{ \tau^{o}_w\bfe(\lambda) \}_{w}$ and $\{ \tau^{\omega}_w\bfe(\lambda) \}_{w}$ is diagonal with non-zero entries; therefore the latter is also free over $\bfC[V]$.
\end{proof}

Define the filtration $F_{\le n}\bfA^{-\infty} = (F_{\le n}\bfA^{-o})\otimes_{\calZ}{\Frac\calZ}\subset \bfA^{-\infty}$.

\begin{coro}\label{coro:filt}
	For each $n$, we have
	\begin{align*}
		F_{\le n}\bfA^{\omega} = F_{\le n}\bfA^{-\infty} \cap \bfA^{\omega}.
	\end{align*}
\end{coro}
\begin{proof}
	Let $F'_{\le n}\bfA^{\omega} =  F_{\le n}\bfA^{-\infty} \cap \bfA^{\omega}$. We have $F_{\le n}\bfA^{\omega} \subset F'_{\le n}\bfA^{\omega}$. Fix $\lambda, \lambda'\in W_S \lambda_0$ and denote $A = \bfe(\lambda')\bfA^{\omega}\bfe(\lambda)$. Put $N = \# \left\{ w\in W_S\;;\; w\lambda = \lambda' \right\}$, then we have $F_{\le n}A = A  = F_{\le n}A'$ for $n\ge N$ by~\autoref{theo:basis}. We prove by induction on $k\in [0, N]$ that $F_{\le N - k}A = F'_{\le N - k}A$. It is already clear for $k = 0$. Suppose $k \ge 1$. Then we have the obvious diagram:
	\[
		\begin{tikzcd}
			0 \arrow{r} & F_{\le N-k}A \arrow{r}\arrow{d}{\varphi} & F_{\le N-k+1}A \arrow{r}\arrow{d}{\psi} & \Gr^{F}_{N-k+1}A \arrow{r}\arrow{d}{\eta} & 0 \\
			0 \arrow{r} & F'_{\le N-k}A \arrow{r} & F'_{\le N-k+1}A \arrow{r}& \Gr^{F'}_{N-k+1}A \arrow{r} & 0.
		\end{tikzcd}
	\]
	The morphism $\psi$ is an isomorphism by the induction hypothesis and $\varphi$ is injective. By the snake lemma, we have $\ker \eta \cong \coker \varphi$. \autoref{theo:basis} implies that $\Gr^{F}_{N-k + 1}A$ is $\bfC[V]$-torsion-free whereas $\coker\varphi$ is a $\bfC[V]$-torsion module. Therefore $\coker \varphi = 0$ and $\varphi$ is an isomorphism. Summing over $\lambda,\lambda'\in W_S\lambda_0$, we obtain $F_{\le n}\bfA^{\omega} = F'_{\le n}\bfA^{\omega}$ for all $n\in \bfN$.
\end{proof}
\begin{rema}
	In view of (the proof of)~\autoref{lemm:braid}, one can define a ``Bruhat filtration'' $\{F_{\calI}\}_{\calI}$ indexed by the order ideals $\calI$ of the affine Weyl group $W_S$ with respect to the Bruhat order, so that $F_{\calI}\bfA^{\omega}$ is spanned by $\bfC[V]\tau^{\omega}_w \bfe(\lambda)$ for $\lambda\in W_S \lambda_0$ and $w\in \calI$. Our filtration by length $\{F_{\le n} \bfA^{\omega}\}_{n\in \bfN}$ can be viewed as part of the Bruhat filtration because we have $F_n \bfA^{\omega} = F_{\calI_n}\bfA^{\omega}$ for $\calI_n = \left\{w\in W_S\;;\; \ell(w)\le n  \right\}$.
\end{rema}

\subsection{The associated graded \texorpdfstring{$\Gr^F\bfA^{\omega}$}{Gr A}}\label{subsec:grA}
We describe in greater detail the structure of the associated graded $\Gr^F\bfA^\omega$. We establish in~\autoref{theo:PBW} a triangular decomposition for $\Gr^F \bfA^{\omega}$, which will be used in the proof of~\autoref{prop:indtf}. The proof of~\autoref{lemm:CM} is technical. The reader is advised to skip this subsection in the first reading. \par

 Recall the extended affine Weyl group $\til W_S = P^{\vee}\rtimes W_R$ defined in~\autoref{subsec:sysaff-fini}. For $\mu\in P^{\vee}$, let $w_{\mu}\in W_R$ be such that $X^{\mu}w_{\mu}$ is the minimal element of the coset $X^\mu W_R$. Define the following map of minimal representatives:
\[
	\theta: P^{\vee}\to \til W_S, \quad \theta(\mu) = X^{\mu}w_{\mu}.
\]
In particular, $\theta(Q^{\vee})\subset W_S$ coincide with the set of minimal representatives for the quotient $W_S / W_R$.\par
\begin{lemm}\label{prop:minimal}
	For each $\mu\in P^{\vee}$, the element $w_\mu$ is characterised by the following property: every positive root $\alpha\in R^+$ satisfies  $w_{\mu}\alpha\in R^-$ if and only if $\langle \alpha, \mu\rangle > 0$. 
\end{lemm}
\begin{proof}
	See~\cite[Proof of 1.4]{cherednik95}
\end{proof}

We consider the nil-Hecke algebra $\bfC[\til W_S]^\nil$ for $\til W_S$: it is the $\bfC$-vector space span by the basis $\left\{ [w]^\nil \right\}_{w\in \til W_S}$ equipped with the following multiplication law
\[
	[w]^{\nil}\cdot[y]^{\nil} = \begin{cases} [wy]^{\nil} & \text{if $\ell(wy) = \ell(w) + \ell(y)$} \\ 0 & \text{otherwise.} \end{cases}
\]
Let $\bfC[W_S]^\nil$ and $\bfC[W_R]^\nil$ be the subspace of $\bfC[\til W_S]^\nil$ spanned by $\left\{ [w]^\nil \right\}_{w\in W_S}$ and $\left\{ [w]^\nil \right\}_{w\in W_R}$ respectively. These are the nil-Hecke algebras for $W_S$ and $W_R$. \par

Let $\rmC_0\subset V^*$ denote the fundamental Weyl chamber and $\ba\rmC_0$ its closure in $V^*$. Let $P^{\vee}_{\dom} = P^{\vee} \cap \ba{\rmC_0}$\index[ch2]{P@$P^{\vee}_+$}\index[ch2]{Q@$Q^{\vee}_+$} (resp. $Q^{\vee}_{\dom} = Q^{\vee} \cap \ba{\rmC_0}$) be the submonoid of $P^{\vee}$ consisting of dominant coweights (resp. dominant coroots). Let $\bfC P^{\vee}_{\dom}$ (resp. $\bfC Q^{\vee}_{\dom}$) denote the monoid algebra of $P^{\vee}_{\dom}$ (resp. $Q^{\vee}_{\dom}$). For $\mu\in P^{\vee}_+$, let $X^{\mu}\in \bfC P^{\vee}_+$ denote the corresponding element.
\par
We define a map
\begin{align*}
		\zeta: \bfC P^{\vee}_{\dom} \to \bfC[\til W_S]^\nil,\quad \zeta(X^{\mu}) = \sum_{\mu'\in W_R \mu } [X^{\mu'}]^\nil.
\end{align*}
\begin{lemm}\label{lemm:CM}
	The following statements hold:
	\begin{enumerate}
		\item\label{lemm:CM-i}
			The map $\zeta$ is a ring homomorphism and yield a left $\bfC P^{\vee}_{\dom}$-module structure on $\bfC[\til W_S]^\nil$ by left multiplication; moreover, $\bfC[\til W_S]^\nil$ is a free $(\bfC P^{\vee}_{\dom}, \bfC[W_R]^\nil)$-bimodule of rank $\# W$ and a basis of which is given by $\left\{ [\theta(b_w)]^\nil \right\}_{w\in W_R}$ with
				\begin{align*}
						b_w = \sum_{\substack{\alpha\in \Delta_0 \\ s_{\alpha}w < w}} w^{-1}w_0\omega^{\vee}_\alpha\in P^{\vee}.
				\end{align*}
		\item\label{lemm:CM-ii}
			The ring $\bfC Q^{\vee}_{\dom}$ is Cohen--Macaulay and the $\bfC P^{\vee}_{\dom}$-module structure on $\bfC[\til W_S]^{\nil}$ restricts to a $\bfC Q^{\vee}_{\dom}$-module structure on $\bfC[W_S]^{\nil}$; moreover, there is a decomposition
			\[
				\bfC [W_S]^\nil = \calE \otimes \bfC [W_R]^\nil,
			\]
			where $\calE\subset\bfC [W_S]^\nil$ is a $\bfC Q^{\vee}_{\dom}$-direct factor and is a Cohen-Macaulay $\bfC Q^{\vee}_{\dom}$-module of maximal dimension.
	\end{enumerate}
\end{lemm}
\begin{proof}
	In view of the length formula~\autoref{prop:lg} for $\til W_S$, the condition $\ell(X^{\mu+\nu}) = \ell(X^{\mu})+\ell(X^{\nu})$ is equivalent to that $\mu'$ and $\nu$ lie in the closure of the same Weyl chamber. Therefore the map $\zeta$ is a ring homomorphism. 		
	Define a decreasing filtration $G^{\bullet}\til W_S$ by 
	\[
		G^{k}\til W_S = \bigcup_{\substack{y\in W_R \\\ell(y) \ge k}} \left\{ v\in \til W_S\;;\; \ell(vy^{-1}) = \ell(v) - \ell(y) \right\}.
	\]
	Since $\bfC[\til W_S]^\nil$ has a canonical basis $\left\{ [w]^{\nil} \right\}_{w\in \til W_S}$, the filtration $G^{\bullet}\til W_S$ induces a filtration on $\bfC[\til W_S]^\nil$, denoted by $G^{\bullet}\bfC[\til W_S]^\nil$.
	\begin{enumerate}
		\item[Step 1.]
			We prove that the map
			\[
				P^{\vee}_{\dom}\times \left\{(w, y)\in (W_R)^2\;;\; \ell(y) = k  \right\}\xrightarrow{\sim} G^k \til W_S \setminus G^{k-1}\til W_S,\quad (\mu, w, y)\mapsto X^{b_w + w^{-1}w_0 \mu}wy
			\]
			is a bijection.

			For $\mu\in P^{\vee}$, let $w_{\mu}\in W_R$ be the element from~\autoref{prop:minimal}. We may partition $P^{\vee}$ into sub-semigroups :
			\[
				P^{\vee} = \bigsqcup_{w\in W_R} P^{\vee}_w,\quad P^{\vee}_w = \left\{ \mu\in P^{\vee}\;;\; w_{\mu} = w \right\}.
			\]
			For $w\in W_R$, we have $b_w\in P^{\vee}_w$ and there is a bijection
			\[
				P^{\vee}_{\dom}\xrightarrow{\sim} P^{\vee}_w,\quad \mu\mapsto b_w + w^{-1}w_0 \mu;
			\]
			we can thus express the set $G^k \til W_S \setminus G^{k-1}\til W_S$ as
			\[
				G^k \til W_S \setminus G^{k-1}\til W_S = \bigsqcup_{\substack{w,y\in W_R \\ \ell(y) = k}} \left\{X^{\mu}w y\;;\; \mu\in P^{\vee}_w  \right\} =  \bigsqcup_{\substack{w,y\in W_R \\ \ell(y) = k}} \left\{X^{b_w + w^{-1}w_0\mu}w y\;;\; \mu\in P^{\vee}_{\dom}  \right\}.
			\]
		\item[Step 2.]
			We prove that for each $\mu\in P^{\vee}$ and $w,y\in W_R$, we have $\zeta(X^{\mu})[X^{b_w}wy]^\nil\in G^{\ell(y)} \bfC[\til W_S]^\nil$ and
			\begin{align}\label{equa:zetaX}
				\zeta(X^{\mu})[X^{b_w}wy]^\nil \equiv [X^{b_w + w^{-1}w_0 \mu}wy]^\nil \mod G^{\ell(y)+1} \bfC [\til W_S]^\nil.
			\end{align}
			Indeed, the defining relations of the nil-Hecke algebra $\bfC[\til W_S]^{\nil}$ yield
			\[
				\zeta(X^{\mu})[X^{b_w}wy]^\nil = \sum_{\substack{\mu'\in W_R \mu \\ \ell(X^{b_w + \mu'}wy) = \ell(X^{b_w}wy) + \ell(X^{\mu'})}} [X^{b_w + \mu'}wy]^\nil
			\]
			in $\bfC [\til W_S]^\nil$. Since
			\[
				\ell(X^{b_w + \mu'}wy) \le \ell(X^{b_w + \mu'}w) + \ell(y) \le  \ell(X^{b_w}w) + \ell(X^{\mu'}) + \ell(y) = \ell(X^{b_w}wy) + \ell(X^{\mu'})
			\]
			(the last equality due to~\autoref{prop:minimal}), the condition
			\begin{align}\label{equa:lsum}
				\ell(X^{b_w + \mu'}wy) = \ell(X^{b_w}wy) + \ell(X^{\mu'}) 
			\end{align}
			implies that $\ell(X^{b_w + \mu'}wy) = \ell(X^{b_w + \mu'}w) + \ell(y)$ and hence $X^{b_w + \mu'}wy\in G^{\ell(y)} \til W_S$. It follows that $[X^{b_w + \mu'}wy]^\nil \in G^{\ell(y)} \bfC [\til W_S]^\nil$ for $\mu'\in W_R \mu$ satisfying~\eqref{equa:lsum} and $[X^{b_w + \mu'}wy]^\nil \in G^{\ell(w)+1} \bfC [\til W_S]^\nil$ unless $w_{b_w + \mu'} =  w$; the latter case happens for the unique element $\mu' = w^{-1}w_0\mu$ in the orbit $W_R \mu$; therefore~\eqref{equa:zetaX} holds.
		\item[Step 3.]
			By Step 1 and Step 2, we see that $G^{k}\bfC[\til W_S]$ is a $\bfC P^{\vee}_{\dom}$-submodules and the successive quotient $G^k\bfC[\til W_S]^\nil / G^{k+1}\bfC[\til W_S]^\nil$ is a free $\bfC P^{\vee}_{\dom}$-module with a basis formed by the congruence classes of $\left\{ [X^{b_w}wy]^\nil \right\}_{y,w\in W_R,\, \ell(y) = k}$. It follows that $\left\{ [X^{b_w}wy]^\nil \right\}_{y,w\in W_R}$ forms a $\bfC P^{\vee}_{\dom}$-basis for $\bfC[\til W_S]^{\nil}$. Since $X^{b_w}w$ is minimal in the coset $X^{b_w}W_R$, we have 
			\[
				[X^{b_w}wy]^\nil = [X^{b_w}w]^\nil\cdot [y]^\nil,\quad \text{for $y\in W_R$};
			\]
			thus $\left\{ [X^{b_w}w]^\nil \right\}_{w\in W_R}$ forms a $(\bfC P^{\vee}_{\dom}, \bfC [W_R]^\nil)$-bimodule basis for $\bfC[\til W_S]^\nil$, whence~\ref{lemm:CM-i}. \par
			\item[Step 4.]
				Let $\calE'\subset \bfC[\til W_S]^\nil$ be the (free) $\bfC P^{\vee}_{\dom}$-submodule generated by $\left\{ [\theta(b_w)]^\nil \right\}_{w\in W_R}$ so that, by~\ref{lemm:CM-i}, there is a decomposition $\bfC[\til W_S]^\nil \cong\calE' \otimes \bfC[W_R]^\nil$.
				Let $\Omega = P / Q$. Define a $\bfC$-linear action of $\Omega$ on $\bfC [\til W_S]^\nil$ by
				\begin{align*}
						\Omega \times \bfC[\til W_S]^\nil&\to \bfC[\til W_S]^\nil \\
						(\beta, [X^{\mu}w]^\nil) &\mapsto e^{2\pi i\langle \beta,\mu\rangle }[X^{\mu}w]^\nil,\quad \beta\in \Omega,\; \mu\in P^{\vee},\; w\in W_R.
				\end{align*} 
				This action preserves the subspace $\calE'\subset \bfC[\til W_S]^\nil$ and fixes $\bfC[W_R]^\nil$ pointwise; hence there is a decomposition of the $\Omega$-fixed subspace 
				\[
					\bfC[W_S]^\nil = (\bfC[\til W_S]^\nil )^{\Omega} \cong \calE \otimes \bfC[W_R]^\nil,\quad \text{where }\calE = (\calE')^{\Omega}. 
				\]
				It remains to show that $\calE$ is a Cohen--Macaulay $\bfC Q^{\vee}_{\dom}$-module of maximal dimension. \par
				Since $\bfC Q^{\vee}_{\dom}$ is integrally closed and $\bfC P^{\vee}_{\dom}$ is regular and an integral ring extension of it, by~\cite[X.2.6,coro 2]{bourbakiAC}, $\bfC P^{\vee}_{\dom}$ is a Cohen--Macaulay $\bfC Q^{\vee}_{\dom}$-module. Thus $\calE'$, being a free $\bfC P^{\vee}_{\dom}$-module, is Cohen--Macaulay of maximal dimension over $\bfC Q^{\vee}_{\dom}$. Since $\calE$ is a direct factor of $\calE'$, so is it Cohen--Macaulay of maximal dimension over $\bfC Q^{\vee}_{\dom}$, whence~\ref{lemm:CM-ii}.
	\end{enumerate}
\end{proof}
Below, we will work with $\gr^F \bfA^{\omega}$ and view the elements $\tau^\omega_a\bfe(\lambda)$ as in $\gr^F \bfA^{\omega}$ for the sake of notational simplicity. View $\bfA^{\omega}$ as $(\bfA^{\omega})^{\op}$-module via the right regular representation. The ring $\End_{(\gr^F\bfA^{\omega})^{\op}}(\gr^F\bfA^{\omega})$ can be viewed as a unital completion of $\bfA^{\omega}$. Define a $\bfC$-linear map 
\begin{align}\label{equa:Theta}
	\Theta: \bfC[W_S]^{\nil}&\to \End_{(\gr^F\bfA^{\omega})^{\op}}(\gr^F\bfA^{\omega}) = \prod_{\lambda\in W_S\lambda_0}\gr^F\bfA^{\omega}\bfe(\lambda)\\
	\Theta([w]^{\nil}) &= \tau^{\omega}_w = \sum_{\lambda\in W_S\lambda_0}\tau^{\omega}_w \bfe(\lambda),\quad w\in W_S \nonumber
\end{align}

\begin{prop}\label{theo:PBW}
	There is a triangular decomposition 
	\begin{align*}
			\Gr^F\bfA^{\omega} \cong \calE\otimes_{\bfC} \left(\bigoplus_{w\in W_R}\bfC \tau^{\omega}_w\right)\otimes_{\bfC} \left(\bigoplus_{\lambda\in W_S \lambda_0}\bfC[V]\bfe(\lambda)\right),
	\end{align*}
	where $\calE\subset \bfC[W_S]^{\nil}$ is the $\bfC Q^{\vee}_{\dom}$-submodule from~\autorefitem{lemm:CM}{ii}.
\end{prop}
\begin{proof}
	From~\autoref{theo:basis}, we see that the $\bfC[W_S]^{\nil}$-action on $\gr^F \bfA^{\omega}$ via $\Theta$ yields a decomposition
	\[
		\bfC[W_S]^{\nil}\otimes \left(\bigoplus_{\lambda\in W_S \lambda_0}\bfC[V]\bfe(\lambda)\right) \xrightarrow{\sim} \gr^F\bfA^{\omega},\quad f\otimes b \mapsto \Theta(f)(b).
	\]
	By~\autorefitem{lemm:CM}{ii}, we can further decompose $\bfC[W_S]^{\nil} = \calE\otimes \bfC[W_R]^{\nil}$. Finally, we have $\Theta(\bfC[W_R]^{\nil}) = \bigoplus_{w\in W_R}\bfC \tau^{\omega}_w$.
\end{proof}

\section{Module categories of \texorpdfstring{$\bfA^{\omega}$}{A omega}}\label{sec:Agmod}
We keep the notation of~\autoref{sec:agha}.  We put a $\bfZ$-grading on $\bfA^{\omega}$ as follows: the generators are homogeneous: $\deg\alpha\,\bfe(\lambda) = 2$ for $ \alpha\in V^*$ and $\deg \tau^{\omega}_a \bfe(\lambda) = \omega_{\lambda}(a) + \omega_{s_a\lambda}(a)$. If $M = \bigoplus_n M_n$ is a graded vector space, denote by $M\langle m\rangle$\index[ch2]{@$\langle n\rangle$} the grading shift given by $M\langle m\rangle_n = M_{m+n}$. For two graded vector spaces $M$ and $N$, we denote by $\Hom(M, N)$ the space of $\bfC$-linear maps of degree $0$ and $\gHom(M, N) = \bigoplus_{k\in\bfZ}\Hom(M, N\langle k\rangle)$\index[ch2]{gHom@$\gHom$}. \par

Below, by ``modules'' we mean left modules. All statements can be turned into those for right modules by means of the anti-involution $\bfA^{\omega} \cong \left(\bfA^{\omega}\right)^{\op}$ defined by $\tau^{\omega}_a\bfe(\lambda) \mapsto \tau^{\omega}_a\bfe(s_a\lambda)$.
\subsection{Graded \texorpdfstring{$\bfA^{\omega}$}{A omega}-modules}\label{subsec:grAmod}
An $\bfA^\omega$-module $M$ is called a weight module if there is a decomposition
\begin{align*}
	M = \bigoplus_{\lambda\in W_S \lambda_0} \bfe(\lambda) M.
\end{align*}
Let $\bfA^\omega\gMod$\index[ch2]{gmod@$\gmod,\gmod_0$} denote the category of graded weight modules of $\bfA^\omega$. Let $\bfA^\omega\gmod\subset \bfA^\omega\gMod$ be the subcategory of compact objects (i.e. $M\in \bfA^{\omega}\gmod$ if $\Hom_{\bfA^{\omega}\gMod}(M, \relbar)$ commutes with filtered colimits) and let $\bfA^{\omega}\gmod_0\subset \bfA^{\omega}\gmod$ be the subcategory of $\frakm_{\calZ}$-nilpotent objects. The following lemma is obvious.
\begin{lemm}\label{prop:gmod}
	For every object $M\in \bfA^{\omega}\gMod$ there exists an index set $J$ and two families of integers $\left\{ a_j \right\}_{j\in J}$ and $\left\{ \lambda_j \right\}_{j\in J}$ such that there exists an epimorphism in $\bfA^{\omega}\gMod$
	\[
		\bigoplus_{j = 1}^r\bfA^\omega\bfe(\lambda_j)\langle a_j\rangle\twoheadrightarrow M. \pushQED{\qed}\qedhere\popQED
	\]
\end{lemm}
We define a homomorphism of graded rings
	\begin{align}\label{equa:centreH}
		\calZ \to \gEnd\left( \id_{\bfA^\omega\gMod} \right)
	\end{align}
	as follows: For every $f\in \bfC[V]^{W_{\lambda_0}}$ and $w\in W_S$, let $f$ acts on $\bfe(w\lambda_0)M$ by multiplication with $(\partial w)(f)\in \bfC[V]^{W_{w\lambda_0}}$. 

\subsection{Intertwiners}
For each $\lambda\in W_S \lambda_0$ and $a\in \Delta$, introduce the following element in $\bfA^\omega$:
\begin{align*}
	\varphi_a \bfe(\lambda) = \begin{cases} ((\partial a) \tau^{\omega}_a + 1 ) \bfe(\lambda) & \omega_{\lambda}(a) = -1 \\ \tau^{\omega}_a \bfe(\lambda) & \omega_{\lambda}(a) \ge 0 \end{cases}.
\end{align*}
It satisfies the following relations:
\begin{align*}
	\varphi_a^2\bfe(\lambda) &=  \begin{cases} \bfe(\lambda) &\omega_{\lambda}(a) = -1\\ \pm (\partial a)^{n_{\lambda, a}}\bfe(\lambda) & \omega_{\lambda}(a) \ge 0 \end{cases}\\
	\varphi_a f \bfe(\lambda)&= s_a(f) \varphi_a\bfe(\lambda)\quad f\in \bfC[V],
\end{align*}
where $n_{\lambda,a} = \max(\omega_{\lambda}(a)+\omega_{s_a\lambda}(-a), 0)$. These elements satisfy the usual braid relations. Thus, we may write $\varphi_w\bfe(\lambda) = \varphi_{a_l}\cdots\varphi_{a_1}\bfe(\lambda)$ by choosing any reduced expression $w = s_{a_l} \cdots s_{a_1}$.
\begin{lemm}\label{lemm:phi-omega}
		Let $w\in W_S$ and $a\in \Delta$. Then the right multiplication by the intertwiner $\varphi_a$ induces an isomorphism of $\bfA^\omega$-modules
		\begin{align*}
			\bfA^\omega\bfe(\lambda) \cong\bfA^\omega\bfe(s_a\lambda)
		\end{align*}
		if $\omega_{\lambda}(a) + \omega_{s_a\lambda}(-a)\le 0$.
	\end{lemm}
	\begin{proof}
		The right multiplication by the element $\varphi_a\bfe(s_a\lambda) = \bfe(\lambda)\varphi_a\bfe(s_a\lambda)$ yields $\bfA^{\omega}\bfe(\lambda)\xrightarrow{\sim} \bfA^{\omega}\varphi_a\bfe(s_a\lambda)\xrightarrow{\sim}  \bfA^{\omega}\varphi_a^2\bfe(s_a\lambda)$. Hence if $\varphi^2_a\bfe(\lambda) = f\bfe(\lambda)\in \bfC[V]\bfe(\lambda)$ for $f\in \bfC[V]$ invertible, then $\varphi^2_a\bfe(\lambda)$ is an isomorphism. The condition that $f$ be invertible is exactly as stated. Clearly, if $\varphi^2_a\bfe(\lambda)$ and $\varphi^2_a\bfe(s_a\lambda)$ are isomorphisms, then so are $\varphi_a\bfe(\lambda)$ and $\varphi_a\bfe(s_a\lambda)$. The statement follows. 
	\end{proof}

	\subsection{Clan decomposition}\label{subsec:clans}
	As in~\autoref{subsec:basis}, we extend $\omega_{\lambda_0}$ to a $W_S$-invariant function $\til\omega_{\lambda_0}:S\to \bfZ_{\ge -1}$ and we suppose that the extension $\til\omega_{\lambda_0}$ has finite support. Consider the following sub-family of hyperplanes
	\begin{align*}
			\frakD^{\omega} = \left\{ H_a\subset E\;;\; a\in S,\; \til\omega_{\lambda_0}(a) \ge 1 \right\}.
	\end{align*}
	The connected components of the following space
	\begin{align*}
		E^{\omega}_{\circ} = E \setminus \bigcup_{H\in \frakD^{\omega}} H
	\end{align*}\index[ch2]{E@$E^{\omega}_{\circ}$} 
	are called clans. Since $\til\omega_{\lambda_0}$ is supposed to be finitely supported, the family $\frakD^{\omega}$ is finite, the set of connected components $\pi_0\left(E^{\omega}_{\circ}\right)$ is finite and there are only a finite number of clans.  \par

	Let $\frakC\subset E^{\omega}_{\circ}$ be a clan. Since $E^{\omega}_{\circ}$ is the complement of a finite hyperplane arrangement, $\frakC$ is a convex polytope. The \textbf{salient cone} of $\frakC$ is defined to be the convex polyhedral cone $\kappa\subset V$ whose dual cone $\kappa^{\vee}$ is the cone of linear functions which are bounded from below on $\frakC$:
	\begin{align*}
		\kappa^{\vee} = \left\{ v\in V^*\;;\; \inf_{x\in \frakC} \langle v,x\rangle > -\infty \right\},\quad \kappa = \kappa^{\vee\vee} = \left\{ x\in V\;;\; \langle v,x\rangle \ge 0,\; \forall v\in \kappa^{\vee} \right\}.
	\end{align*}
	Then $\kappa$ is a convex polyhedral generated by a finite subset of $P^{\vee}$. We say that clan $\frakC\subset E^{\omega}_{\circ}$ is \textbf{generic} if its salient cone is of maximal dimension. \par

	Denote by $\nu_0\in E$ the fundamental alcove associated with the basis $\Delta$.
	\begin{lemm}\label{lemm:clanisom}
		Let $w\in W_S$ and $a\in \Delta$. Then $w^{-1} \nu_0$ and $w^{-1}s_a\nu_0$ are in the same clan if and only if the intertwiner $\varphi_a$ induces an isomorphism of $\bfA^\omega$-modules
		\begin{align*}
			\bfA^\omega\bfe(w \lambda_0) \cong\bfA^\omega\bfe(s_aw \lambda_0).
		\end{align*}
	\end{lemm}
	\begin{proof}
		Using~\autoref{lemm:phi-omega}, we have
		\begin{align*}
			\varphi_a^2\bfe(w\lambda_0) = \bfe(w\lambda_0) \Leftrightarrow \omega_{w\lambda_0}(a) + \omega_{s_aw\lambda_0}(-a) \le 0 \\
			\Leftrightarrow \til\omega_{\lambda_0}(w^{-1}a) + \til\omega_{\lambda_0}(-w^{-1}a) \le 0 \Leftrightarrow H_{wa}\notin \frakD^{\omega}
		\end{align*}
		The last condition is equivalent to that $w^{-1}\nu_0$ and $w^{-1}s_a\nu_0$ belong to the same clan.
	\end{proof}
	The following proposition follows immediately from the above lemma.
	\begin{prop}\label{prop:entrelacement}
		If $w,w'\in W_S$ are such that $w^{-1}\nu_0$ and $w'^{-1}\nu_0$ lie in the same clan, then right multiplication by the intertwiner $\varphi_{w'w^{-1}}\bfe(w\lambda)$ yields an isomorphism $\bfA^{\omega}\bfe(w'\lambda)\to \bfA^{\omega}\bfe(w\lambda)$.
	\end{prop}
	\qed
	\begin{coro}\label{coro:gdim}
		Let $M\in \bfA^{\omega}\gmod$. If $w,w'\in W_S$ are such that $w\nu_0^{-1}$ and $w'\nu_0^{-1}$ lie in the same clan, then multiplication by the intertwiner $\varphi_{w'w^{-1}}\bfe(w\lambda)$ yields an isomorphism of graded $\calZ$-modules $\bfe(w\lambda_0)M\cong \bfe(w'\lambda_0)M$. In particular, in this case there is an equality of graded dimensions 
		\begin{align*}
			\gdim \bfe(w\lambda_0)M = \gdim \bfe(w'\lambda_0)M.
		\end{align*}
	\end{coro}
	\begin{proof}
		Indeed, we have
		\begin{align*}
			\bfe(w\lambda_0) M\cong \Hom_{\bfA^{\omega}}\left( \bfA^{\omega}\bfe(w\lambda_0), M\right)\xrightarrow{\relbar\circ\varphi_{w'w^{-1}}}\Hom_{\bfA^{\omega}}\left( \bfA^{\omega}\bfe(w'\lambda_0), M\right)\cong \bfe(w\lambda_0) M.
		\end{align*}
	\end{proof}

	\begin{exam}\label{exam:clan}
	In the setting of~\autorefitem{exam:Aomega}{iii}, the alcoves in $E$ are of the form $\left] n, n + 1/2\right[$ for $n\in (1/2)\bfZ$ and the fundamental alcove is $\nu_0 = \left]0, 1/2\right[$. We have $\frakD^{\omega} = \left\{ H_{a_0}, H_{a_1} \right\}$, with $\left\{ a_0 = 1 - 2\epsilon, a_1 = 2\epsilon \right\} = \Delta$. The clan decomposition is depicted as follows:
		\usetikzlibrary{plotmarks}
		\[
			\begin{tikzpicture}
				\draw (-5,0) -- (-3,0) node[below]{$\frakC_-$} -- (-1,0) node[above]{$H_{a_1}$} node[below]{$0$} -- (0,0) node[below]{$\frakC_0$} -- (1,0) node[above]{$H_{a_0}$} node[below]{$1/2$} -- (3,0) node[below]{$\frakC_+$} -- (5,0);
				\pgfplothandlermark{\pgfuseplotmark{*}}
				\pgfplotstreamstart
				\pgfplotstreampoint{\pgfpoint{-1cm}{0cm}}
				\pgfplotstreampoint{\pgfpoint{1cm}{0cm}}
				\pgfplotstreamend
			\end{tikzpicture}
		\]
	The clans $\frakC_{-} = {\left]-\infty, 0\right[}$ and $\frakC_{+} = {\left]1/2, +\infty\right[}$ are generic whereas the clan $\frakC_0 = {\left] 0, 1/2\right[} = \nu_0$ is not generic. To each alcove $\nu = w^{-1}\nu_0$ with $w\in W_S$, we attach the element $\lambda_{\nu} = w \lambda_0\in E$
		\[
			\begin{tikzpicture}
				\draw (-5,0) -- (-4,0) node[above]{$5/4$} -- (-3,0) node[below]{$-1/2$} -- (-2,0) node[above]{$-1/4$} -- (-1,0) node[below]{$0$} -- (0,0) node[above]{$1/4$} -- (1,0) node[below]{$1/2$} -- (2,0) node[above]{$-3/4$} -- (3,0) node[below]{$3/2$} -- (4,0)  node[above]{$5/4$} -- (5,0) (5.5,0) node[above]{$\lambda_{\nu}$} node[below]{$\epsilon$} ;
				\pgfplothandlermark{\pgfuseplotmark{|}}
				\pgfplotstreamstart
				\pgfplotstreampoint{\pgfpoint{-3cm}{0cm}}
				\pgfplotstreamend
				\pgfplothandlermark{\pgfuseplotmark{*}}
				\pgfplotstreamstart
				\pgfplotstreampoint{\pgfpoint{-1cm}{0cm}}
				\pgfplotstreampoint{\pgfpoint{1cm}{0cm}}
				\pgfplotstreamend
				\pgfplothandlermark{\pgfuseplotmark{|}}
				\pgfplotstreamstart
				\pgfplotstreampoint{\pgfpoint{3cm}{0cm}}
				\pgfplotstreamend
			\end{tikzpicture}
		\]
	In particular, the alcoves $\nu = \left]1/2,3/2\right[$ and $\nu' = \left]3/2,5/2\right[$ lie in the same clan $\frakC_+$ with $\lambda_{\nu} = -3/4$ and $\lambda_{\nu'} = s_0\lambda_{\nu} = 5/4$. In this case~\autoref{prop:entrelacement} amount to the fact that the intertwiners $\varphi_{a_0}\bfe(\lambda_{\nu'}) : \bfA^{\omega}\bfe(\lambda_{\nu})\to \bfA^{\omega}\bfe(\lambda_{\nu'})$ and $\varphi_{a_0}\bfe(\lambda_{\nu}) : \bfA^{\omega}\bfe(\lambda_{\nu'})\to \bfA^{\omega}\bfe(\lambda_{\nu})$ are isomorphisms and inverse to each other.\par
		The projective $\bfA^\omega$-modules $\bfA^{\omega}\bfe(\lambda_{\nu})$ are indecomposable and they are non-isomorphic for alcoves $\nu$ in the three different clans $\frakC_-, \frakC_0$ and $\frakC_+$. Choose any alcoves $\nu_+\subset \frakC_+$, $\nu_-\subset \frakC_-$ and denote $\lambda_+ = \lambda_{\nu_+}, \lambda_- = \lambda_{\nu_-}$, $P_+ = \bfA^{\omega}\bfe(\lambda_{+})$,$P_0 = \bfA^{\omega}\bfe(\lambda_0)$ and $P_- = \bfA^{\omega}\bfe(\lambda_{-})$.  Their simple quotients, denoted by $L_+$, $L_0$ and $L_-$, form a complete collection of simple objects of $\bfA^{\omega}\gmod$ up to grading shifts. The graded dimension is given by
		\begin{align*}
			\gdim \bfe(\lambda_{\nu}) L_* = \begin{cases} 1 & \nu\subseteq \frakC_* \\ 0 & \nu\not\subseteq \frakC_*\end{cases},\quad *\in \{+, 0, -\}.
		\end{align*}
		In particular, $L_+$ and $L_-$ are infinite-dimensional and $L_0$ is finite-dimensional.
		The cosocle filtrations of $P_+$, $P_0$ and $P_-$ are described as follows: 
		\begin{align*}
			P_+ = \begin{bNiceMatrix}[small]L_+ & & \\ &L_0\langle -1\rangle& \\ L_+\langle -2\rangle& & L_-\langle -2\rangle \\ &L_0\langle -3\rangle& \\ L_+\langle -4\rangle& & L_-\langle -4\rangle \\ &\vdots& \end{bNiceMatrix},\;
			P_0 = \begin{bNiceMatrix}[small] & L_0 & \\ L_+\langle -1\rangle& & L_-\langle -1\rangle \\ &L_0\langle -2\rangle& \\ L_+\langle -3\rangle& & L_-\langle -3\rangle \\&L_0\langle -2\rangle& \\ &\vdots& \end{bNiceMatrix},\;
			P_- = \begin{bNiceMatrix}[small] & & L_- \\ &L_0\langle -1\rangle& \\ L_+\langle -2\rangle& & L_-\langle -2\rangle \\ &L_0\langle -3\rangle& \\ L_+\langle -4\rangle& & L_-\langle -4\rangle \\ &\vdots& \end{bNiceMatrix}.
		\end{align*}
	\end{exam}

	\subsection{Basic properties of graded modules of \texorpdfstring{$\bfA^{\omega}$}{A omega}}\label{subsec:Agmod}
	We choose a finite subset $\Sigma\subset W_S$ such that for every clan $\frakC\subset E^{\omega}_{\circ}$, there exists $w\in \Sigma$ with $w^{-1}\nu_0\subset \frakC$. Set $\bfe_{\Sigma} = \sum_{w\in \Sigma}\bfe(w\lambda_0)$ and $P_{\Sigma} = \bfA^\omega\bfe_\Sigma$. 
\begin{lemm}\label{lemm:PSigma}
	The module $P_{\Sigma}$ is a graded compact projective generator of $\bfA^\omega\gMod$. 
\end{lemm}
\begin{proof}
	For any $y\in W_S$, we can find $w\in \Sigma$ such that $y^{-1}\nu_0$ and $w^{-1}\nu_0$ are in the same clan. By~\autoref{prop:entrelacement}, there exists an isomorphism
	\begin{align*}
		\bfA^\omega\bfe(w\lambda_0)\cong \bfA^\omega\bfe(y\lambda_0)
	\end{align*}
	Since the former is a direct factor of $P_{\Sigma}$, the above isomorphism yields a surjection $P_{\Sigma}\twoheadrightarrow \bfA^\omega\bfe(y\lambda_0)$. Combining this with~\autoref{prop:gmod}, we see that $P_{\Sigma}$ is a graded generator, which is clearly compact projective. 
\end{proof}

Put $A_{\Sigma} = (\gEnd_{\bfA^\omega\gMod} P_{\Sigma})^{\op} = \bfe_\Sigma \bfA^{\omega}\bfe_\Sigma$. It follows from~\autoref{lemm:PSigma} and the Morita theory that there is a graded equivalence 
\begin{align}\label{equa:morita1}
		\gHom_{\bfA^\omega\gMod}(P_{\Sigma}, \relbar):\bfA^\omega\gMod\xrightarrow{\sim} A_{\Sigma}\gMod,
\end{align}
which restricts to an equivalence on the subcategories of compact objects $\bfA^\omega\gmod\xrightarrow{\sim} A_{\Sigma}\gmod$.

\begin{prop}\label{prop:Agmod}
	The following statements hold:
	\begin{enumerate}
		\item\label{prop:gmod-i}
			The category $\bfA^\omega\gmod$ is noetherian and the subcategory $\bfA^\omega\gmod_0$ consists of objects of finite length.
		\item\label{prop:gmod-ii}
			For each $M\in \bfA^{\omega}\gmod$ and each $\lambda\in W_S \lambda_0$, the graded dimension $\gdim \bfe(\lambda)M$ is in $\bfN(\!(v)\!)$. Moreover, $M\in \bfA^{\omega}\gmod_{0}$ if and only if $\gdim \bfe(\lambda)M\in \bfN[v^{\pm 1}]$ for all $\lambda\in W_S \lambda_0$.
		\item\label{prop:gmod-iii}
			Every object of $\bfA^\omega\gmod$ admits a projective cover in the same category.
		\item\label{prop:gmod-iv}
			We have $\Irr(\bfA^\omega\gmod_0)\cong \Irr(\bfA^\omega\gmod)$. 
		\item\label{prop:gmod-v}
			The map~\eqref{equa:centreH} is an isomorphism $\calZ\cong \gEnd\left( \id_{\bfA^\omega\gmod} \right)$.
	\end{enumerate}
\end{prop}
\begin{proof}
	By the graded Morita equivalence~\eqref{equa:morita1}, it suffices to show the corresponding statements for $A_{\Sigma}\gmod$. \par
	Since $A_{\Sigma}$ is of finite rank over the graded polynomial ring $\bfC[V]^{W_{\lambda_0}}$, it is laurentian (i.e. its graded dimension is in $\bfN(\!(v)\!)$) and thus graded semi-perfect. The statements~\ref{prop:gmod-i}--\ref{prop:gmod-iv} result from the laurentian property. \par
			We prove \ref{prop:gmod-v}. Consider the $\bfA^\omega$-module $\Pol_{W_{S}\lambda_0}\in \bfA^\omega\gmod$. Since each factor $\Pol_{\lambda} = \bfC[V]$ is a free $\calZ$-module of finite rank, the sum $\Pol_{W_S\lambda_0}$ is a free $\calZ$-module of infinite rank. Taking base-change to the rational function field $\Frac \calZ$, we get a homomorphism
			\begin{align*}
				\rho:\bfA^{-\infty}\to \bigoplus_{\lambda,\lambda'\in W_S \lambda_0}\Hom_{\Frac\calZ}\left(\Rat_\lambda, \Rat_{\lambda'}  \right),
			\end{align*}
			We claim that $\rho$ is an isomorphism. It is injective since $\Pol_{W_S \lambda_0}$ is a faithful $\bfA^\omega$-module by definition and it remains faithful after localisation. It is easy to see from the definition of $\bfA^\omega$ that for $\lambda\in W_S \lambda_0$ and $a\in \Delta$, the operator $s_a\bfe(\lambda):\Rat_\lambda\to \Rat_{s_a\lambda}$ is in the image of $\rho$. For any $\lambda,\lambda'\in W_S \lambda_0$, let $W_{\lambda,\lambda'} = \left\{ w\in W_S \;;\; w\lambda = \lambda' \right\}$. The family $\{\bfe(\lambda') w \bfe(\lambda)\}_{w\in W_{\lambda,\lambda'}}$ is in the image of $\rho$. The rational function field $\Rat_{\lambda}$ is a Galois extension of $\Frac\calZ$ with Galois group $W_{\lambda}$. It follows from the Galois theory that
			\begin{align*}
				\End_{\Frac\calZ}(\Rat_\lambda)\cong \Rat_\lambda\rtimes \bfC W_{\lambda}.
			\end{align*}
			We have already seen that $\left\{ w\,\bfe(\lambda) \right\}_{w\in W_{\lambda}}$ is in $\image \rho$ the and $\Rat_\lambda$ is also in the image of $\rho$. It follows that $\End_{\Frac\calZ}(\Rat_\lambda)\subset \image \rho$. Let $\lambda,\lambda'\in W_S \lambda_0$ and choose $w\in W_{\lambda, \lambda'}$. Then $w\,\bfe(\lambda)\in \image \rho$ is an isomorphism $w\,\bfe(\lambda): \Rat_{\lambda} \cong \Rat_{\lambda'}$ and the pre-composition yields 
			\begin{align*}
				\relbar \circ w\bfe(\lambda): \End_{\Frac\calZ}(\Rat_\lambda)\cong \Hom_{\Frac\calZ}\left(\Rat_\lambda, \Rat_{\lambda'}  \right).
			\end{align*}
			Thus $\Hom_{\Frac\calZ}\left(\Rat_\lambda, \Rat_{\lambda'}  \right)\subset \image \rho$. We see that $\rho$ is surjective and the claim is proven. \par
			 There is an isomorphism
			\[
			A_{\Sigma}\otimes_{\calZ}\Frac\calZ = \bfe_{\Sigma}\bfA^{-\infty}\bfe_{\Sigma} \cong \End_{\Frac\calZ}\left(\bigoplus_{w\in \Sigma}\Rat_{w\lambda_0} \right)
			\]
			induced by $\rho$. Since the right-hand side is a matrix algebra over a field $\Frac \calZ$, its centre is $\Frac \calZ$. It follows that $\rmZ(A_{\Sigma}) = \Frac\calZ$. Hence
			\begin{align*}
					\gEnd\left( \id_{\bfA^\omega\gmod} \right)&\cong \gEnd\left( \id_{A_{\Sigma}\gmod} \right) = \rmZ\left( A_{\Sigma} \right) = \rmZ\left(A_{\Sigma}\otimes_{\calZ}\Frac\calZ \right)\cap A_{\Sigma}\\
					&= \Frac\calZ\cap A_{\Sigma} = \calZ,
			\end{align*}
			where the last equation follows from the basis theorem~\autoref{theo:basis}.
\end{proof}

\subsection{Basic properties of ungraded \texorpdfstring{$\bfA^{\omega}$}{A}-modules}
Let $U:\bfA^\omega\gmod_0\to \bfA^\omega\mof_0$ be the grading-forgetting functor. We extend it to $U:\bfA^\omega\gmod\to \Pro(\bfA^\omega\mof_0)$ by requiring $U$ to preserve filtered inverse limits. The extended functor is exact. Define the subcategory $\bfA^{\omega}\mof^{\wedge}\subset \Pro(\bfA^\omega\mof_0)$ to be the essential image of this functor. Let $\calZ^{\wedge} = \varprojlim_{N\rightarrow \infty} \calZ / \frakm_{\calZ}^N$. 
\begin{prop}
	Then the following properties are satisfied:
	\begin{enumerate}
		\item\label{prop:Amod-i}
			The functor forgetting the grading $U:\bfA^\omega\gmod\to \bfA^{\omega}\mof$ is exact and it induces $\Irr(\bfA^\omega\gmod) / \langle \bfZ\rangle \cong\Irr(\bfA^\omega\mof^{\wedge})$. Moreover, for all $M, N\in \bfA^\omega\gmod$ and $n\in \bfN$ we have 
			\begin{align*}
				\prod_{k\in \bfZ}\Ext^n(M, N\langle k\rangle) \cong \Ext^n(U M, U N).
			\end{align*}

		\item
			The category $\bfA^{\omega}\mof^{\wedge}$ is noetherian and the subcategory $\bfA^\omega\mof_0$ consists of objects of finite length.
		\item
			Every object of $\bfA^{\omega}\mof^{\wedge}$ admits a projective cover in the same category.
		\item
			We have $\Irr(\bfA^\omega\mof_0)\cong \Irr(\bfA^{\omega}\mof^{\wedge})$. 
		\item
			The ungraded analogue of the map~\eqref{equa:centreH} induces an isomorphism $\calZ^{\wedge}\cong \End\left( \id_{{\bfA^\omega}\mof^{\wedge}} \right)$.
	\end{enumerate}
\end{prop}
These statements follow from~\autoref{prop:Agmod}. \qed \par

\subsection{Induction and restriction}\label{subsec:indres}
Let $\bfA^{\omega}_R\subset \bfA^{\omega}$ be the subalgebra generated by $f\bfe(\lambda)$ and $\tau^{\omega}_{a}\bfe(\lambda)$ for $\lambda\in W_S\lambda_0$, $f\in \bfC[V]$ and $a\in \Delta_0$. For $\lambda_1\in W_S\lambda_0$, denote $\bfe_{R,\lambda_1} = \sum_{\lambda\in W_R \lambda_1}\bfe(\lambda)$ and define $\bfA^{\omega}_{R, \lambda_1} = \bfe_{R, \lambda_1}\bfA^{\omega}_R\bfe_{ R,\lambda_1}$ to be the idempotent subalgebra. In other words, $\bfA^{\omega}_{R, \lambda_1}$ is the subalgebra of $\bfA^{\omega}$ generated by $f\bfe(\lambda)$ and $\tau^{\omega}_{a}\bfe(\lambda)$ for $\lambda\in W_R\lambda_1$, $f\in \bfC[V]$ and $a\in \Delta_0$. \par

For each $\lambda_1\in W_S \lambda_0$, we define the induction, restriction and co-induction functors
\begin{align*}
	&\ind_{R,\lambda_1}^{S}: \bfA^{\omega}_{R, \lambda_1}\gmod \to\bfA^{\omega}\gmod,\quad N\mapsto \bfA^{\omega}\bfe_{R, \lambda_1}\otimes_{\bfA^{\omega}_{R, \lambda_1}}N\\
	&\res_{R,\lambda_1}^{S}: \bfA^{\omega}\gmod \to\bfA^{\omega}_{R, \lambda_1}\gmod,\quad M\mapsto \bfe_{R, \lambda_1}M\cong \gHom_{\bfA^{\omega}}\left( \bfA^{\omega}\bfe_{R, \lambda_1}, M \right) \\
	&\coind_{R,\lambda_1}^{S}: \bfA^{\omega}_{R, \lambda_1}\gmod \to\bfA^{\omega}\gmod,\quad N\mapsto \bigoplus_{\lambda\in W_S \lambda_0}\gHom_{\bfA^{\omega}_{R, \lambda_1}}\left( \bfe_{R, \lambda_1}\bfA^{\omega}\bfe(\lambda), N \right).\\
\end{align*}
\index[ch2]{ind@$\ind^S_R$}
They form a triplet of adjoint functors $\left( \ind_{R,\lambda_1}^{S}, \res_{R,\lambda_1}^{S}, \coind_{R,\lambda_1}^{S} \right)$
\begin{prop}\label{prop:ind}
	The functors $\ind_{R,\lambda_1}^{S}, \res_{R,\lambda_1}^{S}$ and $\coind_{R,\lambda_1}^{S}$ are exact.
\end{prop}
\begin{proof}
	The functor $\res_{R,\lambda_1}^{S}$ is clearly exact. By~\autoref{theo:basis}, we have a decomposition of right $\bfA^{\omega}_{R, \lambda_1}$-module
	\begin{align}\label{equa:decompA}
		\bfA^{\omega}\bfe_{R, \lambda_1}\cong \bigoplus_{w\in W^R} \tau^{\omega}_w \bfA^{\omega}_{R, \lambda_1}
	\end{align}
	where $W^R\subset W_S$ is the set of shortest representatives of the elements in $W_S / W_R$ and $\tau^{\omega}_w = \bigoplus_{\lambda\in W_R \lambda_1}\tau^{\omega}_{a_l}\cdots\tau^{\omega}_{a_1}\bfe\left( \lambda \right)$ for any reduced expression $w = s_{a_l} \cdots s_{a_1}$. Therefore $\bfA^{\omega}\bfe_{R, \lambda_1}$ is a free right $\bfA^{\omega}_{R, \lambda_1}$-module, so $\ind^{S}_{R, \lambda_1}$ is exact. Similarly, $\coind^{S}_{R, \lambda_1}$ is also exact.
\end{proof}

\section{Filtered \texorpdfstring{$\bfA^{\omega}$}{A omega}-modules}\label{sec:Amodfilt}

We consider $\bfA^{\omega}$-modules equipped with filtrations which are compatible with the filtration by length $F$ on $\bfA^{\omega}$. Most results in this section are non-unital version of the classical theory of filtered rings and filtered modules which one can find in~\cite{hotta2007d}. The goal of this section is to introduce (\autoref{subsec:supp}) the support and the Gelfand--Kirillov dimension of an object $M\in \bfA^{\omega}\gmod_0$ and show (\autoref{prop:indtf}) that ``induced modules'' have the full support.

\subsection{Good filtrations on \texorpdfstring{$\bfA^{\omega}$}{A omega}-modules}
Let $M\in \bfA^{\omega}\gmod$.
\begin{defi}\label{defi:FM}
	A good filtration $F$ on $M$ is a sequence $\{ F_{\le n} M\}_{n\in \bfZ}$ of graded $\bfC[V]$-submodules of $M$ satisfying the following properties:
	\begin{enumerate}
		\item
			$F_{\le n-1}\subseteq F_{\le n}$ for all $n\in \bfZ$;
		\item
			for each $n\in \bfZ$, there exists a finite subset $\Sigma_n\subset W_S \lambda_0$ such that\footnote{We require this condition because we work with a non-unital associative algebra.}
			 \begin{align*}
				 F_{\le n}M = \bigoplus_{\lambda\in \Sigma_n}\bfe(\lambda)F_{\le n}M;
			 \end{align*}
		\item
			$F_n M = 0$ for $n \ll 0$;
		\item
			$\bigcup_{n\in \bfZ} F_n M = M$;
		\item
			\begin{align*}
				\left(F_{\le n}\bfA^{\omega}\right)\left(F_{\le m}M\right) \subseteq F_{\le n+m}M, \quad \forall n, m\in \bfN;
			\end{align*}
		\item
			there exists $m_0 \gg 0$ satisfying  
			\begin{align*}
				\left(F_{\le n}\bfA^{\omega}\right)\left(F_{\le m}M\right) = F_{\le n+m}M,\quad \forall n\ge 0,\forall m\ge m_0.
			\end{align*}
	\end{enumerate}
\end{defi}
The following result is standard, see~\cite[D.1.3]{hotta2007d}
\begin{prop}\label{prop:bonfilt}
	Good filtrations exist for the objects of $\bfA^{\omega}\gmod$. If $F$ and $F'$ are two good filtrations on $M\in \bfA^{\omega}\gmod$, then there exists $i_0 \gg 0$ such that
	\[
		F'_{\le n -i_0}M \le F_{\le n}M \le F'_{\le n +i_0}M, \quad \forall n\in \bfZ.\pushQED{\qed}\qedhere\popQED
	\]
\end{prop}
The following lemma is a direct consequence of~\autoref{prop:bonfilt}.
\begin{coro}\label{lemm:FF'}
	If $F$ and $F'$ are good filtrations on $M$, then there exist
	\begin{enumerate}
		\item
			a finite filtration of $\gr^F\bfA^{\omega}$-submodules $F'$ on $\gr^F M$,
		\item
			a finite filtration of $\gr^F\bfA^{\omega}$-submodules $F$ on $\gr^{F'} M$ and
		\item
			an isomorphism of $\gr^F\bfA^{\omega}$-modules $\gr^{F'}\gr^{F}M\cong \gr^{F}\gr^{F'}M$.
	\end{enumerate}
\end{coro}
\begin{proof}
	By~\autoref{prop:bonfilt}, there exists $i_0 \gg 0$ such that $F_{\le n - i_0} M \le F'_{\le n} M \le F_{\le n + i_0} M$ for all $n\in \bfZ$. For $m\in \left[ -i_0, i_0 \right]$, define $F'_{\le n, \le m} M = \left(F'_{\le n}M \cap F_{\le n + m}M\right) + F'_{\le n-1} M$. Then the quotient $\gr^{F'}M$ acquires a filtration
	\begin{align*}
		F_{\le m}\gr^{F'}_n M = F'_{\le n, m} M / F'_{\le n-1} M \subseteq  F'_{\le n} M / F'_{\le n-1} M = \gr^{F'}_n M,
	\end{align*}
	which satisfies $\left(\gr^F_l\bfA^{\omega}\right)\left(F_{\le m}\gr^{F'}_n M\right)\subseteq F_{\le m}\gr^{F'}_{n+l} M$. Hence for each $m\in \left[ -i_0, i_0 \right]$, the quotient $\Gr^{F}_m \Gr^{F'}M = F_{\le m}\Gr^{F'}M / F_{\le m-1}\Gr^{F'}M$ is itself a $\Gr \bfA^{\omega}$-module. Similarly, we put $F_{\le m, \le n} M = \left(F_{\le m}M \cap F'_{\le m + n}M\right) + F_{\le m-1} M$ so that $\gr^F M$ acquires a filtration by $\gr^F \bfA^{\omega}$-modules.  
	Zassenhaus lemma yields $\gr_{m-n}^{F}\gr^{F'}_n M\cong \gr^{F'}_{n - m}\gr_{m}^{F} M$. Therefore, 
	\begin{align*}
		\bigoplus_{n = -i_0}^{i_0}\gr^{F'}_n\gr^F M\cong \bigoplus_{m = -i_0}^{i_0}\gr^{F}_m\gr^{F'} M.
	\end{align*}
\end{proof}

\subsection{Associated graded of good filtrations}

Recall the monoid algebra $\bfC Q^{\vee}_{\dom}$ from~\autoref{subsec:grA}. Given a good filtration $F$ on an object $M\in \bfA^{\omega}\gmod$, the associated graded $\Gr^F M = \bigoplus_{k\in \bfZ} F_{\le k} M / F_{\le k - 1} M$ is a $\Gr^F \bfA^{\omega}$-module. The $\gr^F\bfA^{\omega}$-action on $\gr^F M$ extends to an action of the unital completion introduced in~\autoref{subsec:grA} via the natural inclusion
\[
	\gr^F\bfA^{\omega}\hookrightarrow \End_{(\gr^F\bfA^{\omega})^{\op}}(\gr^F\bfA^{\omega}) \cong \prod_{\lambda\in W_S\lambda_0}\gr^F\bfA^{\omega}\bfe(\lambda).
\]
We obtain a $\bfC Q^{\vee}_{\dom}$-module structure on $\Gr^F M$ via the map~\eqref{equa:Theta}. 
\begin{prop}\label{prop:coherence}
	Let $M\in \bfA^{\omega}\gmod$ and $F$ a good filtration on $M$. Then $\Gr^F M$ is a coherent $\bfC Q^{\vee}_{\dom}\otimes \bfC[V]$-module. Moreover, if $M\in \bfA^{\omega}\gmod_0$, then $\Gr^F M$ is a coherent $\bfC Q^{\vee}_{\dom}$-module.
\end{prop}
\begin{proof}
	We observe that the coherence for $\gr^{F} M$ is independent of the choice of the good filtration $F$. Indeed, if $F'$ is another good filtration on $M$, then by~\autoref{lemm:FF'}, 
	\begin{align*}
		\gr^F M\text{ coherent} \Leftrightarrow \bigoplus_{n = -i_0}^{i_0}\gr^{F'}_n\gr^F M\cong \bigoplus_{m = -i_0}^{i_0}\gr^{F}_m\gr^{F'} M\text{ coherent} \Leftrightarrow \gr^{F'} M\text{ coherent}.
	\end{align*}

	We prove the first assertion. By~\autoref{prop:gmod} and the compactness of $M$, there is a surjection of the form
	\begin{align*}
			p:\bigoplus_{j = 1}^r\bfA^\omega\bfe(\lambda_j)\langle a_j\rangle\twoheadrightarrow M.
	\end{align*}
	Equip the source of $p$ with the length filtration and the target of $p$ with the induced filtration, denoted by $F$, so that $p$ induces a surjection on the associated graded $\gr^F \bfA^{\omega}$-module. The coherence on the source of $\gr^F p$ implies that of $\gr^F M$. Thus we may suppose that $M$ is of the form $M = \bfA^{\omega}\bfe(\lambda_j)$ and equipped with the length filtration. It follows from~\autoref{theo:PBW} that
	\begin{align*}
			\Gr^F\bfA^{\omega}\bfe(\lambda_j) \cong \calE\otimes_{\bfC} \left(\bigoplus_{w\in W_R}\tau^{\omega}_w\bfC[V]\right)\bfe(\lambda_j).
	\end{align*}
	Since $\calE$ is coherent over~$\bfC Q^{\vee}_{\dom}$ by~\autorefitem{lemm:CM}{ii} and $\bigoplus_{w\in W_R} \bfC[V]\tau^{\omega}_w$ is free of finite rank over $\calZ$, it follows that $\Gr^F\bfA^{\omega}\bfe(\lambda_j)$ is coherent over $\bfC Q^{\vee}_{\dom}\otimes \calZ$. \par
	Suppose now $M\in \bfA^{\omega}\gmod_0$ so that $\calZ$ acts via the quotient $\calZ / \frakm_{\calZ}^n$ for some $n\in \bfN$. Since $\calZ / \frakm_{\calZ}^n$ is finite-dimensional, $M$ must be coherent over $\bfC Q^{\vee}_{\dom}$.
\end{proof}

\subsection{Support of \texorpdfstring{$\bfA^{\omega}$}{A omega}-modules of finite length} \label{subsec:supp}
Let $M\in \bfA^{\omega}\gmod_0$. In view of~\autoref{prop:coherence}, we can make the following definition:
\begin{defi}
	The support of $M$, denoted by $\supp M$, is defined to be the support of $\Gr^F M$ as coherent $\bfC Q^{\vee}_{\dom}$-module, for any choice of good filtration $F$ on $M$.
\end{defi}\index[ch2]{supp@$\supp$}
By~\autoref{lemm:FF'}, the definition of $\supp M$ is independent of the choice of a good filtration. \par

We define the Gelfand--Kirillov dimension of a weight module $M$ of $\bfA^{\omega}$ to be the following number: upon choosing a good filtration $F$ on $M$,
\begin{align*}
	\dim_{\GK}M = \limsup_{n\to\infty}\frac{\log\dim F_{\le n}M}{\log n}.
\end{align*}\index[ch2]{dim@$\dim_{\GK}$} 
By~\autoref{prop:bonfilt}, this number does not depend on the choice of $F$.
\begin{prop}\label{prop:dimGK}
	Let $M\in \bfA^{\omega}\gmod_0$. Then the Gelfand--Kirillov dimension $\dim_{\GK}M$ coincides with the Krull dimension of $\supp M$.
\end{prop}
\begin{proof}
	Taking the associated graded, we have
	\begin{align*}
		\dim F_{\le n}M = \dim \bigoplus_{k = -\infty}^n\gr^F_k M.
	\end{align*}
	Notice that $\bfC Q^{\vee}_{\dom}$ is finitely generated graded ring, where $\deg X^{\mu} = \ell(X^{\mu})$, and $\gr^FM$ is a finitely generated graded module over it. Hence $\dim_{\GK} M$ is nothing but the degree of the Hilbert polynomial of $\gr^F M$, which is equal to the Krull dimension of $\supp M$.
\end{proof}

\subsection{Induction of filtered modules}
Recall the subalgebra $\bfA^{\omega}_{R, \lambda_1}\subset \bfA^{\omega}$ from~\autoref{subsec:indres}. Good filtrations on objects of $\bfA^{\omega}_{R, \lambda_1}\gmod$ are defined in a similar manner. \par
Suppose $N\in \bfA^{\omega}_{R, \lambda_1}\gmod$ is equipped with a good filtration $F$ which satisfies $F_{\le k}N = \left( F_{\le k} \bfA^{\omega}_{R, \lambda_1} \right)\left( F_{\le 0} N \right)$ for $k\ge 0$ and $F_{\le -1} N = 0$.  \par

Let $M = \ind^{S}_{R, \lambda_1} N$. The adjunction unit yields an inclusion of $\calZ$-modules $N\hookrightarrow M$. Define a filtration $F_{\le n} M = \left(F_{\le n} \bfA^{\omega}  \right)\left( F_{\le 0} N \right)$.
\begin{lemm}\label{lemm:indfilt}
	The filtration $F$ on $M$ is good and satisfies
	\begin{align*}
		\Gr^F M\cong \left(\Gr^F \bfA^{\omega} \bfe_{R, \lambda_1}\right)\otimes_{\Gr^F\bfA^{\omega}_{R, \lambda_1}}\left(\Gr^F N\right).
	\end{align*}
\end{lemm}
\begin{proof}
	By the hypothesis on $F_{\le n}N$, we have $\Gr^F_{n} N = \left(\Gr^F_{n}\bfA^{\omega}_{R, \lambda_1}\right)\left(\Gr^F_0 N\right)$ and $\Gr^F_{n} M = \left(\Gr^F_{n}\bfA^{\omega}\right)\left(\Gr^F_0 N\right)$. By the decomposition~\eqref{equa:decompA}, we deduce
	\begin{align*}
		& \Gr^F_{k}\bfA^{\omega}\bfe_{R, \lambda_1} = \bigoplus_{j = 0}^{k}\bigoplus_{\substack{w\in W^R \\ \ell(w) = j}}\tau^{\omega}_w\Gr^F_{k-j}\bfA^{\omega}_{R, \lambda_1},
	\end{align*}
	from which
	\begin{align*}
		\left(\left(\Gr^F \bfA^{\omega} \bfe_{R, \lambda_1}\right)\otimes_{\Gr^F\bfA^{\omega}_{R, \lambda_1}}\left(\Gr^F N\right)\right)_n &=
		\bigoplus_{j=0}^n\bigoplus_{\substack{w\in W^R \\ \ell(w) = j}}\tau^{\omega}_w\left(\Gr^F_{n-j}\bfA^{\omega}_{R,\lambda_1} \right)(\Gr^F_0 N) \\
		&= \left(\Gr^F_{n}\bfA^{\omega}_n\right)\left( \Gr^F_0 N \right) = \Gr^F_n M
	\end{align*}
	
\end{proof}
\begin{prop}\label{prop:indtf}
	For any $N\in \bfA^{\omega}\gmod_0$ and $0\neq M'\subset \ind^{S}_{R, \lambda_1} N$, we have $\Supp M' = \Spec \bfC Q^{\vee}_{\dom}$. 
\end{prop}
\begin{proof}
	Let $F$ be a good filtration on $N$ as above and denote $M = \ind^S_{R, \lambda_1} N$, so that $\Gr^F M\cong \Gr^F \left(\bfA^{\omega}\bfe_{R, \lambda_1}\right)\otimes_{\Gr^F\bfA^{\omega}_{R, \lambda_1}}\left(\Gr^F N\right)$ by~\autoref{lemm:indfilt}. By~\autoref{theo:PBW}, we have $\Gr^F \bfA^{\omega}\bfe_{R, \lambda_1} \cong \calE\otimes_{\bfC} \gr^F\bfA^{\omega}_{R, \lambda_1}$; hence 
	\begin{align*}
			\Gr^F M \cong \calE\otimes_{\bfC}\Gr^F N.
	\end{align*}
	By~\autorefitem{lemm:CM}{ii}, $\calE$ and thus $\gr^F M$ are a Cohen--Macaulay module of maximal dimension over $\bfC Q^{\vee}_{\dom}$, so it is torsion-free. For any $0\neq M'\subset M$, the restriction to $M'$ of $F$ is a good filtration and $\Gr^F M'\subset \Gr^F M$. Hence $\Supp M' = \bfC Q^{\vee}_{\dom}$.
\end{proof}
\begin{rema}
	\autoref{prop:indtf} is an analogue of the following basic property for a double affine Hecke algebra $\bbH$: the induced module $\bbH\otimes_{\ubar\bbH} M$ is free over the polynomial part $\bfC[E]\subset \bbH$ for every module $M$ over the graded affine Hecke algebra $\ubar\bbH\subset \bbH$. A similar property for rational Cherednik algebras was used in~\cite{GGOR03} in the proof of the double centraliser property of the KZ functor. Our proof of the double centraliser property~\autoref{theo:bicommutante} also relies on it.
\end{rema}

\section{Quiver Hecke algebras}\label{sec:MQHA} 

We keep the notation of root systems $\left( E,S, \Delta \right)$ and $(V, R, \Delta_0)$. In this section, we introduce an algebra $\bfB^{\omegaa}$, which can be viewed as a variant of quiver Hecke algebras. The relation between the quiver Hecke algebras and $\bfB^{\omegaa}$ in the case where the root system $(V, R)$ is of type A is explained in~\autoref{rema:KLR}.

\subsection{The algebra \texorpdfstring{$\bfB^{\omegaa}$}{B omegaa}}\label{subsec:Bomega}
Define the torus $T = Q^{\vee}\otimes \bfC^{\times}$\index[ch2]{T@$T$} so that the ring of regular functions $\bfC[T]$ is isomorphic to the group algebra $\bfC P$. For any $\alpha\in P$, we denote by $Y^{\alpha}\in \bfC[T]$ the corresponding element. \par

Fix $\ell_0\in T$. Define for each $\ell\in W_R \ell_0$ a polynomial ring $\Pol_\ell = \bfC[V]$ and let $\Pol_{W_R \ell_0} = \bigoplus_{\ell\in W_R\ell_0} \Pol_{\ell}$. For each $\ell$, define $\bfe(\ell):\Pol_{W_R \ell_0}\to \Pol_{\ell}$ to be the idempotent linear endomorphism of projection onto the factor $\Pol_\ell$. Recall that $R_{\idv} = R^+ \setminus 2R$\index[ch2]{R@$R_\idv$} and $R^+_{\idv} = R_\idv \cap R^+$. \par

Choose any $\lambda_0\in \exp^{-1}(\ell_0)$. Then the algebra $\calZ$ from~\autoref{subsec:Z} acts on $\Pol_{\ell}$: for any $w\in W_R$, the element $f\in \calZ = \bfC[V]^{W_{\lambda_0}}$ acts on $\Pol_{w\ell_0}$ by multiplication by $w(f)$. \par

Let $\omegaa = \left\{ \omegaa_{\ell}\right\}_{\ell\in W_R\ell_0}$\index[ch2]{omega@$\omegaa_\ell$} be a family of functions $\omegaa_{\ell}:R^+_\idv\to \bfZ_{\ge -1}$ satisfying the properties:
\begin{enumerate}
	\item
		If $2\alpha\notin R$, then $\omegaa_\ell(\alpha) = -1$ implies $Y^{\alpha}(\ell) = 1$.
	\item
		If $2\alpha\in R$, then $\omegaa_\ell(\alpha) = -1$ implies $Y^{\alpha}(\ell) \in \left\{ 1, -1 \right\}$.
	\item
		For $w\in W_R$ and $\alpha\in R^+_\idv\cap w^{-1}R^+_\idv$ we have $\omegaa_{\ell}(\alpha) = \omegaa_{w\ell}(w\alpha)$. 
\end{enumerate}
	For each $\alpha\in \Delta_0$ and $\ell\in W_R \ell_0$, we define an operator $\tau^\omegaa_\alpha\bfe(\ell):\Pol_{\ell}\to \Pol_{s_\alpha\ell}$ by
\begin{align*}
	\tau^\omegaa_\alpha\bfe(\ell) = 
	\begin{cases} \alpha^{-1}(s_{\alpha} - 1) & \omegaa_{\ell}(\alpha) = -1 \\
		\alpha^{\omegaa_{\ell}(\alpha)} s_{\alpha} & \omegaa_{\ell}(\alpha) \ge 0 \\
	\end{cases}.
\end{align*}\index[ch2]{tau@$\tau^{\omegaa}_\alpha$} 
Here $s_{\alpha}: \bfC[V]\to \bfC[V]$ is the reflection with respect to $\alpha$. 
\begin{defi}
	We define $\bfB^{\omegaa} = \bfB(R, V, \Delta, W_R \ell_0, \omegaa)$\index[ch2]{B@$\bfB^{\omegaa}$} to be the subalgebra of $\End_{\calZ}\left( \Pol_{W_R\ell_0} \right)$ generated by $\bfC[V]\bfe(\ell)$ and $\tau^{\omegaa}_{\alpha}\bfe(\ell)$. 
\end{defi}
All the statements of~\autoref{prop:Agmod} for $\bfA^{\omega}$ hold equally for $\bfB^{\omegaa}$. In particular, the centre of $\bfB^{\omegaa}$ is equal to $\calZ$. \par

\begin{exam}
	\noindent
	\begin{enumerate}
		\item
			If $\ell_0 = 1\in T$ and $\omegaa = \left\{ -1 \right\}_{\ell = \ell_0}$ is the $-1$ constant function, then $\bfB^{\omegaa}$ is the affine nil-Hecke algebra of type $W_R$ and is isomorphic to a matrix algebra over its centre.
		\item
			If $\ell_0 = 1\in T$ and $\omegaa = \left\{ 0 \right\}_{\ell =\ell_0}$ is the zero constant function, then $\bfB^{\omegaa} = \bfC[V]\rtimes W_R$ is the skew tensor product.
	\end{enumerate}
\end{exam}

\begin{rema}\label{rema:KLR}
	In the case where the finite root system $(V, R)$ is of type $A_{n-1}$, the algebra $\bfB^{\omegaa}$ recovers the notion of quiver Hecke algebras. \par
	For any quiver $\Gamma = (I, H)$ with $I\subset \bfC^{\times}$ and a dimension vector $\beta\in \bfN I$ with $|\beta| = n$, the quiver Hecke algebra, denoted by $R_{\beta}(\Gamma)$ according to~\cite{rouquier12}, is generated by three sets: idempotents $\left\{ \bfe(\ell) \right\}_{\ell\in I^{\beta}}$, Hecke operators $\left\{\tau_i\right\}_{i = 1}^{n-1}$, polynomial part $\left\{ x_i \right\}_{i = 1}^n$. By translating suitably the set $I\subset \bfC^{\times}$, we may assume that $\prod_{r\in I}r^{\beta_r} = 1$, so that each sequence $\nu = (\nu_1, \cdots, \nu_n)\in I^{\beta}\subset (\bfC^{\times})^n$ lies in the maximal torus $T\subset(\bfC^{\times})^n$ of $\SL_n(\bfC)$. We put $\omegaa_{\nu}(\alpha_{i,j}) = \#\{(h: i\to j)\in H\}- \delta_{\nu_i = \nu_j}$. Then there is a surjective homomorphism 
	\begin{align*}
		R_{\beta}(\Gamma)&\to \bfB^{\omegaa} \\
		\bfe(\nu) &\mapsto \bfe(\nu) \\
		\tau_i &\mapsto \tau_{\alpha_i},\quad i \in \left\{ 1, \ldots, n-1 \right\}\\
		x_k &\mapsto \frac{1}{n}\left(-\sum_{1\le j < k} j\alpha_j + \sum_{k\le j < n}(n-j)\alpha_j\right)\bfe_R,\quad k\in \left\{ 1, \ldots, n \right\} \\
	\end{align*}
	whose kernel the ideal generated by $x_1 + \cdots + x_n$. .
\end{rema}

\subsection{Basis theorem}
\begin{theo}\label{theo:basisB}
	For any $w\in W_R$, choose a reduced expression $w = s_{a_1} \cdots s_{a_l}$ and put $\tau^{\omegaa}_w = \tau^{\omegaa}_{\alpha_l}\cdots \tau^{\omegaa}_{\alpha_1}$. Then there is a decomposition
	\begin{align*}
		\bfB^{\omegaa} = \bigoplus_{\ell\in W_R \ell_0}\bigoplus_{\substack{w\in W_R \\ \ell(w) = n}} \bfC[V]\tau^{\omegaa}_w\bfe(\ell).
	\end{align*}
\end{theo}
\begin{proof}
	To prove it, we shall apply the results~\autoref{theo:basis-generic} and~\autoref{theo:isomKH} whose proofs do not rely on this theorem. By~\autoref{lemm:choixomega} below, we can choose $\omega = \left\{ \omega_{\lambda} \right\}_{\lambda\in W_S\lambda_0}$ such that $\intt\omega = \omegaa$. Then \autoref{theo:isomKH} implies that upon choosing a good $\gamma\in Q^{\vee}$, there is an isomorphism $\bfB^{\omegaa}\cong \bfe_{\gamma}\bfA^{\omega}\bfe_{\gamma}$ identifying $\tau^{\omegaa}_{\alpha}\bfe(\ell)$ with $\sigma_{\alpha}\bfe(\pre\gamma{\ell})$ and by~\autoref{theo:basis-generic}, the idempotent subalgebra $\bfe_{\gamma}\bfA^{\omega}\bfe_{\gamma}$ has a decomposition in terms of $\sigma_{\alpha}\bfe(\pre\gamma{\ell})$. Hence $\bfB^{\omegaa}$ also has a decomposition as in the statement.
\end{proof}
\begin{lemm}\label{lemm:choixomega}
	Given any family of order functions $\omegaa = \left\{ \omegaa_{\ell} \right\}_{\ell\in W_R\ell_0}$ for $\bfB^{\omegaa}$, there exists a family of order functions $\omega = \left\{ \omega_{\lambda} \right\}_{\lambda\in W_S \lambda_0}$ satisfying the conditions from~\autoref{subsec:Aomega} such that $\intt \omega = \omegaa$, where $\intt\omega$ is defined in~\autoref{subsec:formuleomega}. 
\end{lemm}
\begin{proof}
We choose a point $\lambda_0\in \exp^{-1}(\ell_0)\subset V$. 
	Such a family $\omega = \left\{ \omega_{\lambda} \right\}_{\lambda\in W_S \lambda_0}$ is determined by a $W_{\lambda_0}$-invariant function $\til\omega_{\lambda_0}: S\to \bfZ_{\ge -1}$ and it suffices to construct it. However, one needs to be careful about the condition (i) from~\autoref{subsec:Aomega}. We first define a function $\til\omegaa_{\ell_0}: R\to \bfZ_{\ge -1}$ as follows:
	\begin{enumerate}
		\item
			For any $\alpha\in R^+_\idv$ such that $2\alpha\notin R$, we set $\til\omegaa_{\ell_0}(\alpha) = \omegaa_{\ell_0}(\alpha)$ and $\til\omegaa_{\ell_0}(-\alpha) = \omegaa_{w_0\ell_0}(-w_0 \alpha)$.
		\item
			For any $\alpha\in R^+_\idv$ such that $2\alpha\in R$, if $Y^{\alpha}(\ell_0) = -1$, then we set $\til\omegaa_{\ell_0}(2\alpha) = \omegaa_{\ell_0}(\alpha)$, $\til\omegaa_{\ell_0}(-2\alpha) = \omegaa_{w_0\ell_0}(-w_0 \alpha)$ and $\til\omegaa_{\ell_0}(\alpha) = \til\omegaa_{\ell_0}(-\alpha) = 0$; otherwise, we set $\til\omegaa_{\ell_0}(\alpha) = \omegaa_{\ell_0}(\alpha)$, $\til\omegaa_{\ell_0}(-\alpha) = \omegaa_{w_0\ell_0}(-w_0 \alpha)$ and $\til\omegaa_{\ell_0}(2\alpha) = \til\omegaa_{\ell_0}(-2\alpha) = 0$. 
	\end{enumerate}
	The function $\til\omegaa_{\ell_0}$ is $W_{\ell_0}$-invariant by the assumption (iii) from~\autoref{subsec:Bomega} and has image in $\bfZ_{\ge -1}$. We choose a section of the projection $W_{\lambda_0}\backslash S\to W_{\ell_0}\backslash R$, denoted $f:W_{\ell_0}\backslash R \to W_{\lambda_0}\backslash S$, in such a way that for each $\alpha\in R_\idv$, the condition $f(\alpha)(\lambda_0) = 0$ holds whenever $\alpha(\ell_0) = 0$ . We set $\til\omega_{\lambda_0} = f_* \til\omegaa_{\ell_0}$ so that $\til\omega_{\lambda_0}:S\to \bfZ_{\ge -1}$ is a $W_{\lambda_0}$-invariant function of finite support. The family $\left\{ \omega_{\lambda} \right\}_{\lambda\in W_S \lambda_0}$ is then defined by $\omega_{w\lambda_0}(a) = \til\omega_{\lambda_0}(w^{-1}a)$ for all $w\in W_S$ and $a\in S^+$. \par
\end{proof}

\subsection{Frobenius form on \texorpdfstring{$\bfB^{\omegaa}$}{B omega}}
As observed in~\cite{BGS08}, the basis theorem~\autoref{theo:basisB} implies that the algebra $\bfB^{\omegaa}$ is Frobenius over its centre $\calZ$. 
\begin{lemm}\label{lemm:Frob}
	$\bfB^{\omegaa}$ is a Frobenius algebra over $\calZ$.
\end{lemm}
\begin{proof}
	Consider the filtration by length
	\begin{align*}
		F_{\le n}\bfB^\omegaa = \sum_{\ell\in W_R \ell_0}\sum_{k =0}^n\sum_{\left( \alpha_1, \ldots, \alpha_k \right)\in \Delta_0^k} \bfC[V]\tau^{\omegaa}_{\alpha_1}\cdots \tau^{\omegaa}_{\alpha_k}\bfe(\ell).
	\end{align*}
	We set $N = \#R^+ = \ell(w_0)$ and let $w_0 = s_{\alpha_N}\cdots s_{\alpha_1}$ be any reduced expression for the longest element $w_0\in W_R$ and set $\tau^{\omegaa}_{w_0}\bfe(\ell) = \tau^{\omegaa}_{\alpha_N}\cdots \tau^{\omegaa}_{\alpha_1}\bfe(\ell)$. By~\autoref{theo:basisB}, we have $F_{\le N}\bfB^{\omegaa} = \bfB^{\omegaa}$ and 
	\begin{align*}
		\Gr^F_N \bfB^{\omegaa}\cong \bigoplus_{\ell\in W_R\ell_0}\bfC[V]\tau^{\omegaa}_{w_0}\bfe(\ell). 
	\end{align*}
	Let $R_{\lambda_0} = \left\{ \alpha\in R\;;\; \alpha(\lambda_0) = 0 \right\}$ be the sub-root system associated with $\lambda_0$ and let $\Delta_{\lambda_0}\subset R_{\lambda_0}$ be any basis, which determines a set of positive roots $R^+_{\lambda_0}\subset R_{\lambda_0}$ and a set of Coxeter generators $\{s_a\}_{a\in \Delta_{\lambda_0}}\subset W_{\lambda_0}$. It is well known that $\bfC[V]$ is a symmetric algebra over $\calZ$ with the trace map $f\mapsto \vartheta_{w_0(W_{\lambda_0})}(f)$, where $\vartheta_{w_0(W_{\lambda_0})}$ is a composition of Demazure operators for the longest element $w_0(W_{\lambda_0})$ of the Coxeter group $(W_{\lambda_0}, \Delta_{\lambda_0})$.  Let $\tr$ be the composition
	\begin{align*}
		\bfB^\omegaa &\to \Gr^F_N\bfB^{\omegaa} = \bigoplus_{\ell} \bfC[V]\tau_{w_0}\bfe(\ell) \xrightarrow{\tau_{w_0}\bfe(\ell)\mapsto 1}\bigoplus_{\ell} \bfC[V] \\
		&\xrightarrow{\vartheta_{w_0(W_{\ell})}}\bigoplus_{\ell}\bfC[V]^{W_{\ell}} \cong\bigoplus_{\ell} \calZ \xrightarrow{\sum_{\ell\in W_R\ell_0}} \calZ.
	\end{align*}
	Then $\tr$ is a Frobenius form.
\end{proof}

\section{Knizhnik--Zamolodchikov functor \texorpdfstring{$\bfV$}{V}}\label{sec:KZ}
We resume to the assumptions of~\autoref{sec:agha}. \par

In this section, we introduce a functor $\bfV:\bfA^\omega\gmod\to \bfB^{\omegaa}\gmod$, which is a quotient functor satisfying the double centraliser property. It can be viewed as a generalisation of the monodromy functor of~\cite{VV04} for dDAHAs (which has been reviewed in~\autoref{subsec:monodromie}) to the family of algebras $\bfA^{\omega}$. It is thus expected to satisfy some properties of the monodromy functor. The main results of this article~\autoref{theo:bicommutante} and~\autoref{prop:catV} provide some evidence. We construct $\bfV$ by choosing an idempotent element $\bfe_{\gamma}\in\bfA^{\omega}$ and establish an isomorphism $\bfB^{\omegaa}\cong \bfe_{\gamma}\bfA^{\omega}\bfe_{\gamma}$ in~\autoref{theo:isomKH}. 

\subsection{The idempotent construction}\label{subsec:idemp}
Consider the following exponential map
\begin{align*}
	E\cong V = Q^{\vee}\otimes \bfR &\xrightarrow{\exp} Q^{\vee}\otimes \bfC^{\times} = T\\
	\mu\otimes r&\mapsto \mu\otimes e^{2\pi i r}
\end{align*}
and put $\ell_0 = \exp(\lambda_0)\in T$. Choose an element $\gamma\in Q^{\vee}$ such that
\begin{align}\label{equa:gamma}
	\langle \gamma, \alpha\rangle \ll 0\quad \text{ for all } \alpha\in R^+. 
\end{align}
We define a section of the projection $\dW: W_S\to W_S / Q^{\vee} = W_R$ by 
\begin{align*}
	\pre{\gamma}\bullet: W_R &\to W_S \\
	w&\mapsto X^{\gamma}\,w\, X^{- \gamma} = w X^{w^{-1}\gamma - \gamma}
\end{align*}\index[ch2]{w@$\pre{\gamma}w,\pre{\gamma}\lambda$} 
and a section of the exponential map $W_S\lambda_0\xrightarrow{\exp} W_R \ell_0$ by
\begin{align*}
	\pre{\gamma}\bullet: W_R\ell_0 &\to W_S \lambda_0 \\
	w\ell_0 &\mapsto  X^{\gamma}\,w\lambda_0.
\end{align*}
It is clear that $\pre{\gamma}w\pre{\gamma}\ell = \pre{\gamma}(w\ell)$. The choice of $\gamma$ implies that
\begin{align}\label{equa:gammaell}
	\alpha(\pre{\gamma}\ell) \ll 0\quad \text{ for all } \alpha\in R^+\text{ and } \ell\in W_R\ell_0. 
\end{align}

Given a family of order functions $\left\{ \omega_{\lambda}:S^+\to \bfZ_{\ge -1} \right\}_{\lambda\in W_S \lambda_0}$ satisfying the axioms of~\autoref{subsec:Aomega}, we can associate a family of order functions $\intt\omega = \left\{ \intt\omega_{\ell}:R^+_\idv\to \bfZ_{\ge -1} \right\}_{\ell\in W_R \ell_0}$\index[ch2]{@$\int\omega$}, called the {\itshape integral} of $\omega$ along $\dW$, by setting for each $\ell\in W_R \ell_0$
\begin{align}\label{equa:partialomega}
	\intt\omega_{\ell}(\alpha) = \sum_{\substack{a\in S^+ \\ \partial a \in \{\alpha,2\alpha\}}}\omega_{\pre{\gamma}\ell}(a).
\end{align}
The definition of $\intt\omega$ is independent of the choice of $\gamma$. Denote $\omegaa = \int\omega$. This family of order functions gives rise to an algebra $\bfB^{\omegaa}$ as defined in~\autoref{subsec:Bomega}. \par

For any $\ell$ and $\alpha\in \Delta_0$, we define an operator $\sigma_{\alpha} \bfe(\pre{\gamma}\ell):\Pol_{\pre{\gamma}\ell}\to\Pol_{\pre{\gamma}(s_{\alpha}\ell)}$ by
\begin{align}\label{equa:sigmaalpha}
	\sigma_{\alpha} \bfe(\pre{\gamma}\ell) = 
	\begin{cases} 
		\alpha^{-1} (s_\alpha - 1) & \intt\omega_{\ell}(\alpha) = -1 \\
		\alpha^{\omegaa_{\ell}(\alpha)} s_\alpha &  \intt\omega_{\ell}(\alpha) \ge 0 \\
	\end{cases}.
\end{align}
Define the idempotent
\begin{align}\label{equa:egamma}
	\bfe_{\gamma} = \sum_{\lambda\in \pre{\gamma}(W_R\ell_0)}\bfe(\lambda)\in \bfA^\omega.
\end{align}

The main result is the following, which will be proven in~\autoref{subsec:demo-isomKH}:
\begin{theo}\label{theo:isomKH}
	Upon choosing $\gamma\in Q^{\vee}$ satisfying~\eqref{equa:gamma}, there is an isomorphism of graded $\calZ$-algebras
	\begin{align*}
		i_{\gamma}: \bfB^{\omegaa} &\cong \bfe_{\gamma}\bfA^\omega\bfe_{\gamma} \\
		f\bfe(\ell)&\mapsto f\bfe(\pre{\gamma}\ell) \\
		\tau^{\omegaa}_\alpha\bfe(\ell) & \mapsto \sigma_{\alpha}\bfe(\pre{\gamma}\ell).
	\end{align*}
	Moreover, for any other choice $\gamma'$, the intertwiner 
	\[
		\varphi_{\gamma,\gamma'}:=\sum_{w\in W_R / W_{\ell_0}}\bfe(\pre{\gamma}(w\ell_0))\tau^{\omega}_{X^{w(\gamma - \gamma')}}\bfe(\pre{\gamma'}(w\ell_0))\in \bfe_{\gamma}\bfA^\omega\bfe_{\gamma'}
		\]
		yields a factorisation $i_{\gamma}(f) = \varphi_{\gamma,\gamma'}\cdot i_{\gamma'}(f)\cdot \varphi_{\gamma',\gamma}$ for each $f\in \bfB^{\omegaa}$.
	\end{theo}

\begin{exam}\label{exam:idemp}
	Resume to the setting of~\autorefitem{exam:Aomega}{iii} and~\autoref{exam:clan}. The coroot lattice is given by $Q^{\vee} = \bfZ$, which acts by translation on $E = \bfR$. Recall that $\lambda_0 = 1/4\in E$. We may take $\gamma = s_1s_0 = -1$ so that $\pre{\gamma}{\left(W_R\ell_0\right)} = \left\{ \lambda_+, \lambda_- \right\}$, where $\lambda_+ = s_1s_0\lambda_0 = -3/4$ and $\lambda_- = s_1s_0s_1 \lambda_0 = -5/4$. It follows that $\lambda_- = s_1s_0s_1s_0s_1 \lambda_+$ and 
	\begin{align*}
		\bfe(\lambda_-)\bfA^{\omega}\bfe(\lambda_+) = \bfC[\epsilon] \tau^{\omega}_{a_1}\tau^{\omega}_{a_0}\tau^{\omega}_{a_1}\tau^{\omega}_{a_0}\tau^{\omega}_{a_1}\bfe(\lambda_+),\quad \bfe(\lambda_+)\bfA^{\omega}\bfe(\lambda_-) = \bfC[\epsilon] \tau^{\omega}_{a_1}\tau^{\omega}_{a_0}\tau^{\omega}_{a_1}\tau^{\omega}_{a_0}\tau^{\omega}_{a_1}\bfe(\lambda_-). 
	\end{align*}
	Denote by $s:\bfC[\epsilon]\to \bfC[\epsilon]$ the automorphism $\epsilon \mapsto -\epsilon$. Calculate the products:
	\begin{align*}
		\tau^{\omega}_{a_1}\tau^{\omega}_{a_0}\tau^{\omega}_{a_1}\tau^{\omega}_{a_0}\tau^{\omega}_{a_1}\bfe(\lambda_+) &= \tau^{\omega}_{a_1}\bfe(5/4)\tau^{\omega}_{a_0}\bfe(-1/4)\tau^{\omega}_{a_1}\bfe(1/4)\tau^{\omega}_{a_0}\bfe(3/4)\tau^{\omega}_{a_1}\bfe(-3/4), \\
		&= s\cdot s\cdot (\epsilon s)\cdot s\cdot s = \epsilon s \\
		\tau^{\omega}_{a_1}\tau^{\omega}_{a_0}\tau^{\omega}_{a_1}\tau^{\omega}_{a_0}\tau^{\omega}_{a_1}\bfe(\lambda_-) &= \tau^{\omega}_{a_1}\bfe(3/4)\tau^{\omega}_{a_0}\bfe(1/4)\tau^{\omega}_{a_1}\bfe(-1/4)\tau^{\omega}_{a_0}\bfe(5/4)\tau^{\omega}_{a_1}\bfe(-5/4) \\
		&= s\cdot (-\epsilon s)\cdot s\cdot s\cdot s = \epsilon s.
	\end{align*}
	Let $\alpha = \partial a_1\in \Delta_0$ be the simple root for $(V,R) = A_1$. Denote $\ell_+ = \exp(2\pi i\lambda_+) = i$ and $\ell_- = \exp(2\pi i\lambda_-) = -i$. The family of order functions $\omegaa = \intt \omega$ for $\bfB^\omegaa$ is given by 
	\begin{align*}
		\omegaa_{\ell_+}(\alpha) = \sum_{k\in \bfN} \omega_{\lambda_+}(\alpha + k) = 1,\quad \omegaa_{\ell_-}(\alpha) = \sum_{k\in \bfN} \omega_{\lambda_-}(\alpha + k) = 1. \\
	\end{align*}
	It follows that 
	\begin{align*}
		\sigma_{\alpha}\bfe(\lambda_+) = \alpha^{\omegaa_{\ell_+}(\alpha)}s = \tau^{\omega}_{a_1}\tau^{\omega}_{a_0}\tau^{\omega}_{a_1}\tau^{\omega}_{a_0}\tau^{\omega}_{a_1}\bfe(\lambda_+)\quad \sigma_{\alpha}\bfe(\lambda_-) = \alpha^{\omegaa_{\ell_-}(\alpha)}s = \tau^{\omega}_{a_1}\tau^{\omega}_{a_0}\tau^{\omega}_{a_1}\tau^{\omega}_{a_0}\tau^{\omega}_{a_1}\bfe(\lambda_-)
	\end{align*}
	and therefore there is an isomorphism 
	\begin{align*}
		\bfB^{\omegaa} &\xrightarrow{\sim} \bfe_{\gamma}\bfA^{\omega}\bfe_{\gamma} \\
		\bfe(\ell_+) &\mapsto \bfe(\lambda_+) \\
		\bfe(\ell_-) &\mapsto \bfe(\lambda_-) \\
		\tau^{\omegaa}_\alpha &\mapsto \tau^{\omega}_{a_1}\tau^{\omega}_{a_0}\tau^{\omega}_{a_1}\tau^{\omega}_{a_0}\tau^{\omega}_{a_1}\bfe_{\gamma}.
	\end{align*}
\end{exam}
\begin{rema}
	As we will see in~\autoref{lemm:clans}, the idempotent $\bfe_{\gamma}$ corresponds to generic clans (\autoref{subsec:clans}). The choice of $\bfe_{\gamma}$ is inspired from the sheaf-theoretic study of extension algebras over a cyclically graded simple Lie algebra $\frakg_*$ in~\cite{liu19} and the sheaf-theoretic construction of the KZ functor. In the language of {\it op. cit.} and~\cite{LYI}, each eigenvalue $\lambda\in W_S\lambda_0$ corresponds to the spiral induction of a cuspidal local system $\scrC$ through one spiral of $\frakg_*$. On the other hand, affine Hecke algebras arise as extension algebra of parabolic inductions of $\scrC$ through parabolic subalgebras of $\frakg_*$, which appear also as spiral induction of $\scrC$ through ``generic spirals''. Therefore, the definition of the sheaf-theoretic KZ functor is nothing but picking idempotents of the extension algebra corresponding to those ``generic spirals''. In the algebraic and combinatorial language, they corresponds to alcoves lying in the generic clans, as introduced in~\autoref{subsec:clans}.
\end{rema}

\subsection{A formula for order functions}\label{subsec:formuleomega}
By the hypothesis of finite support for $\til\omega_{\lambda_0}:S\to \bfZ_{\ge -1}$, there exists $M \gg 0$ such that $\til\omega_{\lambda_0}(\alpha + k) = 0$ for all $\alpha\in R$ and $|k| \ge M$. Let $\gamma\in Q^{\vee}$ be an element satisfying~\eqref{equa:gamma}. More specifically, we require that
\begin{align}\label{equa:alphagamma}
	\langle \alpha, \gamma\rangle \le -M,\quad \forall \alpha\in R^+.
\end{align}
We prove a relation between the family $\omega = \left\{ \omega_{\lambda} \right\}_{\lambda\in W_S \lambda_0}$ for $\bfA^{\omega}$ and its integral $\omegaa = \left\{ \omegaa_{\ell} \right\}_{\ell\in W_R \ell_0}$ for $\bfB^{\omegaa}$ defined in~\eqref{equa:partialomega}. 
\begin{lemm}\label{lemm:prod}
	For any $\ell\in W_R \ell_0$ and $w\in W_R$, following formula holds in $\bfC(V)$:
	\begin{align*}
		\prod_{b\in S^+\cap\pre{\gamma}w^{-1} S^-}(-\partial b)^{\omega_\lambda(b)} = \epsilon\cdot\prod_{\beta\in R^+_\idv \cap w^{-1} R^-_\idv}(-\beta)^{\omegaa_\ell(\beta)}, 
	\end{align*}
	where $\lambda = \pre{\gamma}\ell$ and $\epsilon\in \bfC^{\times}$ is constant (which is a power of $2$).
\end{lemm}
\begin{proof}
	We divide the index set of the product on the left-hand side into two
	\[
		S^+\cap \pre{\gamma}w^{-1} S^- = \{b\in S^+\cap \pre{\gamma}w^{-1} S^-\;;\; \partial b\notin R^+\cap w^{-1}R^-\}\sqcup \{b\in S^+\cap \pre{\gamma}w^{-1} S^-\;;\; \partial b\in R^+\cap w^{-1}R^-\} 
	\]
	and treat the two sub-products separately.
	\begin{enumerate}
		\item[Step 1.] 
			We prove that $\omega_{\lambda}(b) = 0$ when $b\in S^+\cap \pre{\gamma}w^{-1} S^-$ and $\beta := \partial b\notin R^+\cap w^{-1}R^-$. \par
			 Write $b = \beta + k$ for $k\in \bfZ$. 
			We deduce
			\begin{align*}
				\omega_{\lambda}(b) = \omega_{\lambda}(\beta + k) = \til\omega_{\lambda_0}(y\beta + k+\langle \beta, \gamma\rangle),\quad \text{where $y\in W_R$ is such that $y\ell = \ell_0$}.
			\end{align*}
			Since $\pre{\gamma}w = w X^{w^{-1}\gamma - \gamma}$, the condition $b\in S^+\cap \pre{\gamma}w^{-1} S^-$ implies that $0\le k < \langle  w\beta - \beta, \gamma \rangle$. There are two cases: $\beta\in R^-$ or $\beta \in w^{-1} R^+$. 
			In the case where $\beta \in w^{-1}R^+$, since $k + \langle \beta, \gamma\rangle < \langle w\beta, \gamma\rangle \le -M$, the hypothesis~\eqref{equa:alphagamma} implies that $\omega_{\lambda}(b)  = 0$. 
			In the case where $\beta \in R^-$, we have $k + \langle \beta, \gamma\rangle \ge  \langle \beta, \gamma\rangle \ge M$, whence $\omega_{\lambda}(b) = 0$ as well.
		\item[Step 2.]
			We prove that
			\begin{align}\label{equa:prod2}
				\prod_{\substack{b\in S^+\cap\pre{\gamma}w^{-1} S^-\\ \partial b \in  R^+\cap w^{-1}R^-}}(-\partial b)^{\omega_\lambda(b)} = \epsilon\cdot\prod_{\beta\in R^+_\idv\cap w^{-1}R^-_\idv}(-\beta)^{\omegaa_{\ell}(\beta)}.
			\end{align}
			for some $\epsilon\in \bfC^{\times}$ which is a power of $2$. We rewrite the left-hand side according to $\partial b$:
			\begin{align}\label{equa:prod3}
				\prod_{\substack{b\in S^+\cap\pre{\gamma}w^{-1} S^-\\ \partial b \in  R^+\cap w^{-1}R^-}}(-\partial b)^{\omega_\lambda(b)} = \epsilon\cdot\prod_{\beta\in R^+_\idv\cap w^{-1}R^-_\idv}\prod_{\substack{b\in S^+\cap\pre{\gamma}w^{-1} S^-\\ \partial b \in \{\beta,2\beta\}}}(-\beta)^{\omega_\lambda(b)}
			\end{align}
			Let $\beta\in R^+_\idv\cap w^{-1}R^-_\idv$. Let $N:= \langle w\beta - \beta, \gamma\rangle$. It follows by the same arguments as Step 1 that $b = \beta + k\in S^+\cap \pre{\gamma} w^{-1} S^-$ for $0\le k \le N$.	For $k\ge N$, we obtain $k+\langle \beta, \gamma\rangle \ge \langle w\beta, \gamma\rangle \ge M$ , thus $\omega_{\lambda}(\beta + k) = \til \omega_{\lambda_0}(y\beta + k + \langle \beta, \gamma\rangle) = 0$ and hence
			\begin{align*}
				\sum_{\substack{b\in S^+\cap \pre{\gamma}w^{-1}S^- \\\partial b = \beta}} \omega_{\lambda}(b) = \sum_{k  = 0}^{N}\omega_{\lambda}(\beta + k) = \sum_{k  \in \bfN}\omega_{\lambda}(\beta + k) = \sum_{\substack{b\in S^+ \\\partial b = \beta}} \omega_{\lambda}(b) .  \\
			\end{align*}
			In the case where $2\beta\in R$, we obtain similarly
			\begin{align*}
				&\sum_{\substack{b\in S^+\cap \pre{\gamma}w^{-1}S^- \\\partial b = 2\beta}} \omega_{\lambda}(b) = \sum_{k= 0}^{N-1}\omega_{\lambda}(2\beta + (2k+1)) = \sum_{k \in \bfN}\omega_{\lambda}(2\beta + (2k+1)) =  \sum_{\substack{b\in S^+ \\\partial b = 2\beta}} \omega_{\lambda}(b) .  \\
			\end{align*}
			Hence
			\begin{align}\label{equa:somm2}
				\sum_{\substack{b\in S^+\cap \pre{\gamma}w^{-1}S^- \\\partial b \in \{\beta, 2\beta\}}} \omega_{\lambda}(b) = \omegaa_{\ell}(\beta).
			\end{align}
			The equation~\eqref{equa:prod2} follows from~\eqref{equa:prod3} and~\eqref{equa:somm2}.
	\end{enumerate}
	Combining the two steps, we obtain the product formula.
\end{proof}

\subsection{Preparatory lemmas}
Let $\gamma\in Q^{\vee}$ be an element satisfying~\eqref{equa:gamma}. Recall the notion of clans and generic clans from~\autoref{subsec:clans} and the fundamental alcove $\nu_0\subset E$. 
\begin{lemm}\label{lemm:clans}
	For $w\in W_R$, the alcove $w^{-1}X^{-\gamma}\nu_0$ is in a generic clan and every generic clan contains at least one such alcove. Moreover, for a different choice $\gamma'\in Q^{\vee}$, the alcoves $w^{-1}X^{-\gamma}\nu_0$ and $w^{-1}X^{-\gamma'}\nu_0$ are in the same clan. 
\end{lemm}
\begin{proof}
	Since the clans are connected components of the complement $E^{\omega}_{\circ}$ of the hyperplanes in $\frakD^{\omega} = \left\{ H_a\subset E\;;\; a\in S,\; \til\omega_{\lambda_0}(a) \ge 1 \right\}$, any two points $x, y\in E^{\omega}_{\circ}$ are in the same clan if $a(x)a(y) > 0$ for all $a\in S$ with $H_a\in \frakD^{\omega}$. Let $\frakC_w\subset E^{\omega}_{\circ}$ be the clan such that $w^{-1}X^{-\gamma}\nu_0\subset \frakC_w$. Take any point $x\in \nu_0$. Set $x_{w}(t) = w^{-1}(x - (1+t)\gamma)$ for $t\in \bfR_{\ge 0}$ so that in particular $x_w(0)\in \frakC_w$. Let $a\in S$ such that $H_a\in \frakD^{\omega}$. Suppose that $w\partial a\in R^{+}$ (resp. $w\partial a\in R^-$); then $\langle w\partial a, \gamma\rangle \ll 0$ (resp. $\langle w\partial a, \gamma\rangle \gg 0$), so we have
	\begin{align*}
			a(x_w(t)) = (w a)(x) - \langle w\partial a, (1+t)\gamma\rangle \gg 0,\quad \forall t\in \bfR_{\ge 0}
	\end{align*}
	(resp. $a(x_w(t)) \ll 0$). Hence $x_w(t)\in \frakC_w$ for all $t\ge 0$. Moreover, we see that the value $a(x_w(t))$ is unbounded when $t\to +\infty$.  Hence every affine root is unbounded on the $\frakC_w$, from which the genericity of $\frakC_w$.  \par
	Conversely, let $\frakC$ be a generic clan, consider the salient cone $\kappa$ defined in~\autoref{subsec:clans}. The genericity of $\frakC$ means that $\kappa$ is of full dimension $\dim V$. Let $\rmC_0\subset V$ be the fundamental Weyl chamber and let $w^{-1}\rmC_0\subset V$ be a Weyl chamber with $w\in W_{R}$ such that $\Int(\kappa)\cap w^{-1}\rmC_0\neq \emptyset$. It is obvious that $w^{-1}X^{-\gamma}\nu_0\subset \frakC$.
\end{proof}
Recall the element $\sigma_{\alpha}$ from~\eqref{equa:sigmaalpha}.
\begin{lemm}\label{lemm:sigmainA}
	We have $\sigma_{\alpha} \bfe(\pre{\gamma}\ell)\in \bfA^\omega$.
\end{lemm}
\begin{proof}
	Denote $\lambda = \pre{\gamma}\ell$. 
	Let $\pre{\gamma}{s_{\alpha}} = s_{a_l}\cdots s_{a_1}$ be any reduced decomposition and denote $\sigma'_{\alpha}\bfe(\lambda) = \tau^\omega_{a_l}\cdots \tau^\omega_{a_1}\bfe(\lambda)$.  
	Applying~\autoref{prop:formula}(ii), we see that 
	\begin{align}\label{equa:tausigma}
		\sigma'_{\alpha}\bfe(\lambda)\equiv s_{a_l}\cdots s_{a_1}\left(\prod_{c\in S^+\cap \pre{\gamma}s_\alpha S^-} (-\partial c)^{\omega_\lambda(c)}\right)\bfe(\lambda)\mod F_{\le l-1}\bfA^{-\infty}
	\end{align}
	and~\autoref{lemm:prod} yields
	\begin{align*}
			\prod_{c\in S^+\cap \pre{\gamma}s_\alpha S^-} (-\partial c)^{\omega_\lambda(c)} = \epsilon\cdot(-\alpha)^{\omega_{\ell}(\alpha)},\quad \epsilon\in \bfC^\times.
	\end{align*}
	Thus the right-hand side of~\eqref{equa:tausigma} is congruent to $\epsilon\sigma_{\alpha}\bfe(\lambda)$ modulo $F_{\le l-1}\bfA^{-\infty}$. Notice that $\sigma_{\alpha}\bfe(\lambda)\in \bfA^o$ by~\autorefitem{prop:Ao}{iii} and $\sigma'_{\alpha}\bfe(\lambda)\in \bfA^{\omega}\subset \bfA^{o}$. Hence by the compatibility of the filtrations by length~\autoref{coro:filt}, we have
	\begin{align*}
		(\sigma'_{\alpha} - \epsilon\sigma_{\alpha})\bfe(\lambda) \in  \bfe(\pre{\gamma}s_{\alpha}\lambda)\left(F_{\le l - 1}\bfA^{-\infty}\cap \bfA^{o}\right)\bfe(\lambda) = \bfe(\pre{\gamma}s_{\alpha}\lambda)\left(F_{\le l - 1}\bfA^o\right)\bfe(\lambda).
	\end{align*}
	We show that in fact $(\sigma'_{\alpha} - \epsilon\sigma_{\alpha})\bfe(\lambda) \in \bfA^\omega$.
	For any different choice $\gamma'$ satisfying~\eqref{equa:gamma}, \autoref{lemm:clans} implies that the intertwiner $\varphi_{\gamma,\gamma'}$ defined in~\autoref{theo:isomKH} satisfies $\varphi_{\gamma,\gamma'}\varphi_{\gamma',\gamma} = \bfe_{\gamma}$, $\varphi_{\gamma',\gamma}\varphi_{\gamma,\gamma'} = \bfe_{\gamma'}$ and 
	\[
		\varphi_{\gamma,\gamma'}\sigma_\alpha \bfe(\pre{\gamma'}\ell)\varphi_{\gamma',\gamma} = \sigma_\alpha \bfe(\pre{\gamma}\ell),\quad \varphi_{\gamma,\gamma'}\sigma'_\alpha \bfe(\pre{\gamma'}\ell)\varphi_{\gamma',\gamma} = \sigma'_\alpha \bfe(\pre{\gamma}\ell);
	\]
	thus the validity of the statement is independent of the choice of $\gamma$. We claim that if we choose $\gamma$ in such a way that $|\langle \alpha,  \gamma\rangle| \ll |\langle \beta,  \gamma\rangle|$ for all $\beta\in \Delta_0\setminus \{\alpha\}$, then there is an inequality of lengths
	\begin{align}\label{equa:longueur}
		l = \ell(\pre{\gamma}s_{\alpha}) \le \ell(\pre{\gamma}w),\quad \forall w\in W_R\setminus\{1\}.
	\end{align}
	We complete the proof provided~\eqref{equa:longueur}. Note that the stabilisers satisfy $\pre{\gamma}W_{\ell} = W_{\lambda}$. There are two cases to be discussed: 
	\begin{enumerate}
		\item
			If $s_{\alpha}\ell \neq \ell$, then by~\eqref{equa:longueur} we have $\ell(w) \ge l$ for all $w\in W_{S}$ such that $w\lambda = \pre{\gamma}s_{\alpha}\lambda$. It follows from~\autoref{theo:basis} that $\bfe(\pre{\gamma}s_\alpha\lambda)\left(F_{\le l - 1}\bfA^o\right)\bfe(\lambda) = 0$. Hence $\sigma_{\alpha}\bfe(\lambda) = \epsilon\sigma'_{\alpha}\bfe(\lambda)\in \bfA^\omega$.
		\item
			If $s_{\alpha}\ell = \ell$, then by~\eqref{equa:longueur} we have $\ell(w) \ge l$ for all $1\neq w\in W_{\lambda}$ and thus by~\autoref{theo:basis}, we see that $\bfe(\lambda)\left(F_{\le l - 1}\bfA^o\right)\bfe(\lambda) = \bfC[V]\bfe(\lambda)= \bfe(\lambda)\left(F_{\le l - 1}\bfA^\omega\right)\bfe(\lambda)$. Thus $(\sigma_{\alpha}- \epsilon^{-1}\sigma'_{\alpha})\bfe\left( \lambda \right)\in \bfA^\omega$ and consequently $\sigma_\alpha\bfe(\lambda)\in \bfA^\omega$. Hence the proof is completed.
	\end{enumerate}

	We prove~\eqref{equa:longueur}. Indeed by~\autoref{prop:lg},
	\begin{align*}
			l &= \ell(\pre{\gamma}s_{\alpha}) \le 1 + \ell(X^{-\langle \alpha, \gamma\rangle \alpha^{\vee}}) \le  1 +  |\langle 2\rho, \alpha^{\vee}\rangle\langle \alpha, \gamma\rangle|\ll |\langle \beta, \gamma\rangle|,\\
	\end{align*}
	while for any $w\in W_R \setminus\left\{ 1, s_{\alpha} \right\}$, there exists $\beta\in R^+_\idv \cap w^{-1}R^-_\idv$ with $\beta\neq \alpha$, so
	\begin{align*}
			\ell(\pre{\gamma}w) \ge \left|\langle \beta, w^{-1}\gamma - \gamma\rangle \right| - \ell(w) = \left|\langle w\beta - \beta,\gamma\rangle \right| - \ell(w) \ge  |\langle \beta, \gamma \rangle|- \ell(w)\gg l;
	\end{align*}
	here, the second-to-last inequality is due to~\eqref{equa:gamma}.
\end{proof}

\subsection{Basis theorem for generic clans}
Let $\gamma\in Q^{\vee}$ be an element satisfying~\eqref{equa:gamma}. Recall the idempotent of generic clans $\bfe_{\gamma}$ from~\eqref{equa:egamma} and the elements $\sigma_\alpha\bfe(\pre{\gamma}\ell)$ from~\eqref{equa:sigmaalpha}.
\begin{theo}\label{theo:basis-generic}
	The idempotent subalgebra $\bfe_{\gamma} \bfA^\omega\bfe_{\gamma}$ is generated by $\bfC[V]\bfe(\lambda)$ and $\sigma_{\alpha}\bfe(\lambda)$ for $\alpha\in \Delta_0$ and $\lambda\in \pre{\gamma}(W_R\ell_0)$. Moreover, if for any $w\in W_R$ we set $\sigma_{w} = \sigma_{\alpha_n}\cdots \sigma_{\alpha_1}$ by choosing any reduced expression $w = s_{\alpha_n} \cdots s_{\alpha_1}$, then there is a decomposition
	\begin{align*}
		\bfe_{\gamma} \bfA^\omega\bfe(\lambda) = \bigoplus_{\lambda\in \pre{\gamma}(W_R\ell_0) }\bigoplus_{w\in W_R} \bfC[V]\sigma_w \bfe(\lambda).
	\end{align*}
\end{theo}
\begin{proof}
	Let $\ell\in W_R\ell_0$ and $w\in W_R$. Denote $\lambda = \pre{\gamma}\ell$. Choose any reduced expressions $w = s_{\alpha_n}\cdots s_{\alpha_1}$ and $\pre{\gamma}w = s_{a_l} \cdots s_{a_1}$ for $\alpha_1, \cdots, \alpha_n\in \Delta_0$ and $a_1, \cdots, a_l\in \Delta$ and set
	\begin{align*}
			\sigma'_w \bfe(\lambda) &= \tau^{\omega}_{a_l}\cdots \tau^{\omega}_{a_1}\bfe(\lambda)\in \bfA^{\omega},\quad \sigma_w \bfe(\lambda) &= \sigma_{\alpha_n}\cdots \sigma_{\alpha_1} \bfe(\lambda).
	\end{align*}
	By~\autoref{lemm:sigmainA}, we see that $\sigma_w\bfe(\lambda)\in \bfA^{\omega}$.
	We claim that
	\begin{align}\label{equa:congruence}
		\sigma'_w\bfe(\lambda) \equiv \epsilon\sigma_w\bfe(\lambda) \mod F_{\le l-1}\bfA^\omega.
	\end{align}
	for some $\epsilon\in \bfC^{\times}$. Recall the rational function matrix algebra $\bfA^{-\infty} = \Frac \calZ\otimes_{\calZ}\bfA^{o}$. By~\autorefitem{prop:formula}{ii} and~\autoref{lemm:prod}, we have
	\begin{align*}
		&\sigma'_w\bfe(\lambda) \equiv s_{\partial a_l}\cdots s_{\partial a_1}\left( \prod_{b\in S^+\cap \pre{\gamma}w^{-1}S^-}\left( -\partial b \right)^{\omega_{\lambda}(b)} \right)\bfe(\lambda) \mod F_{\le l-1}\bfA^{-\infty}  \\
		&\equiv \epsilon\sigma_w\bfe(\lambda) \mod F_{\le l-1}\bfA^{-\infty}.
	\end{align*}
	for some $\epsilon\in \bfC^{\times}$. As $n\le l$, the above congruences yield $(\sigma'_w - \epsilon\sigma_w)\bfe(\lambda)\in \bfA^\omega\cap F_{\le l - 1}\bfA^{-\infty}$. By~\autoref{coro:filt}, we have $\bfA^\omega\cap F_{\le l - 1}\bfA^{-\infty} = F_{\le l-1}\bfA^\omega$, so the claim~\eqref{equa:congruence} is proven. \par
	According to~\autoref{theo:basis}, the family $\{\sigma'_{w}\bfe(\lambda)\}_{w\in W_R}$ form a basis for $\bfe_{\gamma}\bfA^\omega\bfe(\lambda)$. The decomposition of $\bfe_{\gamma}\bfA^\omega\bfe(\lambda)$ follows from the triangularity~\eqref{equa:congruence} of the transition matrix between the basis $\{\sigma'_{w}\bfe(\lambda)\}_{w\in W_R}$ and the family $\{\sigma_{w}\bfe(\lambda)\}_{w\in W_R}$.
\end{proof}

\subsection{Proof of~\autoref{theo:isomKH}}\label{subsec:demo-isomKH}
\begin{proof}
	We define an isomorphism of $\calZ$-modules $\Pol_{W_R \ell_0}\cong \bfe_\gamma\Pol_{W_S\lambda_0}$ straightforwardly by the identification:
	\begin{align*}
			\Pol_{\ell} = \bfC[V] = \Pol_{\pre{\gamma}\ell},\quad \ell\in W_R\ell_0.
	\end{align*}
	It yields a faithful representation of $\bfB^{\omegaa}$ on $\bfe_{\gamma}\Pol_{W_S \lambda_0}$, which by definition of $\bfB^{\omegaa}$ is described by the formula
	\begin{align*}
			f\bfe(\ell)\cdot g = f\bfe(\pre{\gamma}\ell) g,\quad \tau^{\omegaa}_{\alpha}\bfe(\ell)\cdot g = \sigma_{\alpha}\bfe(\pre{\gamma}\ell) g.
	\end{align*}
	By~\autoref{theo:basis-generic}, the image of $\bfB^{\omegaa}$ in $\End_{\calZ}\left( \bfe_{\gamma}\Pol_{W_S \lambda_0} \right)$ coincides with $\bfe_{\gamma}\bfA^{\omega}\bfe_\gamma$ and the map $\bfB^{\omegaa}\to \bfe_{\gamma}\bfA^{\omega}\bfe_{\gamma}$ must be an isomorphism since both sides are free $\bfC[V]$-modules of same rank. Notice that $\deg \tau^{\omegaa}_\alpha\bfe(\ell) = \omegaa_{\ell}(\alpha) = \deg \sigma_{\alpha}\bfe(\pre{\gamma}\ell)$. Hence the map $i_{\gamma}$ is an isomorphism of graded $\calZ$-algebras.\par
For any other choice $\gamma'$, since by~\autoref{lemm:clans}, $w^{-1}X^{\gamma}$ and $w^{-1}X^{\gamma'}$ lie in the same generic clan for each $w\in W_R$, by~\autoref{prop:entrelacement}, the intertwiner $\varphi_{\gamma',\gamma}$ yields isomorphisms of $\bfA^{\omega}$-modules $\bfA^{\omega}\bfe_{\gamma'} \cong \bfA^{\omega}\bfe_{\gamma}$ by right multiplication and hence isomorphisms of algebras
	\begin{align*}
		\bfe_{\gamma'}\bfA^{\omega}\bfe_{\gamma'}\cong \End_{\bfA^{\omega}}\left(\bfA^{\omega}\bfe_{\gamma'}\right)\cong \End_{\bfA^{\omega}}\left(\bfA^{\omega}\bfe_{\gamma}\right)\cong \bfe_{\gamma}\bfA^{\omega}\bfe_{\gamma}.
	\end{align*}
	The factorisation $i_{\gamma} = \varphi_{X^{\gamma'-\gamma}}\circ i_{\gamma'}$ follows from the observation that $\partial(X^{\gamma' - \gamma}) = 1\in W_R$. \par
\end{proof}

\subsection{The functor \texorpdfstring{$\bfV$}{V}}\label{subsec:V}
Choose a $\gamma\in Q^{\vee}$ satisfying~\eqref{equa:gamma} as in~\autoref{subsec:formuleomega}. With~\autoref{theo:isomKH}, we can make the following definition:
\begin{defi}\label{defi:KZalg}
	The Knizhnik--Zamolodchikov (KZ) functor $\bfV$\index[ch2]{V@$\bfV$} is defined by
	\begin{align*}
		\bfV:\bfA^\omega\gmod&\to \bfe_{\gamma}\bfA^\omega\bfe_{\gamma}\gmod\autorightarrow{$i_{\gamma}^*$}{$\cong$}\bfB^{\omegaa}\gmod \\
		M &\mapsto \bfe_{\gamma}M.
	\end{align*}
\end{defi}
By the second assertion of~\autoref{theo:isomKH}, the definition of $\bfV$ is independent of the choice of $\gamma$ up to canonical isomorphism (provided by the intertwiner $\varphi_{\gamma,\gamma'}$).  \par
Since $\bfV$ is defined as an idempotent truncation, it admits left and right adjoint functors
\begin{align*}
	\bfV^{\top}: N\mapsto \bigoplus_{\lambda\in W_S \lambda_0}\gHom_{\bfB^{\omegaa}}\left(\bfe_{\gamma}\bfA^\omega \bfe(\lambda), N  \right)\quad\text{and}\quad
	\pre{\top}\bfV: N\mapsto \bfA^\omega\bfe_{\gamma}\otimes_{\bfB^{\omegaa}}N
\end{align*}
and $\bfV$ is a quotient functor in the sense that the adjoint counit $\bfV\circ \bfV^{\top}\to \id_{\bfB^{\omegaa}}$ is an isomorphism.

\subsection{Support characterisation of \texorpdfstring{$\bfV$}{V}}\label{subsec:suppV}
For $M\in \bfA^{\omega}\gmod$, define the following subset of $E$:
\begin{align*}
	\Spec_E M = \left\{ \lambda\in W_S \lambda_0\;;\; \bfe(\lambda)M\neq 0\right\}.
\end{align*}
For each alcove $\nu\subset E$, there is a unique $w\in W_S$ such that $\nu = w^{-1}\nu_0$; we denote $\lambda_{\nu} = w\lambda_0$. Recall the Gelfand--Kirillov dimension $\dim_{\GK}M$ and the support $\supp M$ from~\autoref{subsec:supp}.
\begin{theo}\label{theo:kerV}
	Let $M\in \bfA^{\omega}\gmod_0$. The following conditions are equivalent:
	\begin{enumerate}
		\item\label{theo:kerV-i}
			$\bfV M = 0$;
		\item\label{theo:kerV-ii}
			for every alcove $\nu$ lying in a generic clan, we have $\bfe(\lambda_{\nu}) M= 0$;
		\item\label{theo:kerV-iii}
			the set $\Spec_E M$ is contained in a finite union of (not-necessarily root) affine hyperplanes of $E$;
		\item\label{theo:kerV-iv}
			$\dim_{\mathrm{GK}} M \le \rk R-1$;
		\item\label{theo:kerV-v}
			$\supp M \neq \Spec \bfC Q^{\vee}_{\dom}$.
	\end{enumerate}
\end{theo}
\begin{proof}
	Since every object of the category $\bfA^{\omega}\gmod$ is of finite length and all the conditions \ref{theo:kerV-i}--\ref{theo:kerV-v} are stable under extensions, we may suppose that $M$ is simple. \par

	\ref{theo:kerV-i} $\Leftrightarrow$ \ref{theo:kerV-ii} follows from the definition $\bfV M = \bfe_{\gamma}M$ and the invariance of dimension of $\bfe(\lambda)M$ for $\lambda$'s in the same clan~\autoref{coro:gdim}. \par

	We prove \ref{theo:kerV-ii} $\Rightarrow$ \ref{theo:kerV-iii}.  By the finiteness of the clan decomposition, it suffices to show that for each non-generic clan $\frakC$, the set $\left\{ \lambda_\nu\;;\; \nu\subseteq \frakC \right\}$ lies in a finite union of affine hyperplanes of $E$. By the non-genericity of $\frakC$, there exists $\alpha\in R$ which is bounded on $\frakC$. Let $\Lambda = \ker\alpha\cap Q^{\vee}$. Notice that $Q^{\vee}$ is a free $\bfZ$-module of rank $\rk R - 1$. Let $\frakA_{\frakC}$ be the set of alcoves contained in $\frakC$. For $\nu, \nu'\in \frakA_{\frakC}$, we write $\nu\sim_{\Lambda} \nu'$ if there exists $\mu\in \Lambda$ such that $\nu + \mu = \nu'$. For any $\nu \in \frakA_{\frakC}$, since $X^{\mu}\lambda_{\nu} = \lambda_{\nu} + \mu$, the set $\left\{ \lambda_{\nu'}\;;\; \nu'\sim_\Lambda \nu \right\}$ is contained in the hyperplane $w\left(\lambda_0 + \Lambda_{\bfR}\right)$ for any $w\in W_S$ such that $\nu = w^{-1}\nu_0$. Since $\alpha$ is bounded on $\frakC$, the quotient $\frakA_{\frakC} / \sim_{\Lambda}$ is a finite set and thus the set 
	\begin{align*}
		\left\{ \lambda_{\nu}\;;\; \nu\subset \frakC \right\} \subset \bigcup_{\nu\in \frakA_{\frakC} / \sim_\Lambda}\left\{ \lambda_{\nu'}\;;\; \nu'\sim_\Lambda \nu \right\}
	\end{align*}
	is contained in a finite union of hyperplanes, whence~\ref{theo:kerV-iii}. \par

	We prove \ref{theo:kerV-iii} $\Rightarrow$ \ref{theo:kerV-iv}. Suppose that $\Spec_{E} M$ is contained in a finite number of hyperplanes. Choose any $\lambda_1\in \Spec_E M$. Let $r = \rk R = \dim E$. Via the identification $E\cong V$ induced by $\Delta_0\subset \Delta$, we view $E$ as an euclidean vector space.  Since
	\begin{align*}
		\Spec_{E} M\subset \bigcup_{w\in W_R}(w\lambda_1 + Q^{\vee})
	\end{align*}
	is contained in a finite union of the intersection of lattices and hyperplanes, we have
	\begin{align*}
		\lim_{n\to \infty} \frac{\#\left\{\lambda\in \Spec_{E} M \;;\; \|\lambda\| < n\right\}}{n^{r-1 + \varepsilon}} = 0,\quad \forall \varepsilon > 0.
	\end{align*}\par
	For every affine simple root $a\in \Delta$, we have $\tau^{\omega}_a\bfe(\lambda)M\subseteq\bfe(s_a\lambda) M$. Moreover, we have $\|s_a\lambda\| \le \|\lambda\| + \delta$ for some constant $\delta$ which depends only on the affine root system $(E, S)$. It follows that if we define for $t\in \bfR_{\ge 0}$ the subspace
	\begin{align*}
		M_{\le t} = \sum_{\substack{\lambda\in \Spec_{E}M \\ \|\lambda \| \le t}} \bfe(\lambda) M,
	\end{align*}
	then $\tau^{\omega}_a M_{\le t} \subset M_{\le t + \delta}$, so $F_{\le 1}\bfA^{\omega}M_{\le t}\subset M_{\le t + \delta}$. By induction on $n\in \bfN$, we see that $\left(F_{\le n}\bfA^{\omega}\right)M_{\le t}\subset M_{\le t + n\delta}$. Since there is only a finite number of clans and since the dimension of $\bfe(w\lambda_0)M_{\lambda}$ for $w^{-1}\nu_0$ in a fixed clan is constant by~\autoref{coro:gdim}, the set $\{\dim \bfe(\lambda)M\;;\; \lambda\in W_S \lambda_0\}$ is bounded. Hence for any finite-dimensional subspaces $L\subset M$, we have
	\begin{align}\label{equa:estime}
		\lim_{n\to \infty} \frac{\dim\left(F_{\le n}\bfA^{\omega}\cdot L\right)}{n^{r - 1 + \varepsilon}} = 0,\quad \forall \varepsilon > 0.
	\end{align}
	Indeed, let $t_0\in \bfR$ be such that $L\subset M_{\le t_0}$, then  $\dim\left(F_{\le n}\bfA^{\omega}\cdot L\right)\le \dim M_{\le t_0 + n\delta} = o\left( n^{r - 1 + \varepsilon}  \right)$. The estimate~\eqref{equa:estime} implies~\ref{theo:kerV-iv}. The equivalence \ref{theo:kerV-iv} $\Leftrightarrow$ \ref{theo:kerV-v} results from~\autoref{prop:dimGK}.\par
	We prove $\neg$\ref{theo:kerV-ii} $\Rightarrow$ $\neg$\ref{theo:kerV-iv}. Suppose there exists a generic clan $\frakC$ and an alcove $\nu\subset \frakC$ such that $\bfe\left( \lambda_{\nu} \right)M \neq 0$. Let $\kappa\subset V$ be the salient cone of $\frakC$ (\cfauto{subsec:clans}). For any $\mu\in \kappa \cap Q^{\vee}$, we have $X^{\mu}\nu\in \frakC$ and by~\autoref{prop:entrelacement}, $\bfe(X^{-\mu}\lambda_{\nu})M\cong \bfe(\lambda_\nu)M\neq 0$. It follows that
	\begin{align*}
		\dim \left(F_{\le n} \bfA^{\omega}\right)\left(\bfe(\lambda_{\nu})M\right) \ge \dim \sum_{\substack{\mu\in \kappa\cap Q^{\vee} \\ \ell(X^{\mu}) \le n }} \bfe(X^{-\mu}\lambda_{\nu})M = \#\{\mu\in \kappa\cap Q^{\vee} \;;\; \ell(X^{\mu}) \le n\} \dim \bfe(\lambda_\nu) M.
	\end{align*}
	By the genericity of $\frakC$, the salient cone $\kappa$ contains an open subset of $V$, so its intersection with a full-ranked lattice $Q^{\vee}$ satisfies 
	\begin{align*}
			\lim_{n\to\infty}\frac{\#\{\mu\in \kappa\cap Q^{\vee} \;;\; \ell(X^{\mu}) \le n\}}{n^r} = c,\quad c > 0.
	\end{align*}
	Hence
	\begin{align*}
		\dim_{\GK} M \ge \lim_{n\to \infty}\frac{\log\dim (F_{\le n} \bfA^{\omega}) \bfe(\lambda_\nu)M }{\log n} \ge \lim_{n\to \infty}\frac{\log  cn^r }{\log n} = r,
	\end{align*}
	whence~\ref{theo:kerV-iv} is not satisfied. 
\end{proof}

\subsection{Double centraliser property}\label{subsec:bicomm}
Recall the parabolic subalgebra $\bfA^{\omega}_{R, \lambda_1}$ from~\autoref{subsec:indres}. 
\begin{lemm}\label{lemm:non-vanishing}
	Let $\lambda_1\in W_S\lambda_0$, $N\in \bfA^{\omega}_{R, \lambda_1}\gmod$ and $L\in \bfA^{\omega}\gmod$. Suppose that $\bfV L = 0$, then $\gHom\left( L, \ind_{R, \lambda_1}^{S}N \right) = 0$.
\end{lemm}
\begin{proof}
	It follows from~\hyperref[theo:kerV]{\autoref*{theo:kerV}~(i)$\Rightarrow$(v)} and~\autoref{prop:indtf}.
\end{proof}
\begin{rema}
	We shall establish in~\autoref{theo:bicommutante} the double centraliser property for the functor $\bfV$. The strategy is close to the case of rational Cherednik algebras in~\cite[5.3]{GGOR03}: the first step consists of showing that ``induced modules'' are torsion-free for the KZ functor. In the case of the dDAHA $\bbH$ discussed in~\autoref{partI}, the parabolic subalgebra $\bfA^{\omega}_{R, \lambda_1}$ plays the r\^ole of be graded affine Hecke subalgebra $\ubar\bbH = \bfC W_R \otimes \bfC[E]\subset \bbH$, whereas $\bfB^\omegaa$ plays the r\^ole of the affine Hecke algebra $\bbK$. In this sense, \autoref{lemm:non-vanishing} is an analogue of the first step in the proof of~{\itshape loc. cit.}
\end{rema}
Let $\left(\bfA^{\omega}/\frakm_{\calZ}\right)\gmod$ be the full subcategory of $\bfA^{\omega}\gmod$ consisting of objects $M$ such that $\frakm_{\calZ}M = 0$. The inclusion $\left(\bfA^{\omega}/\frakm_{\calZ}\right)\gmod\hookrightarrow\bfA^{\omega}\gmod$ has a left adjoint functor $\relbar\otimes_{\calZ}\bfC$, which is right exact. We denote by $\relbar\otimes_{\calZ}^{\rmL}\bfC$ its derived functor.
The next lemma is the method of lifting faithfulness borrowed from~\cite[4.42]{rouquier08}. 
\begin{lemm}\label{lemm:Ext1}
	Let $M\in \bfA^{\omega}\gmod$ be an object satisfying the following properties:
	\begin{enumerate}
		\item
			$M$ is free over the centre $\calZ$;
		\item
			there exists $\lambda_1\in W_S\lambda_0$ and $N\in \bfA^{\omega}_{R,\lambda_1}\gmod_0$ such that $M / \frakm_{\calZ} M\cong \ind^{S}_{R, \lambda_1} N$.
	\end{enumerate}
	Then for any $L\in \bfA^{\omega}\gmod$ such that $\bfV L = 0$, we have $\gHom(L, M) = 0$ and $\gExt^1(L, M) = 0$. 
\end{lemm}
\begin{proof}
	We suppose that $M\neq 0$. Let $K = \rmR \gHom(L, M)$ be in the derived category $\rmD^{+}(\calZ\gMod)$. We suppose that $K$ is a minimal projective resolution. Since 
	\begin{align*}
		K\otimes_{\calZ}\bfC\cong \rmR\gHom\left( L\otimes^{\rmL}_{\calZ}\bfC, M\otimes^{\rmL}_{\calZ}\bfC\right) \cong\rmR\gHom\left( L\otimes^{\rmL}_{\calZ}\bfC,  M / \frakm_{\calZ} M\right)
	\end{align*}
	by the flatness of $M$ over $\calZ$, we have $K\otimes_{\calZ} \bfC\in \rmD^{\ge 0}\left( \bfC \right)$. By the second assumption and~\autoref{lemm:non-vanishing}, we have
	\begin{align*}
		\rmH^{0}(K\otimes_{\calZ}\bfC) = \Hom\left( L\otimes_{\calZ}\bfC, M/\frakm_{\calZ} M \right) = 0.
	\end{align*}
	Consequently $\rmH^{\le 0}(K) = 0$ by Nakayama's lemma. \par
	Suppose that $\rmH^1(K) \neq 0$. Since the localisation $\Frac \calZ\otimes_{\calZ}M$ is a weight module over $\bfA^{-\infty}$, which is semisimple, $\rmH^1(K)$ must be a torsion module over $\calZ$ so $K^0 \neq 0$. However, the minimality of $K$ would imply $\rmH^0(K\otimes_{\calZ}\bfC) \neq 0$, contradiction. Hence $\rmH^{\le 1}(K) = 0$ and so $\gHom(L, M) = 0$ and $\gExt^1(L, M) = 0$ as asserted. \par
\end{proof}

\begin{lemm}\label{prop:fidelite}
	Let $M\in \bfA^{\omega}\gmod$ be an object satisfying~\autoref{lemm:Ext1}. Then the adjoint unit yields an isomorphism $M\cong (\bfV^{\top}\circ\bfV) M$.
\end{lemm}
\begin{proof}
	Set $X = \Cone \left( M\to (\rmR\bfV^{\top}\circ\bfV) M \right)\in \rmD^{+}(\bfA^{\omega}\gmod)$, so that there is a distinguished triangle
	\begin{align}\label{equa:dist}
		M\to (\rmR \bfV^{\top}\circ\bfV) M\to X\to M[1].
	\end{align}
	By the adjunction and the exactness of $\bfV$, we have $\bfV X \cong \Cone(\bfV M\to (\bfV\circ\rmR\bfV^{\top}\circ\bfV) M ) = 0$ and hence 
	\begin{align*}
		\bfV\rmH^k(X) \cong \rmH^k(\bfV X) = 0,\quad k\in \bfZ.
	\end{align*}
	Applying~\autoref{lemm:Ext1} with $L = \rmH^0(X)$ and $L = \rmH^{-1}(X)$, we deduce
	\begin{align*}
		\gHom\left( \rmH^0(X), M \right) = 0,\quad \gHom\left( \rmH^0(X)[-1], M \right) = \gExt^1\left( \rmH^0(X), M \right) = 0  \\
		\gHom\left( \rmH^{-1}(X), M \right) = 0,
	\end{align*}
	whence
	\begin{align}\label{equa:VX=0}
		\gHom(\tau_{\le 0} X, M) = 0,\quad \gHom(\tau_{\le 0} X, M[1]) = \gHom(\tau_{\le 0} X[-1], M) = 0.
	\end{align} \par
	Applying $\rmR\gHom\left(\tau_{\le 0}X, \relbar\right)$ to the distinguished triangle~\eqref{equa:dist}, we obtain the long exact sequence
	\begin{align*}
		\gHom\left( \tau_{\le 0}X, M \right)\to \gHom\left( \tau_{\le 0}X, (\rmR\bfV^{\top}\circ\bfV) M \right)\to \gHom\left( \tau_{\le 0}X, X \right)\to\gHom\left( \tau_{\le 0}X, M[1] \right).
	\end{align*}
	By~\eqref{equa:VX=0}, the first and the last term of the sequence vanish. Hence,
	\begin{align*}
		\gHom\left(\tau_{\le 0} X, X  \right)\cong \gHom\left(\tau_{\le 0} X, (\rmR\bfV^{\top}\circ\bfV) M  \right) \cong\gHom\left(\tau_{\le 0}\bfV X,\bfV M  \right) = 0,
	\end{align*}
	which implies that $\tau_{\le 0} X = 0$. Applying $\rmH^0$ to the distinguished triangle~\eqref{equa:dist}, we deduce that the adjunction unit $M\to (\bfV^{\top}\circ\bfV) M$ is an isomorphism.
\end{proof}

\begin{theo}[Double centraliser property\footnote{Let $A$ and $B$ be unital associative rings. Usually, one says that an $(A, B)$-bimodule $P$ satisfies the double centraliser property if the structural maps $A\to \End_{B^{\op}}(P)$ and $B\to \End_A(P)^\op$ are isomorphisms. The above theorem provides a graded, non-unital version of this property for the $(\bfA^{\omega}, \bfB^{\omegaa})$-bimodule $\bfA^{\omega}\bfe_{\gamma}$.}
	]\label{theo:bicommutante}
	The canonical map
	\begin{align*}
		\bfA^{\omega} \to \bigoplus_{\lambda,\lambda'\in W_S \lambda_0}\gHom_{\bfB^{\omegaa}}\left(\bfV\bfA^{\omega}\bfe\left( \lambda \right), \bfV\bfA^{\omega} \bfe\left( \lambda' \right) \right)
	\end{align*}
	is an isomorphism.
\end{theo}
\begin{proof}
	Observe that for each $\lambda\in W_S \lambda_0$, the module $\bfA^{\omega}\bfe(\lambda)\in \bfA^{\omega}\gmod$ satisfies the conditions of~\autoref{lemm:Ext1}. Indeed, $\bfA^{\omega}\bfe(\lambda)$ is flat over $\calZ$ by~\autoref{theo:basis}. For the second condition, we have $\bfA^{\omega}\bfe(\lambda)\cong \ind^{S}_{R, \lambda_1} \bfA^{\omega}_{R, \lambda}\bfe(\lambda)$, so $\bfA^{\omega}\bfe(\lambda)/\frakm_{\calZ}\cong \ind^{S}_{R, \lambda_1}\left( \bfA^{\omega}_{R, \lambda}\bfe(\lambda)/\frakm_{\calZ}\right)$. Applying~\autoref{prop:fidelite}, we obtain
	\begin{align*}
		\bfA^{\omega} \cong \bigoplus_{\lambda,\lambda'\in W_S \lambda_0}\gHom_{\bfA^{\omega}}\left(\bfA^{\omega}\bfe\left( \lambda \right), \bfA^{\omega}\bfe\left( \lambda' \right)  \right)\xrightarrow{\cong}\bigoplus_{\lambda,\lambda'\in W_S \lambda_0}\gHom_{\bfA^{\omega}}\left(\bfV\bfA^{\omega}\bfe\left( \lambda \right), \bfV\bfA^{\omega}\bfe\left( \lambda' \right)  \right).
	\end{align*}
\end{proof}

\subsection{Categorical characterisation of \texorpdfstring{$\bfV$}{V}}\label{subsec:catV}
We shall exploit the Frobenius structure on $\bfB^{\omegaa}$ introduced in~\autoref{lemm:Frob}. Consider the anti-involution $\bfA^{\omega}\cong \left( \bfA^{\omega} \right)^{\op}$ which fixes pointwise $\bfC[V]\bfe(\lambda)$ for $\lambda\in W_S$ and sends $\tau^{\omega}_a\bfe(\lambda)\mapsto \tau^{\omega}_a\bfe(s_a\lambda)$. The duality 
\begin{align}\label{equa:dual1}
	M\mapsto M^* := \bigoplus_{\lambda\in W_S \lambda_0}\Hom_{\bfC}(\bfe(\lambda)M, \bfC)
\end{align}
yields an equivalence 
\[
	\bfA^{\omega}\gmod_0\cong ((\bfA^{\omega})^{\op}\gmod_0)^\op\cong (\bfA^{\omega}\gmod_0)^\op.
\]
Similarly, the anti-involution $\bfB^{\omegaa}\cong(\bfB^{\omegaa})^\op$ given by $\tau^{\omegaa}_\alpha\bfe(\ell)\mapsto \tau^{\omegaa}_\alpha\bfe(s_\alpha\ell)$ yields
 \[
	 \bfB^{\omegaa}\gmod_0\cong(\bfB^{\omegaa}\gmod_0)^\op. 
 \]
Denote $\ba\bfA^{\omega} = \bfA^{\omega} / \bfA^{\omega}\frakm_{\calZ}$ and $\ba\bfB^{\omegaa} = \bfB^{\omegaa} / \bfB^{\omegaa}\frakm_{\calZ}$.  Notice that the pairing 
	\[
		\bfe_{\gamma}\ba\bfA^{\omega}\times \ba\bfA^{\omega}\bfe_{\gamma}\xrightarrow{(a, b)\mapsto ab} \bfe_{\gamma}\ba\bfA^{\omega}\bfe_{\gamma} = \ba\bfB^{\omegaa}
	\]
	composed with the Frobenius form $\ba\bfB^{\omegaa}\xrightarrow{\tr}\calZ / \frakm_{\calZ} = \bfC$ yields an isomorphism $(\ba\bfA^{\omega}\bfe_{\gamma})^* \cong\ba\bfA^{\omega}\bfe_{\gamma}$.
\begin{lemm}\label{lemm:adj}
	There are canonical isomorphisms $\pre{\top}\bfV\ba\bfB^{\omegaa}\cong \ba\bfA^{\omega}\bfe_{\gamma} \cong \bfV^{\top}\ba\bfB^{\omegaa}$.
\end{lemm}
\begin{proof}
	The first isomorphism is obvious: $\pre{\top}\bfV\ba\bfB^{\omegaa} = \bfA^{\omega}\bfe_{\gamma}\otimes_{\bfB^{\omegaa}}\ba\bfB^{\omegaa} = \ba\bfA^{\omega}\bfe_{\gamma}$. \par
	Observe that $(\ba\bfA^{\omega}\bfe_{\gamma})^*\cong \ba\bfA^{\omega}\bfe_{\gamma}$ implies $\bfV M^* \cong (\bfV M)^*$ for $M\in \ba\bfA^{\omega}\gmod_0$ and hence $\bfV^{\top} N^* \cong (\pre{\top}\bfV N)^*$ for $N\in \ba\bfB^{\omegaa}\gmod_0$. Therefore
	\[
		\bfV^{\top}\ba\bfB^{\omegaa} \cong (\pre{\top}\bfV(\ba\bfB^{\omegaa})^*)^* \cong (\pre{\top}\bfV\ba\bfB^{\omegaa})^* \cong (\ba\bfA^{\omega}\bfe_{\gamma})^*\cong \ba\bfA^{\omega}\bfe_{\gamma}.
	\]
\end{proof}
\begin{theo}\label{prop:catV}
	Let $L\in \bfA^{\omega}\gmod_0$ be a simple object. Then the following conditions are equivalent:
	\begin{enumerate}
		\item\label{prop:catV-i}
			$\bfV L \neq 0$;
		\item\label{prop:catV-ii}
			the injective hull of $L$ in the subcategory $\ba\bfA^{\omega}\gmod$ is projective;
		\item\label{prop:catV-iii}
			the projective cover of $L$ in the subcategory $\ba\bfA^{\omega}\gmod$ is injective.
	\end{enumerate}
\end{theo}
\begin{proof}
	Since $\bfV^{\top}$ preserves injective objects, we see that by \autoref{lemm:adj}, $\ba\bfA^{\omega}\bfe_{\gamma}$ is injective-projective in $\ba\bfA^{\omega}\gmod$. \par
	We prove~\ref{prop:catV-i} $\Leftrightarrow$ \ref{prop:catV-ii}. Let $L\in \bfA^{\omega}\gmod_0$ be any simple object. 
	If $\bfV L = 0$, then by~\autoref{lemm:non-vanishing}, we have $\gHom(L, \ba\bfA^{\omega}) = 0$; hence \ref{prop:catV-ii} fails for $L$. If $\bfV L\neq 0$, since $\bfV$ is a quotient functor, $\bfV L\in \bfB^{\omegaa}\gmod_0$ must be simple. We have $\frakm_{\calZ} \bfV L = 0$, so we may view $L$ as a $\ba\bfB^{\omegaa}$-module. By the self-injectivity of $\ba\bfB^{\omegaa}$, there exists a non-zero map $\iota: \bfV L\hookrightarrow \ba\bfB^{\omegaa}$ and the adjunction yields an injective map $L\to \bfV^{\top}\ba\bfB^{\omega}\cong\ba\bfA^{\omega}\bfe_{\gamma}$, whence~\ref{prop:catV-ii} holds for $L$.  \par
	Finally, since the duality~\eqref{equa:dual1} exchanges the projective and injective objects in $\ba\bfA^{\omega}\gmod$ and preserves $\ker \bfV$, we deduce
	\begin{align*}
		\text{\ref{prop:catV-iii} for $L$} \Leftrightarrow \text{\ref{prop:catV-ii} for $L^*$}\Leftrightarrow \text{\ref{prop:catV-i} for $L^*$}\Leftrightarrow\text{\ref{prop:catV-i} for $L$}.
	\end{align*}
\end{proof}
\begin{exam}
	Resume to the setting of examples~\autorefitem{exam:Aomega}{iii},~\autoref{exam:clan} and~\autoref{exam:idemp}. We have $\bfA^{\omega}\bfe_{\gamma} = P_+ \oplus P_-$ so $\bfV L_0 = 0$, while $\bfV L_+\neq 0$ and $\bfV L_-\ne 0$ are simple objects in $\bfB^{\omega}\gmod$. In regard of~\autoref{theo:kerV}, we have $\dim_{\GK} L_+ = \dim_{\GK} L_- = 1$ while $\dim_{\GK} L_0 = 0$. The cosocle filtration of $\bfV P_+$, $\bfV P_0$ and $\bfV P_-$ are described by the following:
	\begin{align*}
		\bfV P_+ = \begin{bNiceMatrix}[small]\bfV L_+ & \\ \bfV L_+\langle -2\rangle& \bfV L_-\langle -2\rangle\\ \bfV L_+\langle -4\rangle& \bfV L_-\langle -4\rangle \\ \bfV L_+\langle -6\rangle& \bfV L_-\langle -6\rangle\\ \vdots& \end{bNiceMatrix},\;
		\bfV P_0 = \begin{bNiceMatrix}[small] \bfV L_+\langle -1\rangle& \bfV  L_-\langle -1\rangle \\\bfV L_+\langle -3\rangle& \bfV L_-\langle -3\rangle \\\bfV L_+\langle -5\rangle& \bfV L_-\langle -5\rangle \\ \vdots& \end{bNiceMatrix},\;
		\bfV P_- = \begin{bNiceMatrix}[small]& \bfV L_- \\ \bfV L_+\langle -2\rangle& \bfV L_-\langle -2\rangle\\ \bfV L_+\langle -4\rangle& \bfV L_-\langle -4\rangle \\ \bfV L_+\langle -6\rangle& \bfV L_-\langle -6\rangle\\ \vdots& \end{bNiceMatrix}.
	\end{align*}
	From this description it is obvious that the functor $\bfV$ is fully faithful on the projective objects, so $\bfV$ satisfies the double centraliser property~\autoref{theo:bicommutante}. \par
	Consider the quotients 
	\begin{align*}
		P_+/\frakm_{\calZ} = \begin{bNiceMatrix}L_+  \\ L_0\langle -1\rangle \\  L_-\langle -2\rangle  \end{bNiceMatrix},\;
		P_0/\frakm_{\calZ} = \begin{bNiceMatrix} L_0 \\ L_+\langle -1\rangle\; L_-\langle -1\rangle  \end{bNiceMatrix},\;
		P_-/\frakm_{\calZ} = \begin{bNiceMatrix} L_- \\ L_0\langle -1\rangle \\ L_+\langle -2\rangle\end{bNiceMatrix}.
	\end{align*}
	It follows that $P_+ / \frakm_{\calZ}$ (resp. $P_- / \frakm_{\calZ}$) is the injective hull of $L_-\langle -2\rangle$ (resp. $L_+\langle -2\rangle$) in the category $(\bfA^{\omega} / \frakm_{\calZ})\gmod$ while $P_0 / \frakm_{\calZ}$ is not injective. Hence $L_+$ and $L_-$ satisfy the equivalent conditions of~\autoref{prop:catV}.
\end{exam}

\appendix

\section{Category of pro-objects}\label{sec:pro}
\subsection{}
Let $\calA$ be an abelian category. We denote by $\Pro(\calA)$ and $\Ind(\calA)$ the category of pro-objects and ind-objects. The basic reference for these is~\cite[8.6]{KS06}. All the results below are stated for $\Pro(\calA)$ while they all have a dual version for $\Ind(\calA)$. 
An object of $\Pro(\calA)$ is a filtered ``projective limit'' of objects of $\calA$. If 
\begin{align*}
	M^{(i)} = \prolim{j\in \calI^{(i)}} M^{(i)}_j, \quad M^{(i)}_j\in \calA, \quad i \in \{1,2\}
\end{align*}
are two objects of $\Pro(\calA)$, where $\calI^{(i)}$'s are filtrant diagram categories and $M^{(i)}:\calI^{(i)\op}\to \calA$'s are functors, then the Hom-space between them is given by
\begin{align}\label{equa:hom-pro}
	\Hom_{\Pro(\calA)}\left(M^{(1)}, M^{(2)}  \right) = \varprojlim_{j\in\calI^{(2)}}\varinjlim_{i\in\calI^{(1)}}\Hom_{\calA}\left(M^{(1)}_i, M^{(2)}_j  \right).
\end{align}

\subsection{}
For every $M\in \Pro(\calA)$, let $\calA^M$ denotes the category whose objects are pairs $(M', a)$ where $M'\in \calA$ and $a\in \Hom_{\Pro(\calA)}(M, M')$, and whose morphisms are given by
\begin{align*}
	\Hom_{\calA^M}\left( (M_1, a_1),(M_2, a_2) \right) = \left\{ b\in \Hom_{\calA}(M_1, M_2)\;;\; a_2 = b\circ a_1 \right\}.
\end{align*}
Every object $M\in\Pro(\calA)$ can be expressed as the following filtered limit:
\begin{align}\label{equa:limAM}
	M \cong \prolim{(M', a)\in \calA^M} M'.
\end{align} 
\par
Let $\calA^M_{\epi}\subset\calA^M$ be the full subcategory whose objects are the pairs $(M', q)$ with $q$ being an epimorphism. 

\begin{prop}\label{lemm:pro}
	Let $\calA$ be an artinian abelian category. Then the following statements hold:
	\begin{enumerate}
		\item
			$\calA$ is a Serre subcategory of $\Pro(\calA)$.
		\item
			Every object $M\in \Pro(\calA)$ can be written as the following filtered projective limit
			\begin{align*}
					M \cong \prolim{(M', a)\in \calA^M_{\epi}}M'.
			\end{align*}
		\item
			$\calA$ is the full subcategory of artinian objects in $\Pro(\calA)$.
		\item
			If $\varphi:N\to M$ is a morphism in $\Pro(\calA)$ such that for every $(M',q)$ in $\calA^M_{\epi}$, the composite $q\circ\varphi$ is an epimorphism, then $\varphi$ is an epimorphism.
	\end{enumerate}
\end{prop}
\begin{proof}

We first prove that $\calA\subset \Pro(\calA)$ is closed under taking sub-objects. \par

Let $M\in \Pro(\calA)$. Suppose that there exists $\til M\in \calA$ and a monomorphism $\iota:M\rightarrowtail \til M$. We can consider the full subcategory $\calA^M_1\subset \calA^M$ of pairs $(M', a)$ with $a$ being monomorphism. The subcategory $\calA^M_1$ is cofinal. Indeed, if $(M', a)\in \calA^M$, then $\left(M'\times\til M, (a, \iota)\right)\in \calA^M_1$. Let $\calA^M_2\subset \calA^M_1$ be the full subcategory of objects which are minimal, in the sense that if there is $(M'', b)\in \calA^M_1$ with a monomorphism $\varphi\in \Hom_{\calA}(M'', M')$ such that $\varphi \circ b = a$, then $\varphi$ is an isomorphism. By the minimality of the objects of $\calA^M_2$, it is easy to see that the Hom-space $\Hom_{\calA^M_2}\left( (M', a), (M'', b) \right)$ consists of exactly one element for every $(M', a), (M'', b)\in \calA^M_2$. It follows that any object $(M', a)\in \calA^M_2$ yields an isomorphism $a: M\cong M'$. As $\calA$ is artinian, $\calA^M_2$ cannot be empty, whence $M\in \calA$. \par
To prove (ii), in view of~\eqref{equa:limAM}, it suffices to show that $\calA^M_{\epi}$ is cofinal. The previous paragraph shows that for $(M', a)\in \calA^M$, the image $\im(a)$ is in $\calA$. Consider the factorisation $M\xrightarrow{\pi_a}\im(a)\xrightarrow{\bar a} M'$. Then $(\im(a), \pi_a)\in \calA^M_{\epi}$ and there is a morphism $\bar a: (\im(a), \pi_a)\to \left( M', a \right)$ in $\calA^M$. Thus $\calA^M_{\epi}$ is cofinal in $\calA^M$. \par
We prove (iii). Let $M\in \calA$. Since $\calA\subset \Pro(\calA)$ is closed under taking sub-objects, every descending chain of sub-objects of $M$ is in the subcategory $\calA$, which by assumption must stabilise. Thus $M$ is artinian in $\Pro(\calA)$. Suppose that $M\in \Pro(\calA)$ is artinian. There must be a minimal sub-object $M'\subset M$ such that $M / M'$ lies in $\calA$, meaning that the category $\calA^M_{\epi}$ has an initial object. By (ii), $M$ being the projective limit over $\calA^M_{\epi}$ must lie in $\calA$, whence (iii). The assertion (i) follows immediately from (iii). \par
We prove (iv). Let $c:M\to \coker\varphi = C$ be the cokernel. Suppose that $C\neq 0$. Since $C\in \Pro(\calA)$, there exists an epimorphism $p:C\to C'$ with $0\neq C'\in \calA$. Since $p\circ c: M\to C'$ is epimorphism, the composite $ p\circ c\circ\varphi$ is also an epimorphism by hypothesis. However, as $c\circ\varphi = 0$, we see that $C' = 0$, contradiction. Thus $C = 0$ and $\varphi$ is an epimorphism.
\end{proof}

\subsection{}\label{subsec:profunc}
Let $\calA$ and $\calB$ be abelian categories and $F:\calA\to \calB$ an additive functor. We define the extension of $F$:
\begin{align*}
	F:\Pro(\calA)\to \Pro(\calB),\quad F(M) = \prolim{(M',a)\in \calA^M}F(M').
\end{align*}
According to~\cite[8.6.8]{KS06}, if $F$ is exact, then the extended functor $F:\Pro(\calA)\to \Pro(\calB)$ is also exact.
\subsection{}
Suppose that $\calA$ is noetherian-artinian. We define an endo-functor 
\begin{align*}
	\hd:\Pro(\calA)\to \Pro(\calA), \quad
	\hd(M) = \prolim{(M',q)\in\calA^M_{\epi}} \hd(M')
\end{align*}
where $\hd(M')$ is the largest semisimple quotient of $M'$ in $\calA$.
For every $M\in \Pro(\calA)$, there is a canonical map $\pi_M:M\to \hd(M)$.
\begin{prop}[Nakayama's lemma]\label{lemm:hd}
	Let $\calA$ be a noetherian-artinian abelian category. Let $\varphi:N\to M$ be a morphism in $\Pro(\calA)$. Suppose that the composite $N\xrightarrow{\varphi} M\xrightarrow{\pi_M} \hd(M)$ is an epimorphism. Then $\varphi$ is an epimorphism. 
\end{prop}
\begin{proof}
	We first prove the statement in the case where $M\in \calA$. In this case, since $\coker \varphi$ is a quotient of $M$, we have an epimorphism $\hd(M)\twoheadrightarrow \hd(\coker\varphi)$. As the composite $N\to \hd(M)\to \hd(\coker\varphi)$ is zero and is an epimorphism, it implies that $\hd(\coker\varphi) = 0$. As $\calA$ is noetherian, it follows that $\coker\varphi = 0$, so $\varphi$ is surjective. \par
	In general, let $M\in\Pro(\calA)$. Let $(M',q)$ be any object of $\calA^M_{\epi}$. Then $\pi_{M'}\circ q\circ\varphi$ is an epimorphism. By the previous paragraph, $q\circ \varphi$ is also an epimorphism. Then~\autoref{lemm:pro}~(iv) implies that $\varphi$ is an epimorphism. 
\end{proof}

\subsection{}
\begin{prop}\label{lemm:envproj}
	Suppose $\calA$ is an essentially small noetherian-artinian abelian category. Let $M\in \calA$ be a simple object. Then there exists a projective cover $P_M\in \Pro(\calA)$. 
\end{prop}
\begin{proof}
	We construct an object $P^{(n)}\in \Pro(\calA)$ for $n\in \bfN$ by induction. Let $P^{(0)} = M$. For $n > 0$, let
	\begin{align*}
			0\to \proprod{\substack{L\in \Irr(\calA) /\sim \\ \gamma\in \Ext^1_{\calA}\left( P^{(n-1)}, L \right)}}L\to P^{(n)}\to P^{(n-1)}\to 0
	\end{align*}
	be the short exact sequence corresponding to the tautological class
	\begin{align*}
			\Delta = (\gamma)_{L,\gamma}\in \prod_{\substack{L\in \Irr(\calA) /\sim \\ \gamma\in \Ext^1_{\calA}\left( P^{(n-1)}, L \right)}}\Ext^1_{\calA}\left( P^{(n-1)}, L \right).
	\end{align*}
	Put $P = \prolim{n\to \infty} P^{(n)}$. Then $P$ is a projective since we have 
	\begin{align*}
		\Ext^1_{\Pro(\calA)}(P, L) = 0
	\end{align*}
	by construction and since $\calA$ is noetherian-artinian. Let $p:P\to M$ be the obvious epimorphism. \par
	Now, let $\calA^P_M$ be the category whose objects are triples $(\pi, Q, \pi')$, where 
	\begin{itemize}
		\item
			$Q\in \calA$
		\item
			$\pi\in \Hom_{\Pro(\calA)}(P, Q)$ is an epimorphism and 
		\item
			$\pi'\in \Hom_{\calA}(Q, M)$
	\end{itemize}
	such that 
	\begin{itemize}
		\item
			$\pi'\circ\pi = p\in \Hom_{\Pro(\calA)}(P, M)$ and 
		\item
			$\pi'$ induces an isomorphism $\hd(Q)\cong M$. 
	\end{itemize}
	The morphisms are defined by
	\begin{align*}
		\Hom_{\calA^P_M}\left( (\pi_1, Q_1, \pi'_1), (\pi_2, Q_2, \pi'_2) \right) = \left\{ \varphi\in \Hom_{\calA}\left( Q_1, Q_2 \right)\;;\; \varphi\circ\pi_1 = \pi_2 \right\}.
	\end{align*}
	Put
	\begin{align*}
		P_M = \prolim{(\pi, Q,\pi')\in  \calA^P_M} Q.
	\end{align*}
	Then the obvious morphism $P_M\to M$ is a projective cover.
\end{proof}

\section*{Index of notation}
\thispagestyle{plain}
\printindex[ch1]
\printindex[ch2]
\thispagestyle{plain}

\printbibliography
\end{document}